\documentclass{gtart}
\usepackage{amsmath,amsthm,amsfonts,amssymb,amscd,enumerate}
   
\theoremstyle{plain}
\newtheorem{theorem}{Theorem}[section]
\newtheorem{proposition}[theorem]{Proposition}
\newtheorem{conjecture}[theorem]{Conjecture}
\newtheorem{corollary}[theorem]{Corollary}
\newtheorem{lemma}[theorem]{Lemma}
\newtheorem*{decomp}{The Decomposition Theorem}
\newtheorem*{theoremA}{The Main Theorem}
\newtheorem*{Theorem}{Theorem}
\newtheorem*{solomonthm}{Solomon's Theorem}
\newtheorem*{Conjecture}{Conjecture}
\newtheorem*{Question}{Question}
		\mathchardef\MYCOL=\mathcode`:   
    \mathcode`:="8000            
     {\catcode`:=\active
      \gdef:{\MYCOL\thinspace}            
     } 
\newcommand{\rkn}[1]{{\bf #1}\stdspace}
\newtheorem*{I(n)}{$\I(n)$}
\newtheorem*{III(2k+1)}{$\III(2k+1)$}
\newtheorem*{III'(2k+1)}{$\III'(2k+1)$}
\newtheorem*{V(n)}{$\V(n)$}

\theoremstyle{definition}
\newtheorem{step}{Claim}
\newtheorem{definition}[theorem]{Definition}
\newtheorem{example}[theorem]{Example}

\newtheorem{remark}[theorem]{Remark}

\newtheorem*{notation}{Notation}
\newtheorem*{notes}{Note}
\newtheorem*{claim}{Claim}
\newtheorem*{Remark}{Remark}
\newcommand{\ca}{\mathcal {A}}
\newcommand{\cb}{\mathcal {B}}
\newcommand{\cac}{\mathcal {C}}

\newcommand{\ch}{\mathcal {H}}
\newcommand{\ci}{\mathcal {I}}

\newcommand{\cn}{\mathcal {N}_{\mathbf q}}
\newcommand{\cN}{\mathcal {N}}
\newcommand{\cp}{\mathcal {P}}
\newcommand{\car}{\mathcal {R}}
\newcommand{\cs}{\mathcal {S}}

\newcommand{\cu}{\mathcal {U}}

\newcommand{\nn}{{\mathbb N}}

\newcommand{\zz}{{\mathbf Z}}

\newcommand{\R}{\mathbf{R}_{\mathbf q}[W]}
\newcommand{\ltwo}{L^2_{\mathbf {q}}}
\newcommand{\Ltwo}{L^2_q}
\newcommand{\muq}{\mu_{\mathbf {q}}}
\newcommand{\nuq}{\nu_{\mathbf {q}}}
\newcommand{\angt}{\langle T\rangle}

\newcommand{\bH}{{\mathbf H}}
\newcommand{\bR}{{\mathbf R}}
\newcommand{\bZ}{{\mathbf Z}}
\newcommand{\bz}{{\mathbf z}}

\newcommand{\bC}{{\mathbf C}}
\newcommand{\ba}{{\mathbf a}}
\newcommand{\be}{{\mathbf e}}

\newcommand{\bq}{{\mathbf q}}
\newcommand{\bt}{{\mathbf t}}
\newcommand{\bone}{{\mathbf 1}}
\newcommand{\bdelta}{{\mathbf d}}
\newcommand{\ta}{{\tilde{a}}}  
\newcommand{\tih}{{\tilde{h}}}
\newcommand{\te}{{\tilde{e}}}

\newcommand{\wh}{\widehat}
\newcommand{\ol}{\overline}

\def\clap#1{\hbox to 0pt{\hss#1\hss}}

\def\mathclap{\mathpalette\mathclapinternal}

\def\mathclapinternal#1#2{%
\clap{$\mathsurround=0pt#1{#2}$}}

\newcommand{\norm}[1]{\lVert#1\rVert}
\newcommand{\normg}[1]{\lVert#1\rVert_g}
\newcommand{\normf}[1]{\lVert#1\rVert_f}
\newcommand{\comment}[1]{}
\newcommand{\osum}{\oplus}
\newcommand{\ga}{\alpha}
\newcommand{\gb}{\beta}
\newcommand{\gd}{\delta}
\newcommand{\geps}{\varepsilon}
\newcommand{\gf}{\varphi}
\newcommand{\gi}{\iota}

\newcommand{\gs}{\sigma}

\newcommand{\gl}{\lambda}
\newcommand{\go}{\omega}

\newcommand{\gG}{\Gamma}

\newcommand{\gO}{\Omega}
\newcommand{\gL}{\Lambda}
\newcommand{\gS}{\Sigma}

\newcommand{\In}{\operatorname{In}}

\newcommand{\Card}{\operatorname{Card}}
\newcommand{\Hom}{\operatorname{Hom}}

\newcommand{\Ker}{\operatorname{Ker}}
\newcommand{\Ima}{\operatorname{Im}}
\newcommand{\Lk}{\operatorname{Lk}}
\newcommand{\can}{\operatorname{can}}
\newcommand{\tr}{\operatorname{tr}}

\newcommand\mapright[1]{\smash{\mathop{\longrightarrow}\limits^{#1}}}

\newenvironment{enumeratei}{\begin{enumerate}[\upshape (i)]}%
        {\end{enumerate}}

\newenvironment{enumeratea}{\begin{enumerate}[\upshape 
(a)]}{\end{enumerate}}
\newenvironment{enumeratea'}{\begin{enumerate}[\upshape 
(a)$'$]}{\end{enumerate}}

\def\H{\operatorname{\mathfrak h}^\bq}
\def\Hdown{\operatorname{\mathfrak h}_\bq}
\def\Hq{\operatorname{\mathfrak h}^q}
\def\Hqminus{\operatorname{\mathfrak h}_{\bq^{-1}}}
\def\Hqq{\operatorname{\mathfrak h}^{\bq^{-1}}}

\DeclareMathOperator{\I}{\bf I}
\DeclareMathOperator{\III}{\bf III}
\DeclareMathOperator{\V}{\bf V}

\numberwithin{equation}{section}
\begin{document}

\title{Weighted $L^2$-cohomology of Coxeter groups}


\primaryclass{20F55}
\secondaryclass{20C08, 20E42, 20F65, 20J06, 
46L10, 51E24, 57M07, 58J22.}

\shortauthors{Davis, Dymara, Januszkiewicz, Okun}
\author[M.W. Davis]{Michael W. Davis}

\author[J. Dymara]{Jan Dymara}

\author[T. Januszkiewicz]{Tadeusz Januszkiewicz}
\author[B. Okun]{Boris Okun}
\def\nnnn{\noexpand \\}
\addresses{\rm Email:\stdspace \tt mdavis@math.ohio-state.edu \ \ dymara@math.uni.wroc.pl \\
 tjan@math.ohio-state.edu\ \ okun@uwm.edu}

\date{}

\keywords{Coxeter group, Hecke algebra, von Neumann algebra, building, 
$L^2$-cohomology}

\begin{abstract}
Given a Coxeter system $(W,S)$ and a positive real multiparameter $\bq$, we study the ``weighted $L^2$-cohomology groups,'' of a certain 
simplicial complex $\Sigma$ associated to $(W,S)$.  These cohomology 
groups are Hilbert spaces, as well as modules over the Hecke algebra 
associated to $(W,S)$ and the multiparameter $\bq$.  They have a ``von 
Neumann dimension'' with respect to the associated ``Hecke - von Neumann 
algebra,'' $\cn$.  The dimension of the $i^{\mathrm{th}}$ cohomology group 
is denoted $b^i_\bq(\Sigma)$.  It is a nonnegative real number which 
varies continuously with $\bq$.  When $\bq$ is integral, the 
$b^i_\bq(\Sigma)$ are the usual $L^2$-Betti numbers of buildings of type 
$(W,S)$ and thickness $\bq$.  For a certain range of $\bq$, we calculate 
these cohomology groups as modules over $\cn$ and obtain explicit formulas 
for the $b^i_\bq(\Sigma)$.  The range of $\bq$ for which our calculations 
are 
valid depends on the region of convergence of the growth series of $W$.  
Within this range, we also prove a Decomposition Theorem for $\cn$, 
analogous to a theorem of L. Solomon on the decomposition of the 
group algebra of a finite Coxeter group.
\end{abstract}

\maketitle

\section{Introduction}\label{s:intro}
Suppose $(W,S)$ is a Coxeter system.  (The precise definition will 
be given in Section~\ref{s:coxeter}.   For now,  it suffices to say that 
$W$ is a group and $S$ is a set of involutions 
which generate $W$.)
Associated to $(W,S)$ there is a certain
contractible simplicial complex $\gS$ on which $W$ acts properly and 
cocompactly.  (The definition of $\gS$ can be found in 
\cite{d83,d98,d02,davisbook, dm,moussong}, as well as in Section~\ref{s:spaces} 
below.)  
Let  $i:S\to I$ be a function to some index set $I$ so that $i(s)=i(s')$ 
whenever $s$ and $s'$ are conjugate.    
Given an $I$-tuple 
$\bq=(q_i)_{i\in I}$ of positive real numbers,  the 
second author  \cite{dym} defined certain ``weighted $L^2$-cohomology 
spaces, '' here denoted $\ltwo \ch^i(\gS)$.   The weighted 
$L^2$-cochain complex, $\ltwo C^*(\gS)$, 
is a subcomplex  of the complex $C^*(\gS;\bR)$ of ordinary 
cellular cochains.  It consists of those cochains which are square summable with 
respect to an inner product defined via a weight function depending on 
the multiparameter $\bq$.  
As we explain in 
Sections~\ref{s:h-neumann} and \ref{s:weighted},  to each of the 
Hilbert spaces $\ltwo \ch^i(\gS)$ one can attach 
a ``von Neumann dimension.''   It 
is a nonnegative real number,  denoted $b^i_\bq(\gS)$ and called the 
\emph{$i^{th}$ $\ltwo$-Betti number of $\gS$}. 

Our principal interest in the weighted $L^2$-cohomology of $\gS$ lies in 
the fact that
it computes the $L^2$-cohomology of buildings of type $(W,S)$.  
Here $\bq$ is an $I$-tuple of positive 
integers called the ``thickness vector'' of the building.  
(So, for  buildings, only $\bq$ with integral components can 
occur.)   

The theory of the weighted $L^2$-cohomology of $\gS$ is closely tied to 
several other topics, for example, growth series of Coxeter groups,  
decompositions of ``Hecke - von Neumann algebras'' and the Singer 
Conjecture.    Moreover, as $\bq$ 
goes from $0$ to $\infty$, $\ltwo \ch^*(\gS)$ interpolates between 
ordinary cohomology and cohomology with compact supports.  For these 
reasons, we believe that the study of weighted $L^2$-cohomology of Coxeter 
groups has intrinsic interest, 
independent of its connection to buildings.

Let $\bt:=(t_i)_{i\in I}$ be an $I$-tuple
of indeterminates.  Write $t_s$ instead of $t_{i(s)}$.  If $s_1\cdots 
s_k$ is a reduced expression for an element $w\in W$, then the monomial
$t_w:=t_{s_1}\cdots t_{s_k}$ is independent of the choice of 
reduced expression for $w$.  
The \emph{growth series} 
for $W$ is the power series in $\bt$ defined by
	\[
	W(\bt):=\sum_{w\in W} t_w.
	\]
It is a  rational function of $\bt$ (\cite{bourbaki,serre}).   We 
give several explicit formulas for it in Lemma~\ref{l:growth} of 
Section~\ref{s:growth}.  (In the case where $I$ is a singleton, so that
$\bt$ is a single indeterminate $t$, we 
have $t_w=t^{l(w)}$, where $l(w)$ denotes the word length of $w$.  
So, in the case of a single indeterminate, 
$W(t)=\sum t^{l(w)}$ is the usual growth series.)  

Let $\bR^{(W)}$ denote the vector space of finitely supported, 
real-valued functions on $W$ and let $(e_w)_{w\in W}$ be its 
standard basis.  

As we  explain in Section~\ref{s:hecke}, associated 
to each multiparameter $\bq$, there is a deformation of 
the group algebra of $W$ called the ``Hecke algebra'' 
(or sometimes the ``Iwahori-Hecke algebra'')  
of $W$.  We denote it by $\R$.     
When $\bq=\bone$ (the $I$-tuple with all components equal to $1$),  
$\R$ is the group algebra of $W$.  
(No matter what $\bq$ is, the underlying vector space of $\R$ is always 
$\bR^{(W)}$.)  

Also associated to $\bq$, there is an inner product $\langle\ ,\ 
\rangle_\bq$ on $\bR^{(W)}$ defined by $\langle e_w,e_{w'}\rangle_\bq 
=q_w\gd_{ww'}$, where $\gd_{ww'}$ is the Kronecker delta.  The 
completion of $\bR^{(W)}$ with respect to this inner product 
is denoted  $\ltwo (W)$ or simply 
$\ltwo$ when $W$ is understood.  $\ltwo$ is an $\R$-bimodule.
There is an anti-involution on $\R$, 
denoted by $x\to x^*$ and defined by $(e_w)^*:=e_{w^{-1}}$.  
Moreover, $\langle yx,z\rangle_\bq = \langle 
y,zx^*\rangle_\bq$, i.e., right translation by $x^*$ is the adjoint 
of right translation by $x$.
As is explained in 
\cite{dym} and Proposition~\ref{p:hilbertalgebra} below, this makes 
$\R$ into a ``Hilbert algebra'' in the sense of 
Dixmier \cite{dix81}.  It follows that there is an associated von Neumann 
algebra  
$\cn$ acting on $\ltwo$ from the right.  It can be defined as the algebra 
of bounded 
linear operators on $\ltwo$ which commute with the 
left $\R$-action.  $\cn$ is the \emph{Hecke - von Neumann algebra} 
associated to $\bq$.  ($\cn$ is a completion of $\R$ acting from the right 
on $\ltwo$.)  
As in the case of a von Neumann algebra associated to 
a group algebra, $\cn$ is 
equipped with a trace which one can use to define the 
``dimension'' of any 
$\R$-stable closed subspace $V$ of a finite direct sum of copies of 
$\ltwo$.  

Suppose $W$ acts  as a reflection group on a CW complex 
$\cu$ with a strict fundamental domain $Z$.  Assume 
further that for each $s\in S$ 
there is a subcomplex $Z_s\subseteq Z$, called a ``mirror'' of $Z$, so 
that $s$ acts on $\cu$ as a reflection across $Z_s$.  Then $\cu$ is formed 
by gluing together copies of $Z$, one for each element of $W$.  In other 
words, $\cu\cong (W\times Z)/\sim$, where the equivalence relation $\sim$ 
is defined in an obvious fashion.  (See Section~\ref{s:spaces}.)  
The complex $\gS$ can be described in this manner:  the 
fundamental chamber for $W$ on $\gS$ is denoted by $K$ instead of $Z$.  

$C^i_c(\cu)$ is the space of finitely supported, real-valued, 
cellular $i$-cochains on 
$\cu$.  For each oriented $i$-cell $\gs$ of $\cu$, let $e_\gs$ be its 
characteristic function.  So, $\{e_\gs\}_{\gs\in\{\text{$i$-cells}\}}$ is 
a basis for $C^i_c(\cu)$.   As in \cite{dym}, there is a definition of
an inner product on $C^i_c(\cu)$ 
similar to the definition  of $\langle \ ,\ \rangle_\bq$ on $\bR^{(W)}$.  
The $e_\gs$ form an orthogonal basis; however, the norm of $e_\gs$ need 
not be $1$.  Instead, one uses $\bq$ 
to weight the inner product so that 
$\langle e_\gs,e_\gs\rangle_\bq=q_w$, where $w$ is the shortest element of 
$W$ such that  $\gs\subseteq wZ$.  Let $\ltwo C^i(\cu)$  denote
the completion of
$C^i_c(\cu)$ with respect to this inner product.  

As  explained in \cite{dym},  
as well as in Section~\ref{s:weighted}, $\ltwo C^i(\cu)$  can be 
identified 
with a $\R$-stable subspace of $\oplus \ltwo$.  The 
coboundary maps are $\R$-equivariant.  So, the (reduced) cohomology 
group $\ltwo\ch^i(\cu)$ is a closed $\R$-stable subspace of $\oplus \ltwo$ 
and therefore, has a well-defined von Neumann dimension, $b^i_\bq(\cu)$.  
The alternating 
sum of the  $b^i_\bq(\cu)$ is denoted $\chi_\bq (\cu)$ and called the 
\emph{$\ltwo$-Euler characteristic} of $\cu$. It is proved in \cite{dym} 
(and in Proposition~\ref{p:euler}) that  $\chi_\bq 
(\gS)=1/W(\bq)$.  (Recall $W(\bt)$ is a rational function.) Moreover,  
the Betti numbers $b^i_\bq(\cu)$ are continuous functions of $\bq$ 
(Theorem~\ref{t:cont}).

Let $\car$ denote the region of convergence of $W(\bt)$ and let
$$\car^{-1}:=\{\bq\mid \bq^{-1}\in \car\},$$ 
where $\bq^{-1}:=(q_i^{-1})_{i\in I}$.  The closures of these regions are 
denoted $\ol{\car}$ and $\ol{\car^{-1}}$, respectively.
(When $I$ 
is a singleton, we write $q$ instead of $\bq$ and $t$ instead of 
$\bt$.  In this case, $W(t)$ is a power series in one variable.  As such, it 
has a radius of convergence $\rho$ and $\car=(0,\rho)$.) 

The main result of this paper, Theorem~\ref{t:decoupled}, 
is a calculation of $\ltwo \ch^i(\cu)$ (as a $\cn$-module) 
for $\bq\in \ol{\car}\cup \ol{\car^{-1}}$.  
It also gives a formula for the $b^i_\bq(\cu)$ in this range of $\bq$. 
Roughly speaking, the answer is that for $\bq\in 
\ol{\car}$,  $\ltwo \ch^*(\cu)$ looks like ordinary cohomology while for 
$\bq\in \ol{\car^{-1}}$, it looks like
cohomology with compact supports.  Before stating the result precisely, we need 
to set up some notation and recall some background.

Given $T\subseteq S$, the subgroup $W_T$ 
generated by $T$ is called a \emph{special subgroup}.  
It is also a Coxeter group.  The subset $T$ is \emph{spherical} 
if $W_T$ is finite.  Let $\cs$ denote the poset of spherical subsets of 
$S$.  Given an element $w\in W$, set 
$\In (w):=\{s\in S\mid l(ws)<l(w)\}$,  
i.e., $\In (w)$ is the set of letters in $S$ with which a reduced 
expression for $w$ can end.  It turns out that for any $w\in W$, $\In (w)$ 
is always a spherical subset of $S$.  
For each $T\in \cs$, let $W^T:=\{w\in W\mid \In (w)=T\}$ and 
let $\bZ (W^T)$ denote the free abelian group of  finitely supported 
functions on $W^T$.  For any $U\subseteq S$, $Z^U$ denotes the union of 
those mirrors $Z_s$,  with $s\in U$.

\begin{enumeratea}
\item
The homology of $\cu$ is computed in \cite{d87}.  The answer is
\[
H_*(\cu)\cong \bigoplus_{T\in \cs} H_*(Z,Z^T)\otimes \bZ (W^T).
\]
(This implies, in particular, that $\gS$ is acyclic.)
The answer for cohomology is similar, except that it is necessary to 
replace $\bZ (W^T)$ by the abelian group of all functions $W^T\to \bZ$.
\item
The cohomology with compact supports of $\cu$ can be computed as in 
\cite{d98}.  The answer is
\[
H^*_c(\cu)\cong \bigoplus_{T\in \cs} H^*(Z,Z^{S-T})\otimes \bZ (W^T).
\]
\end{enumeratea}
Another proof of the formulas in (a) and (b) is given in \cite{ddjo06, davisbook}.

Given $U\subseteq S$, in Section~\ref{s:h-neumann},  we define idempotents 
$a_U$ and $h_U$ in $\cn$ by
    \begin{align*}
    a_U&:=\frac{1}{W_U(\bq)}\sum_{w\in W_U}e_w,\\
    h_U&:=\frac{1}{W_U(\bq^{-1})}\sum_{w\in W_U}\geps_w q_w^{-1}e_w,
    \end{align*}
where $\geps_w:=(-1)^{l(w)}$.
These idempotents are defined provided $\bq\in \car_U$ 
in the case of $a_U$ and provided 
$\bq\in\car^{-1}_U$ in the case of $h_U$. ($\car_U$ denotes the 
region of 
convergence for $W_U(\bt)$.)   Let $A_U\subseteq \ltwo$ stand for 
$\Ima a_U$ if $\bq\in \car_U$ 
and for the $0$-space, otherwise.  $A_U$ is a  closed $\R$-stable subspace of 
$\ltwo$.  Another closed $\R$-stable subspace is defined by
\[
D_U:=A_{S-U}\cap\Big(\sum _{V\subset U}A_{S-V}\Big)^\perp.
\]
(Throughout this paper we will denote inclusion of a subset by $\subseteq$ 
and use  $\subset$ for inclusion of  a proper subset.)

Here is the precise statement of our calculation of $\ltwo$-cohomology.  
(Compare it with statements (a) and (b) above.)

\begin{theoremA}[Theorem~\ref{t:decoupled} in Section~\ref{s:decouple}]\ 
\begin{enumeratea}
\item
If $\bq\in \ol{\car}$, then
\[
\ltwo \ch^{\ast}(\cu)\cong\bigoplus _{T\in \cs}H^{\ast}(Z,Z^T)\otimes D_T.
\]
\item
If $\bq\in \ol{\car^{-1}}$, then
\[
\ltwo \ch^{\ast}(\cu)\cong\bigoplus _{T\in \cs}H^{\ast}(Z,Z^{S-T})\otimes
D_{S-T}.
\]
\end{enumeratea}
\end{theoremA}

\noindent
(To compare this with the previous answers for ordinary cohomology 
and cohomology with compact supports,  we note that, by Theorem~\ref{t:dgl}, for 
$\bq\in \car$,  $\{e_wh_Ta_{S-T}\}_{w\in W^T}$ spans a  dense 
subspace of $D_T$; while for $\bq\in \car^{-1}$, $\{e_wh_{S-T}a_T)\}_{w\in 
W^T}$ spans a dense subspace of $D_{S-T}$.)

The proof of the Main Theorem depends on the following result.

\begin{decomp}[Theorem~\ref{t:decomp} in Section~\ref{s:decomp}]\ \nopagebreak 
\begin{enumeratea}
\item
If $\bq\in \ol{\car}$, then
\[
\sum_{T\in \cs} D_T
\]
is a direct sum decomposition and  a dense subspace of $\ltwo$.
\item
If $\bq\in \ol{\car^{-1}}$, then 
\[
\sum_{T\in \cs} D_{S-T}
\]
is a direct sum decomposition and  a dense subspace of $\ltwo$.
\end{enumeratea}
\end{decomp}

In the case when $W$ is finite and $\bq=\bone$ (i.e.,  when the 
Hecke algebra is the group algebra) a similar result was proved by 
Solomon 
\cite{s} in 1968.  In Section~\ref{s:solomon} we give a version of the 
Decomposition Theorem (namely, Theorem~\ref{t:gensolomon}) which is more 
transparently a generalization of Solomon's Theorem than the version 
stated above.  The Decomposition 
Theorem is also compatible with the theory of representations of Hecke 
algebras developed by Kazhdan--Lusztig \cite{dl} (cf. Remark~\ref{r:kl}).

Although the Main Theorem is a consequence of the Decomposition Theorem, 
our proof of the Decomposition Theorem ultimately is based on a special 
case of the Main Theorem from \cite{dym}.  The result of \cite{dym} 
states 
that, for $\bq\in\car$, the $\ltwo$-homology of $\gS$ vanishes except in 
dimension $0$.  (N.B. To calculate homology, $\ltwo\ch_*(\gS)$,  
from $\ltwo C_*(\gS)$ 
one does not use the usual boundary map but rather, the adjoint of the usual 
coboundary map.) In Section~\ref{s:ruins} we apply this vanishing result 
to show that, for $\bq\in \car$, the relative $\ltwo$-homology of certain 
pairs of 
subcomplexes of $\gS$ vanishes except in the bottom dimension.  (These 
pairs of subcomplexes are dubbed ``ruins'' in Section~\ref{s:spaces}.) 
For $\bq\in\car$, these vanishing results are essentially an equivalent 
version of the Decomposition Theorem.  One then uses
a certain isomorphism $j:\cn\to \cN_{\bq^{-1}}$ to 
convert the statement of Decomposition Theorem for $\bq\in \ol{\car}$ into 
its 
statement  for $\bq\in \ol{\car^{-1}}$.

The key role played by the case $\bq\in \car$ in this sketch of the proof  
is probably the most compelling reason for studying 
weighted $L^2$-cohomology with $\bq$ an $I$-tuple of arbitrary positive 
real numbers.  
When $W$ is infinite, the vector
$\bq\in \car$ never has all its components equal to positive integers.  
So, on the face of it, the case $\bq\in \car$ of the 
Main Theorem would never seem to be applicable to nonspherical buildings.  
However,  because of various dualities 
(such as the $j$-isomorphism) which switch 
$\bq$ with $\bq^{-1}$, the results for $\bq\in \car$ are equivalent to 
results for $\bq\in \car^{-1}$ and these are applicable to buildings.

For $\bq\in \ol{\car}\cup \ol{\car^{-1}}$,
the Main Theorem (in particular, its version as 
Theorem~\ref{t:decoupledsigma}) gives a complete calculation of 
$\ltwo \ch^*(\gS)$.  
On the other hand, our knowledge about 
what happens for $\bq\notin \ol{\car}\cup \ol{\car^{-1}}$ is fragmentary.
For example, suppose $\gS$ is an $n$-manifold.  By the 
Main Theorem, $\ltwo \ch^*(\gS)$ is concentrated in dimension  $0$ 
for $\bq\in \ol{\car}$ and in dimension $n$ for $\bq\in \ol{\car^{-1}}$.  
We note that when $\gS$ is a manifold (without boundary), $W$ is infinite 
and so, $\bone\notin \car\cup\car^{-1}$.
When $\bq=\bone$, a version of the Singer Conjecture asserts that the 
weighted $L^2$-cohomology of $\gS$ vanishes except in dimension 
$\frac{n}{2}$.  (When $n$ 
is odd, this is to be interpreted as meaning that $L^2_\bone\ch^*(\gS)$ 
vanishes in all dimensions.)  In \cite{do} the first and fourth authors explained some evidence for this conjecture.  For a general $\bq$, in the case where 
$\gS$ is a $n$-manifold, there is a version of Poincar\'e duality 
which exchanges $\bq$ with $\bq^{-1}$ (as well as dimension $k$ 
with dimension $n-k$); see \cite{dymara} or 
Proposition~\ref{p:duality} below.  So, when $\gS$ is a manifold, knowledge 
of $\ltwo \ch^*(\gS)$ for $\bq\leq\bone$ also determines it for $\bq>\bone$. 
In Section~\ref{s:sphere} we explain that the right 
generalization of this version of the Singer Conjecture for $\bq=\bone$ is the 
following.  

\begin{Conjecture}[Conjecture~\ref{conj:singer}]
Suppose $\gS$ is an $n$-manifold.  Then 
\[
\ltwo \ch^k(\gS)=0 \text{ for $k> \frac{n}{2}$ and 
$\bq\le\bone$.}
\]
\end{Conjecture}

In Section~\ref{s:rightsphere}, by modifying the arguments of \cite{do}, 
we prove it in the case 
where $W$ is right-angled and $n\le 4$.
In the same section, we  give  examples where $\gS$ is a 
$4$-manifold and where for certain 
$\bq\notin \ol{\car}\cup\ol{\car^{-1}}$, the 
$\ltwo$-cohomology fails to be concentrated in a single dimension  
(it is nonzero in both dimension $2$ and $3$.)

Next, we make  a few remarks concerning 
buildings.  Buildings come in different types, where the 
``type''  of a building is a Coxeter system $(W,S)$.
In the case of a classical building associated  to an algebraic group, 
its type is always a  spherical or Euclidean reflection group.  The 
simplest example of a Euclidean reflection group is when $W$ is the 
infinite dihedral group acting on the real line and  $S$ consists of 
the two reflections about the endpoints of a fundamental interval.  A 
building of this type is a tree.  (See Example~\ref{ex:trees} for more 
details.)  Many other types for  buildings are  
possible.  Most of our interest in this paper lies with these nonclassical 
types.

Roughly speaking, a building of type $(W,S)$ consists of a set $\Phi$ of 
``chambers'' 
and a family, indexed by $S$, of
``adjacency relations'' on $\Phi$.  An example of a building
is $W$ itself -- the adjacency relation corresponding to $s\in
S$ is defined by calling two distinct elements of $W$ 
\emph{$s$-adjacent} if they form to the same coset of $W_{\{s\}}$.

To define the ``geometric realization'' $X$ of a building, one first 
declares the geometric realization of any chamber to be isomorphic to the 
fundamental chamber $K$ of $\gS$.  One then amalgamates copies of $K$, 
one for each element of $\Phi$, by gluing together  chambers corresponding 
to $s$-adjacent elements of $\Phi$ along the mirror corresponding to 
$s$.  Details of this construction can be found in \cite{d98a}, 
as well as in Section~\ref{s:buildings}.  (N.B. When 
$W$ is an irreducible Euclidean reflection group, $K$ is a simplex and $X$
has the structure of a  simplicial complex 
in which  the top-dimensional simplices are the chambers; however, in  the
general case, this is not the correct picture of the geometric realization 
of a building.)

A group $G$ of automorphisms of a building is \emph{chamber transitive} if 
it acts transitively on $\Phi$.  If the building admits a chamber 
transitive automorphism group, then, for any given $\gf_0\in \Phi$, the 
number of chambers which are $s$-adjacent to $\gf_0$ is independent of the 
choice of $\gf_0$.  We denote this number by $q_s$, i.e., 
\[
q_s=\Card \{\gf\in\Phi \mid \gf \text{ is $s$-adjacent to } \gf_0 \text{ 
and } \gf\neq \gf_0\}.
\]
Moreover, if $s$ and $s'$ are conjugate in $W$, then $q_s=q_{s'}$.  We 
assume throughout that the building has finite thickness, i.e., that each 
$q_s$ is finite.  We then get a well-defined $I$-tuple of integers 
$\bq:=(q_i)_{i\in I}$, called the \emph{thickness vector} of the building, 
where $I$ is the set of conjugacy classes in $S$ and where $q_i:=q_s$ for 
any representative $s$ for $i$.  For example, the thickness vector of $W$ 
(considered as a building) is $\bone$.

How do Hecke algebras arise in the theory of buildings?  Let $\Phi$ be a 
building of finite thickness with a chamber transitive automorphism group 
$G$ and with thickness vector $\bq$. Fix a base chamber $\gf_0\in \Phi$.  
Using the ``$W$-distance'' from $\gf_0$, one gets a retraction $r:\Phi \to 
W$.  Let $C_c(\Phi)$ denote the space of finitely supported, real-valued 
functions on $\Phi$.  It is an algebra with  product  given by convolution.
Consider the subspace $J\subseteq C_c(\Phi)$ consisting of those functions 
which are constant on the fibers of $r$.  It is a subalgebra.  As a vector 
space, $J$ can be identified with $\bR^{(W)}$; however, the product is not 
the usual one for the group algebra.  
As the reader has probably guessed, $J$ 
is identified with the Hecke algebra $\R$, where the multiparameter $\bq$ 
is the thickness of $\Phi$.  

Let $X$ denote the geometric realization of the building $\Phi$.  The 
retraction $r:\Phi\to W$ induces a topological retraction $X\to \gS$, 
which 
we continue to denote by $r$. 
This induces an inclusion 
$r^*:C^*_c(\gS)\to C^*_c(X)$ of (finitely supported) cellular cochains.  
The standard inner product on $C^*_c(X)$ restricts to the  inner product 
$\langle\ ,\ \rangle_\bq$ on $C^*_c(\gS)$.  In this way, $\ltwo\ch^*(\gS)$ 
is identified with a subspace of $\ch^*(X)$, the ordinary reduced
$L^2$-cohomology of $X$.

Since $\Phi$ has finite thickness, $G$ is locally compact and hence, has a 
Haar measure $\mu$, which we normalize by the condition, $\mu (B)=1$, 
where $B$ denotes the stabilizer of a chamber.  Given $\mu$, 
we have the Hilbert space $L^2(G,\mu)$ of square integrable functions on 
$G$ 
and  a von Neumann algebra $\cN (G)$.  Since 
$\ch^i(X)$ is an $\cN (G)$-module, it has a ``dimension'' with respect 
to  $\cN (G)$.
This number is called the \emph{$i^{\mathrm{th}}$ $L^2$-Betti number} and 
denoted $b^i(X;G)$.   
It is proved in \cite{dym} (under slightly stronger 
hypotheses), as well as in 
Theorem~\ref{t:bettibuilding} of Section~\ref{s:buildings}, that $b^i(X;G)=b^i_\bq(\gS)$.

In \cite{dj} the second and third authors calculated $\ch^*(X)$ under 
the assumption that the thickness vector $\bq$ is very large.  The result 
of \cite{dj} is similar to statement (b) of our Main Theorem: it says that 
for $\bq >>\bone$,
\[
\ch^*(X)\cong \bigoplus_{T\in \cs} H^*(K,K^{S-T})\otimes \wh{D}_{S-T}.
\]
where $\wh{D}_{S-T}$ is a specific subrepresentation of $L^2(G/B)$ 
analogous to the subspace $D_{S-T}\subset\ltwo$.  (The 
notation in \cite{dj} is different.)
In Theorem~\ref{t:decoupledsigma'} of Section~\ref{s:buildings} we use  
our Main Theorem to show that, in fact, this formula is valid for all 
$\bq\in \ol{\car^{-1}}$.

If $W$ is a Euclidean reflection group, then the radius of convergence of 
$W(t)$ is $1$ (cf. Proposition~\ref{p:subexp}).  
It follows that  $\bq\in \ol{\car^{-1}}$ 
whenever $\bq\geq \bone$.  From this we deduce the following (known) 
result.

\begin{Theorem}[Corollary~\ref{c:eucbldg} in Section~\ref{s:sphere}]
Suppose $X$ is an affine building (i.e.,  of Euclidean type) and that
its automorphism group is chamber transitive.  Then $\ch^*(X)$ is 
concentrated in the top dimension.
\end{Theorem}

We want to emphasize three points: 
\begin{enumeratei}
\item
There are many interesting classes of buildings which are neither spherical nor affine.
\item
The rational function $W(\bq)$ is completely explicit (cf. Lemma~\ref{l:growth}) and is easily calculated in any given case.
\item
In practice, calculation shows few $I$-tuples of positive integers lie outside $\car^{-1}$.
\end{enumeratei}

With regard to (i), when $W$ is right-angled, there is a building $\Phi$ of type $(W,S)$ and thickness $\bq$ for any $I$-tuple of positive integers $\bq$ (cf. Example~\ref{ex:rabldg}).  
When $(W,S)$ satisfies the ``crystallographic condition'' (that all $m_{ij}=2,3,4,6$ or $\infty$), Tits proved the existence of ``Kac-Moody groups'' over finite fields.  These give locally finite buildings with chamber transitive automorphism groups.  (In this Kac-Moody case, all $q_s$ are restricted to be the order of a finite field.)

With regard to (iii), consider the following concrete example.  Suppose $W$ is the group generated by the reflections across the faces of a right-angled dodecahedron in hyperbolic $3$-space.  Let $I$ be a singleton and $\bq=q$, a positive integer.  Suppose $X_q$ is the (dodecahedral) building of type $W$ and thickness $q$.
If $q\ge 7$, our results give that $\ch^*(X_q)$ is concentrated in the top dimension ($=3$). Previously, this was only known for $q>10^{60}$ (cf. \cite{dj}.)   On the other hand, by Theorem~\ref{t:10.4.1}, for $2\le q\le 6$, $\ch^*(X_q)$ is concentrated in dimension $2$ (and for $q=1$ it vanishes identically).

The results of this paper raise more questions than they answer.  Here are 
two such:
\begin{itemize}
\item
Is there a version of this theory for weighted differential forms?
\item
Is there a version of this for groups other than Coxeter groups?
\end{itemize} 
The short answer to both is ``yes.''  In both 
cases a good deal of foundational work remains to be done.

As for the first question, 
there  exists a literature on weighted $L^2$ de Rham 
cohomology on a Riemannian manifold $M$, for example, \cite{bueler}. 
The inner product on the vector space 
of compactly supported, smooth forms on $M$ is 
modified via a weight function of the form $x\to q^{d(x)}$, where $q\in 
(0,\infty)$, $x\in M$ and $d(x)$ is the distance from a basepoint.   As 
one 
would expect, when $M=\bR^n$,  the   
weighted $L^2$-cohomology is concentrated in dimension $0$ if $q<1$ and in 
dimension $n$ if $q>1$.
In Section~\ref{s:rightsphere}, we use a 
version of this weighted de Rham theory on hyperbolic space equipped with 
an isometric action of a group $W$ generated by reflections across the 
faces of a fundamental polytope $K$.  This time the weight 
function is a step function of the form $x\to q^{l(w)}$, where $w\in W$ 
is such that $x\in wK$.  In this case, the de Rham version and the 
cellular 
version of weighted $L^2$-cohomology are canonically isomorphic.

As for the second question, given a discrete group $\gG$, a CW complex $X$ 
equipped with a cellular $\gG$-action and a positive real number $q$, one 
can deform the standard inner product on $C^*_c(X)$ via a weight function 
of the form $\gamma\to q^{l(\gamma)}$ and then define the weighted $L^2$ 
(cellular) 
cohomology groups of $X$.  As before, as $q$ varies from $0$ to $\infty$, 
these groups interpolate between ordinary cohomology and cohomology with 
compact supports.  The missing feature for a general group $\gG$ (as 
opposed to a
Coxeter group) is that we do not have a deformation of the group 
algebra analogous to the Hecke algebra.   
We will say more about this question in 
Section~\ref{s:final}.  We believe that this topic also has an intrinsic 
interest and we hope to write more about it in the future.

The first author was partially supported by NSF grant DMS 0104026.
The second author was partially supported by KBN grant 2 PO3A 017 25 and 
a scholarship from the Foundation for Polish Science. 
The third author was partially supported by KBN grant 2 P03A 017 25.

\section{Coxeter systems}\label{s:coxeter}
A \emph{Coxeter matrix} on a set $S$ is an $S\times S$ symmetric
matrix $M=(m_{st})$ with entries in $\nn \cup \{\infty\}$ such that each
diagonal entry is 1 and each off-diagonal entry is $\geq 2$.  The matrix
$M$ gives a presentation for an associated \emph{Coxeter group} $W$:  the
set of generators is $S$ and there is a relation
    \begin{equation*}
    (st)^{m_{st}}=1,
    \end{equation*}
for each pair $(s,t)$ of elements in $S$ with
$m_{st}\neq \infty$.  The purpose of this section is to recall some standard 
facts about such groups.  Proofs of most of these facts can be found in 
\cite{bourbaki}.  

The natural map $S\to W$
is injective and henceforth, we identify $S$ with its image in $W$.
Moreover, each element of $S$ has order 2 in $W$ and the order of $st$ in
$W$ is $m_{st}$.  The pair $(W,S)$ is a \emph{Coxeter system}.

Given an element $w\in W$, $l(w)$ denotes its word length.  An expression
for $w$ as a word in $S$, $w=s_1\cdots s_l$,  is a \emph{reduced 
expression}
if $l=l(w)$.

Given $T\subseteq S$,  $W_T$ denotes the subgroup generated by
$T$.  Such a $W_T$ is a \emph{special subgroup} of $W$.  The pair 
$(W_T,T)$ is the Coxeter system whose Coxeter matrix is given by
the restriction of $M$ to $T$ (\cite[Theorem 2 (i), p.
12]{bourbaki}).  The subset $T$ is
\emph{spherical} if $W_T$ is finite.

For $T\subseteq S$ and $w\in W$, the
coset $wW_T$ contains a unique element of minimum length.  This element
is said to be \emph{$(\emptyset,T)$-reduced}.  Moreover, an element $w$ is
$(\emptyset,T)$-reduced if and only if $l(wt)>l(w)$ for all $t\in T$.
(See  \cite[Ex. 3, pp. 31--32]{bourbaki}.)  Let $X_T$ denote the set of
$(\emptyset,T)$-reduced elements of $W$. 

If $W_T$ is finite, then it contains a unique element $w_T$ of maximum
length, called the \emph{element of longest length}.  It is characterized
by the property that $l(w_Tt)<l(w_T)$ for all $t\in T$
(\cite[Ex. 22, p. 40]{bourbaki}).

Given $w\in W$, set
$\In w:=\{s\in S\mid l(ws)<l(w)\}$.
It follows from the ``Exchange Condition'' (cf. \cite[p. 7]{bourbaki}) that
$s\in \In w$ if and only if $w$ has a reduced expression with final letter $s$.
Thus, $\In w$ is the set of letters with which a reduced expression for $w$
can end.  A key fact (\cite[Lemma 7.12]{d83})
is that $\In w$ is always a spherical subset of $S$.

For any spherical subset $T\subseteq S$, define
    \begin{equation}\label{e:defWT}
    W^T:=\{w\in W\mid \In w=T\}.
    \end{equation}

\rk{The simplicial complex $\gS$} 
Given a poset $\cp$ and an element $p\in \cp$, define
$\cp_{>p}:=\{x\in \cp \mid x>p\}$.  Subposets $\cp_{<p}$,
$\cp_{\geq p}$ and $\cp_{\leq p}$ are defined similarly.  Associated to any
poset $\cp$ there is a simplicial complex $|\cp|$, called its 
\emph{geometric
realization}; its vertex set is $\cp$ and a nonempty finite subset of $\cp$
spans a simplex if and only if it is totally ordered.

Let $\cs$ denote the set of spherical subsets of $S$, partially ordered by
inclusion and let
    \begin{equation}\label{e:cs^i}
    \cs^{(i)}=\{T\in\cs\mid \Card (T)=i\}.
    \end{equation}
$\cs$ has a minimum element, namely, $\emptyset$.
$\cs_{>\emptyset}$ is the poset of simplices of a simplicial complex
denoted  by $L(W,S)$ (or $L$ for short) and called the \emph{nerve} of 
$(W,S)$.
(In other words, the vertex set of $L$ is $S$
and a nonempty subset $T\subseteq S$ spans a simplex if and only if it is
spherical.)  $\cs^{(i)}$ is the set of $(i-1)$-simplices in $L$.
\comment{
Define $ \geps (T):=(-1)^{\Card (T)}$. Conversely, given a
simplicial complex $L$ we let $S=S(L)$ to be its vertex set, and define the
\emph{face poset} $\cs=\cs (L)$ of $L$ to be the poset of simplices of $L$,
including the empty simplex.
}

We are also interested in $W\cs$, the \emph{poset  of spherical cosets}.
It is defined as the disjoint union of the sets $W/W_T$,  $T\in\cs$.
Thus, a typical element of $W\cs$ is a coset $wW_T$ for some $T\in \cs$.
The partial order is inclusion.

The geometric realization of $\cs$ is denoted $K$ and the geometric
realization of $W\cs$ by $\gS$.  The group $W$ acts properly on the
simplicial complex $\gS$; the orbit space is the finite complex $K$.
The most important property of $\gS$ is that it is contractible
(\cite[Theorem 10.3 and Section 14]{d83}).

\section{Growth series}\label{s:growth}
Suppose we are given a Coxeter system $(W,S)$, an index set $I$ and a
function $i:S\to I$ so that $i(s)=i(s')$ whenever $s$ and $s'$ are
conjugate in $W$.  (The largest possible choice for $\Ima i$ is the set of 
conjugacy
classes of elements in $S$ and the smallest possible choice 
is a singleton.)
Let $\bt=(t_i)_{i\in I}$ stand for an $I$-tuple of indeterminates.
Write $t_s$ for $t_{i(s)}$.  If $s_1\cdots s_l$
is a reduced expression for $w$, then define $t_w$ to be the monomial
    $t_w:=t_{s_1}\cdots t_{s_l}$.
It  follows from Tits' solution of the word problem for Coxeter groups
(see \cite{tits} or \cite{brown})  that $t_w$ is independent of the choice of
reduced expression for $w$.

For any subset $X$ of $W$, define a power series in $\bt$
    \begin{equation}\label{e:growth2}
    X(\bt):=\sum_{w\in X} t_w.
    \end{equation}
$W(\bt)$ is the \emph{growth series} of $W$ and, for any subset $T$
of $S$, $W_T(\bt)$ is the growth series of the special subgroup
$W_T$.  

\begin{notation}
The region of convergence of $W_T(\bt)$ in $(0,+\infty)^I$ is  denoted
$\car_T$.  Write $\car$ instead of 
$\car_S$.  Put $\car_T^{-1}:=\{\bz\in (0,+\infty) ^I\mid \bz^{-1}\in\car_T\}$.
Denote the closure of the region of convergence by $\ol{\car}$ and put 
$\partial\car:=\ol{\car}-\car$.  Define $\ol{\car^{-1}}$ and $\partial 
\car^{-1}$ similarly.
\end{notation}

From the fact that all the coefficients in $W(\bt)$ are nonnegative real
numbers, we immediately get the following lemmas.

\begin{lemma}\label{l:region}
If $U\subseteq T\subseteq S$, then $\car\subseteq \car _T \subseteq \car
_U$.
\end{lemma}

\begin{lemma}\label{l:boundaryregion}
Suppose $\bq\in (0,\infty)^I$.  Then the following two conditions
are equivalent:
\begin{enumeratea}
\item
$\bq\in\partial \car$,
\item
$1/W(\bq)=0$ and $1/W(\gl \bq)>0$ for all $\gl\in (0,1)$.
\end{enumeratea}
\end{lemma}

Note that if $T$ is spherical, then $W_T(\bt)$ is a polynomial in $\bt$ and
so $\car_T=(0,+\infty)^I$.  If, for each $i\in I$, $t_i \geq 1$, then $\bt \notin
\car_T$ whenever $W_T$ is infinite.  

Define $ \geps (T):=(-1)^{\Card (T)}$.

\begin{lemma}[{\cite[Ex. 26, pp. 42--43]{bourbaki}, 
\cite[Prop. 26]{serre}, \cite{steinberg}}]\label{l:growth}\hfil
\begin{enumeratei}
\item
Suppose $W$ ($=W_S$) is finite and  let $t_S=t_{w_S}$ be the
monomial corresponding to the element of longest length in $W$.
Then
    \begin{enumeratea}
    \item\label{l:Wfinite}
    \(W(\bt)=t_SW(\bt^{-1}).\)
    \item\label{c:Wfinite}
    \[
    \frac{t_S}{W(\bt)}=\sum_{T \subseteq S}
    \frac{\geps (T)}{W_{T}(\bt)}.
    \]
    \end{enumeratea}
\item\label{l:Treduced}
As in Section~\ref{s:coxeter}, for each $T\subseteq S$, suppose $X_T$ 
denotes
the set of $(\emptyset,T)$-reduced elements in $W$.  Then
\[
W(\bt)= X_T(\bt)W_T(\bt).
\]
\item
As in (\ref{e:defWT}), for each spherical subset $T$ of $S$, suppose
$W^T$ denotes the set of $w\in W$ with $\In (w)=T$.  Then
    \begin{enumeratea}
    \item\label{l:W^T}
    \[
    \frac{W^T(\bt)}{W(\bt)}=\sum_{T' \subseteq T}
    \frac{\geps (T-T')}{W_{S-T'}(\bt)}.
    \]
    \item\label{l:W^Ta}
    \begin{equation*}
    \frac{W^T(\bt)}{W(\bt)}=\sum_{U \in\cs_{\geq T}}
    \frac{\geps(U-T)}{W_{U}(\bt^{-1})}.
    \end{equation*}
    \end{enumeratea}
\item\label{c:Wt-1}
\begin{equation}\notag
\frac{1}{W(\bt^{-1})}=\sum_{T\in\cs} \frac{\geps (T)}{W_T(\bt)}.
\end{equation}
\end{enumeratei}
\end{lemma}

\begin{corollary}\label{c:rational}
\[
W(\bt)= \frac{f(\bt)}{g(\bt)},
\]
where $f,g\in \bZ[\bt]$ are polynomials with integral coefficients.
\end{corollary}

The next lemma follows immediately from the definitions.

\begin{lemma}\label{l:product}
Suppose $(W,S)$ decomposes as a product $(W_1\times W_2, S_1\cup
S_2)$. Then $W(\bt_1,\bt_2)=W_1(\bt_1)W_2(\bt_2)$.  Moreover,
$\car=\car_1\times \car_2$, where $\car$, $\car_1$ and $\car_2$
are the regions of convergence for $W(\bt_1,\bt_2)$, $W_1(\bt_1)$
and $W_2(\bt_2)$, respectively.
\end{lemma}

\begin{example}[The infinite dihedral group]\label{ex:dihedral}
Suppose $S=\{s_1,s_2\}$ and
$m_{s_1s_2}=\infty$, so that $W$ is the infinite dihedral group
$D_{\infty}$.  Its nerve is the $0$-sphere.  Also, suppose $I=\{1,2\}$ and
that $S\to I$ sends $s_j$ to $j$. Using Lemma~\ref{l:growth}(iv), we
compute:
\[
\frac{1}{W(\bt)}=\frac{1-t_1t_2}{(1+t_1)(1+t_2)}.
\]
So, $\car=\{(z_1,z_2)\mid |z_1||z_2|<1\}$.  In particular,
$(0,1)^2\subset \car$.
\end{example}

\begin{example}\label{ex:productdihedral}
Suppose $W=(D_\infty)^n$, the $n$-fold product of infinite dihedral groups.
Its nerve
$L$ is then the $n$-fold join of copies of $S^0$, i.e., it is the boundary
complex of an $n$-octahedron.  By Lemma~\ref{l:product} and
Example~\ref{ex:dihedral},
$(0,1)^I\subset \car$.
\end{example}

\rk{The case of a single indeterminate}  Suppose $I$ is a singleton.
Then $\bt$ is a single indeterminate, call it $t$, the monomial $t_w$
is just $t^{l(w)}$ and $W(t)$ is the
usual growth series.  Let $\rho$ denote its radius of convergence.
An immediate corollary to Lemma~\ref{l:boundaryregion} is the following.

\begin{corollary}\label{l:smallroot}
$1/W(\rho)=0$ and $\rho =\min \{|t|\mid t\in\bC \text{ and } 1/W(t)=0\}$.
\end{corollary}

\comment{
\begin{proof}
Write
    \[
    W(t)=\sum_{n=0}^{\infty}a_nt^n.
    \]
Suppose $|t|\leq \rho$.
Since the $a_n$ are all positive, $|W(t)|\leq \sum a_n|t|^n=W(|t|)$.  If
$|t|<\rho$, then $W(t)$ converges absolutely;
hence, $W(t)\neq \infty$, i.e., $1/W(t)\neq 0$ for $|t|<\rho$.

Since $\rho$ is the radius of convergence, we can find a $t\in \bC$
with $|t|=\rho$ and with $W(t)=\infty$.  Since $|W(t)|\leq W(|t|)$
whenever $|t|\leq \rho$, we
have $W(\rho)=\infty$, i.e., $1/W(\rho)=0$.
\end{proof}
}

A corollary to Lemma \ref{l:product} is the following.

\begin{corollary}\label{l:productsingle}
Suppose $(W,S)$ decomposes as a product $(W_1\times W_2, S_1\cup
S_2)$. Then $W(t)=W_1(t)W_2(t)$ and  $\rho=\min(\rho_1, \rho_2)$,
where $\rho$, $\rho_1$ and $\rho_2$ are the radii of convergence
for $W(t)$, $W_1(t)$ and $W_2(t)$, respectively.
\end{corollary}

In the next proposition we list six other conditions which are equivalent 
to 
the condition that the radius of convergence of $W(t)$ be $1$.

\begin{proposition}\label{p:subexp}
The following conditions on a Coxeter system $(W,S)$ are equivalent.
\begin{enumeratei}
\item
$W$ is amenable.
\item
$W$ does not contain a free group on two generators.
\item
$W$ does not virtually map onto the free group on two generators $F_2$
(i.e., $W$ does not have a finite index subgroup $\gG$ which maps onto
$F_2$).
\item
$W$ is virtually abelian.
\item
$(W,S)$ decomposes as $(W_0\times W_1,S_0\cup S_1)$ where $W_1$ is finite 
and
$W_0$ is a cocompact Euclidean reflection group.
\item
$\rho =1$.
\item
$W$ has subexponential growth.
\end{enumeratei}
\end{proposition}

\begin{proof}
then the
(i)$\implies$(ii) is a standard fact.

(ii)$\implies$(iii).  Suppose for some subgroup $\gG$ of $W$ we have a
surjection $f:\gG\to F_2$ where $F_2$ is the free group on $\{x_1,x_2\}$.
Choose $\gamma_1\in f^{-1}(x_1)$, $\gamma_2\in f^{-1}(x_2)$.  Then
$\langle\gamma_1,\gamma_2\rangle$ is a free subgroup of $W$.

(iii)$\implies$(iv).  It is proved in \cite{mv} (and independently
in \cite{cgon}) that  when $W$ is not virtually abelian there is a
subgroup $\gG$ of finite index in $W$ which maps onto a
nonabelian free group.

(iv)$\implies$(v).  Moussong \cite{moussong} proved
$\gS$ has a $CAT(0)$ metric (so $W$ is a ``$CAT(0)$ group'').
This implies that any abelian subgroup of $W$ is finitely generated.  So,
if $W$ is virtually abelian, then it is virtually free abelian.
We suppose that has $W$ has a rank $n$ free abelian subgroup of finite 
index.
Then $W$ is a virtual $PD^n$-group.  By \cite[Theorem B]{d87},
$W$ decomposes as in (v), where the complex $\gS_0$ for $(W_0,S_0)$ is
a $CAT(0)$ homology $n$-manifold.  By the Flat Torus Theorem (\cite{bh}),
the ``min set '' of the free abelian subgroup on $\gS_0$ is isometric to
$\bR^n$.  Hence, $\gS_0=\bR^n$ and  $W_0$ acts as an isometric
reflection group on it.

(v)$\implies$(vi).  Since a Euclidean reflection group is virtually free
abelian, it has polynomial growth and therefore, the radius of convergence
of its growth series is $1$.  (In fact, the poles of its growth series are
all roots of unity; see Remark~\ref{r:Euclidean} below.)

(vi)$\implies$(vii) is obvious.

(vii)$\implies$(i) by the F{\o}lner condition for amenability.
\end{proof}

\begin{remark}\label{r:Euclidean}
Suppose $W$ is a (cocompact) Euclidean
reflection group.  First consider the case where
$(W,S)$ is irreducible.  Let $W'$ be the
finite linear reflection group obtained by quotienting out the translation
subgroup of $W$ and let $m_1,\dots,m_n$ be the exponents of $W'$.
According to \cite[Ex. 10, p. 245]{bourbaki}, the growth series of $W$ is
given by the following formula of Bott:
    \begin{equation*}
    W(t)=\prod_{i=1}^n\frac{1+t+\cdots +t^{m_i}}{1-t^{m_i}}.
    \end{equation*}
In particular, all the poles of $W(t)$ are roots of unity.
We can reach the same conclusion without the
assumption of irreducibility, since 
the growth series of $(W,S)$ is the product of the
growth series of its irreducible factors.
\end{remark}

\begin{notes} In the case where $\bt$ is a single indeterminate, most of 
the
results of this section come from \cite[Ex. 26, pp. 42--43]{bourbaki}.
The idea of extending the results from this exercise to an $I$-tuple
of indeterminates
comes from \cite{serre}.
Lemma~\ref{l:growth}(iii)(a) is  from \cite[Ex. 26 d), p. 43]{bourbaki},
while (iii)(b) is due to Steinberg \cite{steinberg}.
\end{notes}

\section{Hecke algebras}\label{s:hecke}
Let $A$ be a commutative ring with unit.  Denote by $A^{(W)}$ the free 
$A$-module 
on $W$ consisting of all 
finitely supported functions $W\to A$ and denote by
$A[W]$ this $A$-module equipped with its structure as the group ring of 
$W$.   
Let $(e_w)_{w\in W}$ be the  
standard basis for $A^{(W)}$.  We are primarily interested in the case 
where
$A=\bR$, the field of real numbers.

As in the previous section, $i:S\to I$ is a function such that
$i(s)=i(s')$ whenever $s$ and $s'$ are conjugate.  Let
$\bq=(q_i)_{i\in I}\in A^I$ be a fixed $I$-tuple.  As before,  write $q_s$
for $q_{i(s)}$.  By \cite[Exercise 23, p. 57]{bourbaki}, there is a unique 
ring structure on $A^{(W)}$ such that
\begin{equation}\label{e:heckeleft}
e_se_w=
    \begin{cases}
    e_{sw}, &\text{ if $l(sw)>l(w)$;}\\
    q_se_{sw}+(q_s-1)e_{w}, &\text{ if $l(sw)<l(w)$,}
    \end{cases}
\end{equation}
for all $w\in W$.
We will use the notation $A_{\bq}[W]$ to denote $A^{(W)}$ with this ring
structure.
Note that if $\bq$ is the constant $I$-tuple $\bone:=(1,\dots,1)$, then
$A_{\bq}[W]=A[W]$.  So, 
$A_{\bq}[W]$ is a deformation of the group ring. 
It is called the \emph{Hecke algebra} of $W$ associated to
the \emph{multiparameter} $\bq$.

  From (\ref{e:heckeleft}) it follows that
    \begin{align*}
    e_ue_v&=e_{uv},\text{ for all $u,v\in W$ with $l(uv)=l(u)+l(v)$,
    and}\\
    e_s^2&=(q_s-1)e_s+q_s.
    \end{align*}

The function $e_w\to e_{w^{-1}}$ induces
a linear involution $\ast$ of $A^{(W)}$, i.e.,
    \begin{equation}\label{e:ast}
    \left( \sum a_we_w\right)^{\ast}:=\sum a_{w^{-1}}e_w.
    \end{equation}

\begin{lemma}\label{l:ast}
Formula (\ref{e:ast}) defines an anti-involution of the ring
$A_{\bq}[W]$.  In other
words, for all $x,y\in A_{\bq}[W]$, $(xy)^{\ast}=y^{\ast}x^{\ast}$.
\end{lemma}

\begin{proof}
For each $w\in W$, let $L_w$ (resp. $R_w$) denote left (resp. right)
translation
by $e_w$ defined by $L_w(x)=e_wx$ (resp. $R_w(x)=xe_w$).  A quick 
calculation 
using (\ref{e:heckeleft}) gives:  $R_s=\ast L_s \ast$, for
all $s\in S$.  
If $s_1\cdots s_l$ is a reduced expression for $w$, then
\(
R_w=R_{s_l}\cdots R_{s_1}=\ast L_{s_l}\cdots L_{s_1} \ast =\ast
L_{w^{-1}}\ast
\).  
Therefore,
\(
xe_w=R_w(x)= \ast L_{w^{-1}}\ast (x)=(e_{w^{-1}}x^{\ast})^{\ast}
\).  
Hence, $R_w=\ast L_{w^{-1}} \ast$, for all $w\in W$.
So,
\(
(xe_w)^{\ast}=(e_{w^{-1}}x^{\ast})^{\ast\ast}=e_{w^{-1}}x^{\ast}=
e_w^{\ast}x^{\ast}
\).
The lemma follows.
\end{proof}

Using the involution $\ast$, we deduce the following right hand version of
(\ref{e:heckeleft}):
\begin{equation}\label{e:heckeright}
e_we_s=
    \begin{cases}
    e_{ws}, &\text{ if $l(ws)>l(w)$;}\\
    q_se_{ws}+(q_s-1)e_{w}, &\text{ if $l(ws)<l(w)$.}
    \end{cases}
\end{equation}
\begin{proof}[Proof of (\ref{e:heckeright})]
Apply $\ast$, to get
\begin{align}
(e_we_s)^{\ast}&=e_se_{w^{-1}}\notag\\
&=
    \begin{cases}
    e_{sw^{-1}},    &\text{ if $l(sw^{-1})>l(w^{-1})$;}\\
    q_se_{sw^{-1}}+(q_s-1)e_{w^{-1}}, &\text{ if $l(sw^{-1})<l(w^{-1})$.}
    \end{cases}\notag
\end{align}
Hence,
\begin{align}
e_we_s&=(e_we_s)^{\ast\ast}\notag\\
&=
    \begin{cases}
    e_{ws}, &\text{ if $l(sw)>l(w)$;}\\
    q_se_{ws}+(q_s-1)e_{w}, &\text{ if $l(sw)<l(w)$.}
    \end{cases}\notag
\end{align}
\end{proof}

For each $w\in W$, define $q_w$ by the same formula used to define
$t_w$ 
, i.e.,   if $s_1\cdots s_l$
is a reduced expression for $w$, then
    \begin{equation}\label{e:qw}
    q_w:=q_{s_1}\cdots q_{s_l}.
    \end{equation}
Also, set
    \begin{equation}\label{e:gepsw}
    \geps_w:=(-1)^{l(w)}.
    \end{equation}
The maps $e_w\to q_w$ and $e_w\to \geps_w$ extend  linearly to ring
homomorphisms $A_{\bq}[W]\to A$.

Following Kazhdan--Lusztig \cite{kl}, define an isomorphism
$j_{\bq}:A_{\bq}[W] \to A_{\bq^{-1}}[W]$ by  the formula:
    \begin{equation}\label{e:jiso}
    j_{\bq}(e_w):=\geps_w q_we_w.
    \end{equation}
It is easily checked that $j_{\bq}$ is an algebra homomorphism and that
    \(
    (j_{\bq})^{-1}=j_{\bq^{-1}}.
    \)
Hence, $j_{\bq}$ is an isomorphism of Hecke algebras.  It is called the
\emph{$j$-isomorphism} and denoted simply by $j$ when there is no
ambiguity.

\begin{notes}
Most of the material in this section is taken from
\cite[Exercise 23, p. 57]{bourbaki}.
\end{notes}

\section{Hecke -- von Neumann algebras}\label{s:h-neumann}
From now on $A$ is the field of real numbers $\bR$ and 
$\bq=(q_i)_{i\in I}$ is an $I$-tuple of positive reals.
Define an inner product
$\langle\ ,\ \rangle _{\bq}$ on $\bR^{(W)}$ ($=\R$) by
    \begin{equation}\label{e:innerproduct}
    \langle \sum a_we_w,\sum b_we_w \rangle _{\bq}:=\sum a_wb_wq_w,
    \end{equation}
where $q_w$ was defined by (\ref{e:qw}).
As in \cite{l83}, sometimes it is convenient to normalize 
$(e_w)_{w\in W}$ to an orthonormal basis for $\R$ by setting
    \begin{equation}\label{e:thecke0}
    \te_w:=q_w^{-1/2}e_w.
    \end{equation}

The completion of $\bR^{(W)}$ with respect to the inner
product $\langle\ ,\ \rangle _{\bq}$ is denoted  $\ltwo (W)$, or
simply $\ltwo$, when there is no ambiguity.

\begin{proposition}[{\cite[Proposition 2.1]{dym})}]\label{p:hilbertalgebra}
The inner product defined by
(\ref{e:innerproduct}),  multiplication defined by
equations (\ref{e:heckeleft}) and the anti-involution $\ast$  defined by
(\ref{e:ast}),
give $\R$ a Hilbert algebra structure in the sense of
\cite[A.54]{dix}.  This means, in particular, that
\begin{enumeratei}
\item
$(xy)^{\ast}=y^{\ast}x^{\ast}$,
\item
$\langle x,y\rangle _{\bq}=\langle y^{\ast},x^{\ast}\rangle _{\bq}$,
\item
$\langle xy,z\rangle _{\bq}=\langle y,x^{\ast}z\rangle _{\bq}$,
\item
for any $x\in \R$, left translation by $x$,  $L_x:\R\to\R$, defined
by $L_x(y)=xy$, is continuous,
\item
the products $xy$ over all $x,y\in \R$ are dense in $\R$.
\end{enumeratei}
\end{proposition}

Since the action of $\R$ on
itself by multiplication is continuous, $\ltwo$ is a $\R$-bimodule.

An element $x\in\ltwo$ is \emph{bounded} if right
multiplication by $x$ is bounded on $\R$ (or equivalently, if left
multiplication by $x$ is bounded).  Let $\bR_{\bq}^b[W]$ be the
set of  all bounded elements.

As in \cite{dix} there are two von Neumann algebras associated with this
situation.  They are denoted by  $\cn [W]$ and $\cn'[W]$ or simply by $\cn$
and $\cn'$ when there is no ambiguity.  $\cn$ acts from the right on
$\ltwo$ and $\cn'$ from the left.  Here are two equivalent definitions
of $\cn$:
\begin{enumeratei}
\item
$\cn$ is the algebra of all bounded linear endomorphisms of $\ltwo$ which
commute with the left $\R$-action.
\item
$\cn$ is the weak closure  of $\bR_{\bq}^b[W]$ acting from the right on
$\ltwo$.
\end{enumeratei}
If we interchange the roles of left and right  in the above, we get the two
equivalent definitions of $\cn'$.

\begin{lemma}\label{l:jan1}
If $T\subset S$, then the inclusion $\bR_{\bq}[W_T]\hookrightarrow \R$
induces inclusions $\bR_{\bq}^b[W_T]\hookrightarrow \bR_{\bq}^b[W]$ and
$\cn [W_T]\hookrightarrow \cn$.
\end{lemma}

\begin{proof}
Let $\ltwo (wW_T)\subset \ltwo (W)$ denote the subspace of functions 
which are supported on the coset $wW_T$.  
Then $\ltwo (W)$ decomposes as an orthogonal direct sum
of spaces of the form $\ltwo (wW_T)$.
Suppose $\lambda \in \cn [W_T)$.
Right multiplication by $\lambda$ preserves the summands,
and acts on each summand in the same way. The norm
in the space $\ltwo (wW_T)$ is the norm in $\ltwo (W_T)$ rescaled by
a factor of $q_w^{-1/2}$, so that the operator
norms of right multiplication by $\lambda$ on each of these subspaces is
bounded hence, $\lambda \in \cn$.
\end{proof}

\rk{The $j$-isomorphism}  It follows from the definitions that
the isomorphism $j:\R\to \bR_{\bq^{-1}}[W]$ defined by (\ref{e:jiso})
takes the orthonormal basis $(\te_w)$ for $\ltwo$, defined by
(\ref{e:thecke0}),
to the orthonormal basis $(\te_w)$ for $L^2_{\bq^{-1}}$.
So, it is an isometry.  Therefore, it extends to an isometry of Hilbert
spaces $j:\ltwo \to L^2_{\bq^{-1}}$.
From this, it is obvious that $j$ takes a bounded element of $\ltwo$ to a
bounded element of $L^2_{\bq^{-1}}$.  Hence, it extends to an
isomorphism of von Neumann algebras $j:\cn\to \mathcal{N}_{\bq^{-1}}$.

\rk{The von Neumann trace}
Define the \emph{trace} of an element  $\gf\in \cn$ by
    \begin{equation*}\label{e:tr1}
    \tr_{\cn}(\gf):=\langle e_1\gf,e_1\rangle _{\bq},
    \end{equation*}
where $e_1$ denotes the basis element of $\ltwo$ corresponding to the 
identity
element of $W$.  If
    \(
    \Phi:\bigoplus_{i=1}^n \ltwo \to \bigoplus_{i=1}^n \ltwo
    \)
is a bounded linear map of left $\R$-modules,  then we can represent
$\Phi$ as right multiplication by an $n\times n$ matrix $(\gf _{ij})$ with
entries in $\cn$.  Define
    \begin{equation*}\label{e:tr2}
    \tr_{\cn}(\Phi):=\sum_{i=1}^n \tr_{\cn}(\gf _{ii}).
    \end{equation*}

\rk{Hilbert $\cn$-modules and von Neumann dimension}
\begin{definition}\label{d:hilbert}
A subspace $V$ of a finite orthogonal direct sum of copies of $\ltwo$ is 
called
a \emph{Hilbert $\cn$-module} if it is a closed subspace and if it is
stable under the diagonal left action of $\R$.
\end{definition}

A \emph{map} of Hilbert $\cn$-modules means a bounded linear map of left
$\R$-modules.  A map is \emph{weakly surjective} if it has
dense image; it is a \emph{weak isomorphism} if it is injective and weakly
surjective.

Let $V\subseteq \bigoplus_{i=1}^n \ltwo$ be a Hilbert $\cn$-module and let
$p_V:\bigoplus_{i=1}^n \ltwo \to \bigoplus_{i=1}^n \ltwo$ be the orthogonal
projection onto $V$.  The \emph{von Neumann dimension} of $V$ is the
nonnegative real  number defined by
    \begin{equation}\label{e:defdim}
    \dim_{\cn}V=\tr_{\cn}(p_V).
    \end{equation}
As usual, one shows that $\dim_{\cn}V$ doesn't depend on the choice of
embedding of $V$ into a finite direct sum of copies
of $\ltwo$.  If a subspace $V\subseteq \bigoplus \ltwo$ is $\R$-stable but 
not
necessarily closed, one defines $\dim_{\cn}V:=\dim_{\cn}\ol{V}$.
This dimension function satisfies the usual list of properties:
\begin{enumeratei}
\item
$\dim_{\cn}V=0$ if and only if $V=0$.
\item
For any two Hilbert $\cn$-modules $V$ and $V'$,
\[
\dim_{\cn}(V \oplus V')=\dim_{\cn}V + \dim_{\cn}V'.
\]
\item
$\dim_{\cn}\ltwo=1$.
\item
If $f:V\to V'$ is a weak isomorphism of Hilbert $\cn$-modules, then
$\dim_{\cn}V =\dim_{\cn}V'$.
\item
Suppose that $(W',S')$ and $(W'',S'')$ are Coxeter systems, that
$S'\to I'$ and $S''\to I''$ are indexing functions, that  
$\bq'$ and $\bq''$ are multiparameters,
$S=S'\cup S''$ and $I=I'\cup I''$ are disjoint unions 
that $(W,S)=(W'\times W'',S'\cup S'')$ and that $\bq$
is the multiparameter for $(W,S)$ formed by combining $\bq'$ and $\bq''$.
Let $V'$ (resp. $V''$) be a Hilbert $\cN_{\bq'}[W']$-module (resp.
$\cN_{\bq''}[W'']$-module).  Then the completed tensor product $V:=
V'\otimes V''$ is naturally a Hilbert $\cn$-module and
    \[
    \dim_{\cn [W]}(V'\otimes V'')=(\dim_{\mathcal {N}_{\bq'}[W']} V')
    (\dim_{\mathcal {N}_{\bq''}[W'']} V'').
    \]
\item
Suppose that $T\subset S$ and that $V_T$ is a Hilbert
$\cn [W_T]$-module.
The \emph{induced} Hilbert $\cn$-module $V$ is defined to be
the completed tensor product
    \begin{equation*}
    V:=\ltwo (W)\otimes _{\bR_{\bq}[W_T]} V_T.
    \end{equation*}
Its dimension is given by
    \[
    \dim_{\cn}V=\dim_{\mathcal {N}_{\bq}[W_T]} V_T.
    \]
\end{enumeratei}

\rk{Idempotents in $\cn$ and growth series}  Given a subset
$T$ of $S$,
recall $\car _T$ denotes the region of convergence for $W_T(\bt)$.

\begin{lemma}\label{l:aT}
Given $T\subseteq S$ and $\bq \in \car_T$, there is an
idempotent $a_T \in \cn$ defined by
    \begin{equation*}
    a_T:=\frac{1}{W_{T}(\bq)}\sum_{w\in W_{T}}e_w.
    \end{equation*}
\end{lemma}

\begin{proof}
Define
    \begin{equation}\label{e:aT1}
    \tilde{a}_{T}=\sum_{w\in W_{T}}e_w.
    \end{equation}
Then $\langle \ta_T,\ta_T\rangle _{\bq}=\sum q_w=W_T(\bq)$, so $\ta_T\in 
\ltwo
(W_T)$ if and only if $\bq \in \car_T$.   Assume this.
Recall that for each $s\in S$, $X_s$ denotes the set of
$(\emptyset,\{s\})$-reduced elements in $W$.  Using (\ref{e:heckeright}),
we calculate that for each $s\in T$,
    \begin{align}
    \tilde{a}_{T}e_s&=\sum_{w\in X_s\cap
    W_{T}}e_we_{s}+e_{ws}e_s\notag\\
    &=\sum e_{ws}+q_se_w+(q_s-1)e_{ws}\notag\\
    &=q_s\tilde{a}_{T}\notag
\end{align}
Hence, for $w\in W_{T}$,
    \begin{equation}\label{e:aT3}
    \tilde{a}_{T}e_w=q_w\tilde{a}_{T} \text{ and }
   \tilde{a}_{T}\te_w=q_w^{1/2}\tilde{a}_{T}
    \end{equation}
Therefore,
    \begin{equation}\label{e:aT4}
    (\tilde{a}_{T})^2=W_{T}(\bq)\tilde{a}_{T}.
    \end{equation}
We claim  $\ta_T$ is a bounded element of $\ltwo (W_T)$ (hence, by
Lemma~\ref{l:jan1}, it lies in $\cn$).  To see this, note that
if $x=\sum x_w\te_w \in \bR_{\bq}[W_T]$, then (\ref{e:aT3}) can be 
rewritten
as
    \[
    \ta_T\sum x_w\te_w=\left(\sum x_wq_w^{1/2}\right)\ta_T
    \]
and hence, $\|\ta_T x\|_{\bq}\leq \|\ta_T\|_{\bq}\|x\|_{\bq}$.
So, we get an idempotent defined by
    \begin{equation}\label{e:aT5}
    a_T=\frac{\tilde{a}_{T}}{W_{T}(\bq)}.
    \end{equation}
\end{proof}

\begin{lemma}\label{l:hT}
Given a subset $T$ of $S$ and an $I$-tuple $\bq \in \car_T^{-1}$,
there is an idempotent $h_T \in \cn$ defined by
\begin{equation}\notag 
h_T:=\frac{1}{W_{T}(\bq^{-1})}\sum_{w\in W_{T}}\geps _wq_w^{-1}e_w
\end{equation}
(where $q_w$ and $\geps_w$ are defined by (\ref{e:qw}) and (\ref{e:gepsw}), respectively).
\end{lemma}

\begin{proof}
The proof is similar to the previous one.  Define
    \begin{equation}\label{e:1}
    \tilde{h}_{T}:=\sum_{w\in W_{T}}\geps_wq^{-1}_we_w.
    \end{equation}
Then $\langle \tih_T,\tih_T\rangle _{\bq}=\sum q_w^{-1}=W_T(\bq^{-1})$,
so $\tih_T\in \ltwo (W_T)$ if and only if $\bq ^{-1}\in \car_T$.  
Assume this.
For $s\in T$, we calculate
\begin{align}
\tilde{h}_{T}e_s&=\sum_{w\in X_s\cap W_{T}}
\geps_wq^{-1}_we_we_s+\geps_{ws}q_{ws}^{-1}e_{ws}e_{s}\notag\\
&=\sum\geps_wq^{-1}_we_{ws}+\geps_{ws}q^{-1}_wq^{-1}_s
(q_se_w+(q_s-1)e_{ws})\notag\\
&=-\sum\geps_wq^{-1}_we_w+\geps_{ws}q^{-1}_wq^{-1}_se_{ws}\notag\\
&=-\tilde{h}_{T}.\notag
\end{align}
Therefore, for $w\in W_{T}$,
\begin{align}
\tilde{h}_{T}e_w&=\geps_w\tilde{h}_{T}\text{ and}\label{e:hT3}\\
(\tilde{h}_{T})^2&=\sum_{w\in W_{T}}\geps_wq^{-1}_w\tilde{h}_{T}e_w
=W_{T}(\bq^{-1})\tilde{h}_{T}.\label{e:hT4}
\end{align}
As before, it follows that $\tih_T\in \bR_{\bq}^b[W_T]$ and hence, that
$\tih_T\in\cn$.  So, by (\ref{e:hT4}), we get an idempotent defined by
    \begin{equation}\label{e:hT5}
    h_T:=\frac{\tilde{h}_{T}}{W_{T}(\bq^{-1})}.
    \end{equation}
\end{proof}

Using (\ref{e:heckeright}) we get the following right hand versions of
(\ref{e:aT3}) and (\ref{e:hT3}) for $T\subseteq S$ and $w\in W_T$:
    \begin{align}
    e_wa_T&=q_wa_T \text{ and}\label{e:aT3'}\\
    e_wh_T&=\geps_wh_T.\label{e:hT3'}
    \end{align}

What is the effect of the $j$-isomorphism on these idempotents?  It
follows immediately from definitions (\ref{e:jiso}), (\ref{e:aT1}) and
(\ref{e:1}) that
    \[
    j(\tilde{a}_T)=\tilde{h}_T \text{ and }
    j(\tilde{a}_T)=\tilde{h}_T.
    \]
Hence, by the definitions in Lemmas~\ref{l:aT} and \ref{l:hT},
    \begin{equation}\label{e:jajh}
    j(a_T)=h_T \text{ and } j(h_T)=a_T.
    \end{equation}

Using (\ref{e:aT3}), (\ref{e:hT3}), (\ref{e:aT3'}) and (\ref{e:hT3'}), we
easily calculate that for any $U\subseteq T\subseteq S$:
\begin{align}
a_Ua_T&=a_T=a_Ta_U \text{ whenever } \bq\in \car_U,\label{e:aTaU}\\
h_Uh_T&=h_T=h_Th_U \text{ whenever } \bq\in 
\car_U^{-1}.\label{e:hThU}
\end{align}

If $s_1\cdots s_l$ is a reduced expression for $w$, then $s_l\cdots s_1$
is a reduced expression for $w^{-1}$.  It follows that
    \begin{equation*}\label{e:w-1}
    q_{w^{-1}}=q_w \text{ and }\geps _{w^{-1}}=\geps _w.
    \end{equation*}
So, 
    \begin{equation*}
    a_T^{\ast}=a_T \text{ and } h_T^{\ast}=h_T,
    \end{equation*}
whenever the idempotents $a_T$ and $h_T$ make sense.  In other words, the
maps $x\to xa_T$ and $x\to xh_T$ are orthogonal projections from $\ltwo$
onto Hilbert submodules.

\begin{Remark}
The ``a'' in $a_T$ is for ``average,'' while the ``h'' in $h_T$ is for
``harmonic.''
\end{Remark}

\begin{definition}\label{d:characters}
For each $T\subseteq S$, let $\ga_T:\bR[W_T]\to \bR$ and
$\gb_T:\bR_{\bq}[W_T]\to
\bR$ be the algebra homomorphisms defined by $e_w\to q_w$
and $e_w\to \geps _w$, respectively.  $\ga_T$ is the \emph{symmetric character} and
$\gb_T$ is the \emph{alternating character}.
\end{definition}

The next lemma follows immediately from equations (\ref{e:aT3}) and
(\ref{e:hT3}).

\begin{lemma}\label{l:characters}\hfil
\begin{enumeratei}
\item
Supposing $\bq\in\car_T$, the action of  $\bR_{\bq}[W_T]$ on $\ltwo a_T$ by
right multiplication is via the character $\ga_T$.
\item
Supposing $\bq ^{-1}\in\car_T$, the action of  $\bR_{\bq}[W_T]$ on $\ltwo 
h_T$
by right multiplication is via the character $\gb_T$.
\end{enumeratei}
\end{lemma}

\rk{Some Hilbert $\cn$-submodules of $\ltwo$}  
To simplify notation, for each $s\in S$, write $a_s$ and $h_s$ for
the idempotents $a_{\{s\}}$ and $h_{\{s\}}$.  Let $A_s=\ltwo
a_s$ and $H_s=\ltwo h_s$ be the corresponding Hilbert $\cn$-submodules of
$\ltwo$.

\begin{lemma}\label{l:boris0}
For each $s\in S$, the subspaces $A_s$ and $H_s$ are the orthogonal
complements of each other in $\ltwo$.
\end{lemma}

\begin{proof}
\begin{align}
a_s+h_s&=\frac{1}{1+q_s}(1+e_s)+\frac{1}{1+q_s^{-1}}(1-q_s^{-1}e_s)\notag\\
&=1.\notag
\end{align}
So, $a_s$ and $h_s$ are orthogonal projections onto complementary 
subspaces.
\end{proof}

For each $T\subseteq S$, set
\begin{equation}\label{e:ATHT}
A_T:=\bigcap _{s\in T} A_s \text{ and } H_T:=\bigcap _{s\in T} 
H_s.
\end{equation}

For any subspace $E\subset \ltwo$, let $E^{\perp}$ denote its orthogonal
complement.  Since $\perp$ takes sums  to intersections and intersections 
to
closures of sums:
\begin{align}
\Big(\sum_{s\in T}A_s\Big)^{\perp}&=H_T,&
\Big(\sum_{s\in T}H_s\Big)^{\perp}&=A_T,\label{e:sumperp}\\
\ol{\sum_{s\in T}A_s}&=(H_T)^{\perp},&
\ol{\sum_{s\in T}H_s}&=(A_T)^{\perp}.\label{e:intersectperp}
\end{align}

\begin{lemma}\label{l:AS}
Let $A_S$ be the subspace of $\ltwo$ defined in (\ref{e:ATHT}).
\begin{enumeratei}
\item
For all $x\in A_S$ and $w\in W$, $xe_w=q_wx$.
\item
If $\bq \notin\car$, then $A_S=0$.
\item
If $\bq \in \car$, then $A_S$ is the line spanned by $a_S$
and $\ltwo a_S=A_S$.  Hence,
\[
\dim_{\cn}A_S=\frac{1}{W(\bq)}.
\]
\end{enumeratei}
\end{lemma}

There is also a version of this lemma for $H_S$.

\begin{lemma}\label{l:HS}
Let $H_S$ be the subspace of $\ltwo$ defined in (\ref{e:ATHT}).
\begin{enumeratei}
\item
For all $x\in H_S$ and $w\in W$, $xe_w=\geps_wx$.
\item
If $\bq ^{-1}\notin\car$, then $H_S=0$.
\item
If $\bq ^{-1}\in \car$, then $H_S$ is the line spanned by $h_S$
and $\ltwo h_S=H_S$.  Hence,
\[
\dim_{\cn}H_S=\frac{1}{W(\bq ^{-1})}.
\]
\end{enumeratei}
\end{lemma}

We prove only the first version, the proof of the second version being
entirely similar.

\begin{proof}[Proof of Lemma~\ref{l:AS}]
(i)  As in Definition~\ref{d:characters}, $\ga_T:\bR[W_T]\to \bR$ denotes 
the
symmetric character.  The $\ga_{\{s\}}$-eigenspace of
$\bR_{\bq}[W_{\{s\}}]$  on
$\ltwo$ is $\Ker (q_s-e_s)=\Ker h_s =\ltwo a_s=A_s$.  Since the
subalgebras $\bR_{\bq}[W_{\{s\}}]$ generate $\R$,
the intersection of the $A_s$, $s\in S$,  is the
$\ga_S$-eigenspace for $\R$.

(ii)  If $x=\sum x_we_w \in A_S$, then
\[
q_wx_w=\langle e_w,x\rangle _{\bq}=\langle 1,xe_w^{\ast}\rangle _{\bq}= 
\langle
1,xe_{w^{-1}} \rangle _{\bq} =\langle 1,q_wx\rangle _{\bq}=q_wx_1.
\]
In other words, the coefficients $x_w$ are all equal.  Hence, $\langle
x,x\rangle_{\bq}=x_1^2W(\bq)$.  So, if $\bq\notin \car$, $x\notin
\ltwo$ unless $x=0$ and if $\bq\in \car$, $x$ must be a scalar multiple
of $a_S$.

(iii)  By Lemma~\ref{l:characters}, if $\bq\in\car$, then $\ltwo a_S
\subseteq A_S$.  Since $a_S\neq 0$ and $\dim_{\bR}A_S=1$, the inclusion is 
an equality.  Hence,
\[
\dim_{\cn}A_S=\dim_{\cn}\ltwo a_S=\tr_{\cn}a_S=\frac{1}{W(\bq)}.
\]
\end{proof}

\begin{corollary}\label{c:ATHT2}
For any $T\subseteq S$:
\begin{enumeratei}
\item
$A_T=\ltwo a_T$ if $\bq\in \car_T$ and $A_T=0$ if $\bq\notin \car_T$.
\item
$H_T=\ltwo h_T$ if $\bq^{-1}\in \car_T$ and  $H_T=0$ if $\bq^{-1}\notin 
\car_T$.
\end{enumeratei}
\end{corollary}

\begin{proof}
Since $A_T$ and $H_T$ are induced from Hilbert $\bR_{\bq} [W_T]$-modules,
this follows from Lemmas~\ref{l:AS} and \ref{l:HS}.
\end{proof}

\section{Some cell complexes}\label{s:spaces}
\rk{The basic construction}
Suppose we are given the following data:
\begin{itemize}
\item
a Coxeter system $(W,S)$,
\item
a CW complex $Z$ and
\item
a family of subcomplexes $(Z_s)_{s\in S}$.
\end{itemize}
The $Z_s$ are called the \emph{mirrors} of $Z$.
Given these data there is a classical construction of a 
CW complex $\cu=\cu (W,Z)$
with a $W$-action so that $Z$ is a strict fundamental domain.  We recall 
the construction.  

For each subset $T$ of $S$, set
    \begin{align}
    Z_T&:=Z\cap \bigcap_{s\in T}Z_s,\notag\\
    Z^T&:=\bigcup_{s\in T} Z_s.\label{e:Z^T}
    \end{align}
For each cell $c$ of $Z$ and each point $z\in Z$, set
    \begin{align}
    S(z)&:=\{s\in S\mid z\in Z_s\},\label{e:Sz}\\
    S(c)&:=\{s\in S\mid c\subseteq Z_s\}.\label{e:Sc}
    \end{align}
Define $\cu (W,Z):=(W\times Z)/\sim$ where $\sim$ is the equivalence relation
defined by:  $(w,z)\sim (w',z')$ if and only if $z=z'$  and the
cosets $wW_{S(z)}$ and $w'W_{S(z)}$ are equal.  Write $[w,z]$ for the
image of $(w,z)$ in $\cu$.  The group $W$ acts on $\cu$ via $w\cdot
[w',z]=[ww',z]$.  The
orbit space is $Z$.   Identifying  $Z$ with the image of $1\times Z$
in $\cu$, we see that $Z$ is a strict fundamental domain.
$wZ$, the translate of $Z$ by $w$, is identified with the image of $w\times Z$.
The CW structure
on $\cu$ is defined by declaring  the family $(wc)$, with $w\in W$
and $c$ a cell
of $Z$, to be the set of cells in $\cu$.  (Note that $wc$ is the image
of $w\times c$ in $\cu$.)

The setwise stabilizer of a cell $c$ of $Z$ is the special subgroup
$W_{S(c)}$.  Moreover, $W_{S(c)}$ fixes each point of $c$.

The family $(Z_s)_{s\in S}$ is \emph{$W$-finite} if $Z_T=\emptyset$
whenever $W_T$ is infinite.  This  condition insures that each isotropy
subgroup is finite.  It is
equivalent to the condition that $W$ act properly on $\cu$.  We shall assume it throughout.

\rk{The complex $\gS$} The complex $\gS$ can
be described in terms of the basic construction.  As in  Section~\ref{s:coxeter}, denote the
geometric realization of the poset $\cs$ by $K$ and the geometric 
realization of $W\cs$ by $\gS$. 
For each $s\in S$, let 
$K_s$  be the geometric realization of the subposet
$\cs_{\geq \{s\}}$.  It is a subcomplex of $K$.
The space $\cu (W,K)$ is naturally a simplicial
complex.  The natural map $W\times \cs\to W\cs$, defined by $(w,T)\to
wW_T$, induces a map of geometric realizations $W\times K\to \gS$ and this
descends to W-equivariant  map $\cu (W,K) \to \gS$.  As in \cite{d83}, it 
is easily  seen  that this map is a simplicial isomorphism, i.e.,
    \begin{equation}\label{e:gscu}
    \gS\cong \cu (W,K).
    \end{equation}

\rk{Cellulation of $\gS$ by Coxeter cells}
As is explained in \cite{moussong, d98a, d02} and below, $\gS$ has
another cell structure: its  cellulation by 
``Coxeter cells.'' 

Suppose, for the moment, that $W$ is finite and $\Card (S)=n$.  Associated
to $(W,S)$ there is
a $n$-dimensional convex polytope $P$ called the
\emph{Coxeter cell} of type $W$.  $P$ is defined as the convex
hull of a generic $W$-orbit in the canonical representation of $W$ on
$\bR^n$.  $W$ acts
simply transitively on the vertex set of $P$; moreover, a subset of vertices spans
a face if and only if it has the form $wW_Tv_0$ for some special 
coset $wW_T$ and  for a given choice
of base vertex $v_0$ in the interior of the fundamental simplicial cone.
This identifies 
the face poset of $P$ with $W\cs$.  In other words,
it gives a simplicial isomorphism between
$\gS$ and the barycentric subdivision of $P$.

Returning to the case where $(W,S)$ is arbitrary, note that for any
element $wW_T\in W\cs$, the poset $W\cs_{\leq wW_T}$ is identified
with the face poset of $P_T$, the Coxeter cell of type $W_T$.  So, the
subcomplex $\vert W\cs_{\leq wW_T}\vert$ of $\gS$ is identified
with the barycentric subdivision of $P_T$.  This defines the cell
structure on $\gS$:  each simplicial subcomplex $\vert W\cs_{\leq wW_T}\vert$ is
identified with a Coxeter cell of type $W_T$.  So, the vertex
set of $\gS$ is $W$ and a subset of $W$ is the vertex set of a cell if and
only if it is a coset  $wW_T$ for some $w\in W$ and $T\in \cs$.
We shall use the notation $\gS_{cc}$ to denote $\gS$ equipped
with this cell structure, where the subscript $cc$ stands for ``Coxeter
cell.''  (In \cite{dym} this cell structure is denoted $\gS_d$, where the
subscript $d$ stood for ``dual cell.'')  The poset of cells of $\gS_{cc}$ 
is $W\cs$.

Suppose $U\subseteq S$.  Let $\cs (U):=\{T\in \cs\mid T\subseteq U\}$.   
Define $\gS (U)$ to be the subcomplex of $\gS_{cc}$ consisting of all 
Coxeter
cells of type $T$, with $T\in \cs (U)$.  If $K(U):=\gS (U)\cap K$, then
it is not difficult to see that
    \begin{equation}\label{e:sigmaU}
    \gS (U)=\cu (W,K(U))=W\times _{W_U} \cu (W_U,K(U)).
    \end{equation}
Moreover, $\cu (W_U,K(U))$ equivariantly deformation retracts onto
the complex $\gS_{W_U}$ associated to $(W_U,U)$.

\rk{Ruins}
Given $U\subseteq S$ and $T\in \cs (U)$, define three subcomplexes of $\gS
(U)$:
\begin{align}
\gO (U,T): &\text{ the union of closed cells of type $T'$, with $T'\in 
\cs
(U)_{\geq T}$,}\notag\\
\wh{\gO}(U,T): &\text{ the union of closed cells of type $T''$, 
$T''\in \cs
(U)$, $T''\notin \cs (U)_{\geq T}$,}\notag\\
\partial \gO (U,T): &\text{ the cells of $\gO (U,T)$ of type $T''$,
with $T''\notin \cs (U)_{\geq T}$}.\notag
\end{align}
$\gO (U,T)$ is the union of all cells of type $T''$,
where $T''\leq T'$ for some $T'\in \cs(U)_{\geq T}$.
So,
\begin{align}
\partial \gO (U,T)&=\gO (U,T)\cap \wh{\gO} (U,T)\text{ and}\notag\\
\gS (U)&=\gO (U,T)\cup \wh{\gO} (U,T).\notag
\end{align}
The pair $(\gO (U,T),\partial\gO (U,T))$ is called the  
\emph{$(U,T)$-ruin}.
For example, for $T=\emptyset$, we have $\gO (U,\emptyset)=\gS (U)$ and
$\partial\gO (U,\emptyset)=\emptyset$.  The key step in our proofs of the
main results in Sections~\ref{s:decomp} and \ref{s:decouple} is the computation of certain homology groups of 
$(U,T)$-ruins.

Similarly, define
\begin{align*}
K(U,T)&:=\gO (U,T)\cap K, \\
\partial K(U,T)&:=\partial\gO (U,T)\cap K,\\
\wh{K}(U,T)&:=\wh{\gO} (U,T)\cap K.
\end{align*}
So, that
\begin{align*}
\gO (U,T)&=\cu (W,K(U,T)),\\
\partial\gO (U,T)&=\cu (W,\partial K(U,T)),\\
\wh{\gO} (U,T)&=\cu (W,\wh{K}(U,T)).
\end{align*}

\section{Weighted $L^2$-(co)homology}\label{s:weighted}
Notation is as in the previous section: $Z$ is a CW complex, $(Z_s)_{s\in 
S}$ is a $W$-finite family of subcomplexes and $\cu=\cu (W,Z)$.

We begin by defining a chain complex of Hilbert $\cn$-modules for the CW
complex $\cu$.  In this case, each orbit of cells contributes an
$\cn$-module of the form $A_T$ for some $T\in \cs$.  Next  we define 
chain complexes of $\cn$-modules in the cases of the cellulation of $\gS$
(and its subcomplexes of ruins) by Coxeter cells.  In these cases each orbit
of cells contributes a $\cn$-module of the form $H_T$, $T\in \cs$.

\rk{Weighted (co)chain complexes for $\cu (W,Z)$}
Orient the cells of $Z$ arbitrarily and then orient the remaining cells
of $\cu$ so that for each positively oriented cell $c$ of $Z$ and each
$w\in W$, $wc$ is positively oriented.

As usual, $\bq$ is an $I$-tuple of positive real numbers.
Given a cell $c$ of $Z$ define a measure $\mu_{\bq}$ on its orbit $Wc$ by
    \begin{equation}\label{e:muq}
    \mu_{\bq}(wc):=q_u,
    \end{equation}
where $u$ is the shortest element in the coset
$wW_{S(c)}$ (i.e., $u$
is the $(\emptyset,S(c))$-reduced element in this coset).
This extends in a natural way to a measure, also denoted by $\mu
_{\bq}$, on $\cu^{(i)}$ (where $\cu^{(i)}$ denotes
the entire set of $i$-cells in $\cu$).  As in
\cite{dym}, define the \emph{$\bq$-weighted $L^2$-(co)chains} on $\cu$
(in dimension $i$) to be the Hilbert space
    \begin{equation}\label{e:(co)chains}
    \ltwo C_i(\cu)=\ltwo C^i(\cu):=L^2(\cu ^{(i)},\mu_{\bq}).
    \end{equation}
We have coboundary and boundary maps, $\gd^i:\ltwo C^i(\cu)\to
\ltwo C^{i+1}(\cu)$ and $\partial_i:\ltwo C_i(\cu)\to \ltwo C_{i-1}(\cu)$
defined by the usual formulas:
    \begin{align}
    \gd^i(f)(\gamma)&:=\sum[\gb:\gamma]f(\gb)\
    \label{e:coboundary}\\
    \partial_i(f)(\ga)&:=\sum[\ga:\gb]f(\gb),\label{e:partial}
    \end{align}
where the first sum is over all $i$-cells $\gb$ incident to the
$(i+1)$-cell $\gamma$ while the second is over all $\gb$ whose boundary
contains the $(i-1)$-cell $\ga$.  In contrast to  the standard situation (where $\bq=\bone$), the maps
$\gd^i$ and $\partial_{i+1}$ are not adjoint to one another.  
Define $\partial_i^{\bq}:\ltwo C_i(\cu)\to \ltwo 
C_{i-1}(\cu)$ by
\begin{align}
\partial_i^{\bq}(f)(\ga)&:=\sum[\ga:\gb]\muq (\gb)\muq (\ga)^{-1}f(\gb).
\label{e:partialq}
\end{align}
A quick calculation (cf. \cite[Section 1]{dym}) then shows that
$\gd ^{\ast}=\partial ^{\bq}$.  
Since $\gd^2=0$, taking  adjoints, we get $(\partial^{\bq})^2=0$.
Hence, $(\ltwo C_{\ast}(\cu),\partial ^{\bq})$ is also a chain complex.

One defines the \emph{$\bq$-weighted $L^2$-(co)homology} of $\cu$
in dimension $i$ by
\begin{align*}
\ltwo H^i(\cu)&:=H^i((\ltwo C^{\ast}(\cu),\gd))=\Ker \gd^i/\Ima
\gd^{i-1},\\
\ltwo H_i(\cu)&:=H_i((\ltwo C_{\ast}(\cu),\partial^{\bq}))=
\Ker \partial^{\bq}_i/\Ima \partial^{\bq}_{i+1}.
\end{align*}
Notice that while we are using the ordinary coboundary map $\gd$,  the
boundary map $\partial^{\bq}$ is not the usual one:  it is modified by
coefficients depending on $\bq$.
There is a standard problem with these (co)homology  groups:  the 
quotients need not be Hilbert spaces.  To remedy this, define 
\emph{reduced} weighted $L^2$-(co)homology by
\begin{align*}
\ltwo \ch^i (\cu)&:=\Ker \gd^i/\ol{\Ima \gd^{i-1}},\\
\ltwo \ch_i (\cu)&:=\Ker \partial^{\bq}_i/\ol{\Ima \partial^{\bq}_i}.
\end{align*}
Since $\gd^{\ast}=\partial ^{\bq}$ and $(\partial ^{\bq})^{\ast}=\gd$, we 
have
the Hodge decomposition:
\[
\ltwo C^i(\cu)=(\Ker \gd^i\cap\Ker \partial^{\bq}_i)\oplus \ol{\Ima
\gd^{i-1}}\oplus \ol{\Ima \partial ^{\bq}_{i+1}}.
\]
It follows that both $\ltwo\ch^i(\cu)$ and $\ltwo\ch_i(\cu)$ can be
identified with the space $\Ker \gd^i\cap\Ker \partial^{\bq}_i$ of
harmonic cochains.  In particular, $\ltwo\ch^i(\cu)\cong\ltwo\ch_i(\cu)$.

\begin{lemma}\label{l:qtoq-1}
The chain complexes $(\ltwo C_{\ast}(\cu),\partial ^{\bq})$ and
$(L^2_{\bq^{-1}} C_{\ast}(\cu),\partial)$ are isomorphic.
\end{lemma}

\begin{proof}
For a chain $f$ on $\cu$, define another chain $\theta (f)$ by
$\theta (f) (\gb) :=\muq (\gb)f(\gb)$ and
note that $\theta (f)$ is $\bq$-square summable if and only if $f$
is $\bq^{-1}$-square summable.  Hence, it defines a linear isomorphism
$\theta=\theta_\bq:L^2_{\bq^{-1}} C_{\ast}(\cu)\to \ltwo C_{\ast}(\cu)$.  Using
(\ref{e:partial}) and (\ref{e:partialq}),  computation shows that 
$\theta\circ
\partial = \partial^{\bq}\circ \theta$.  So, $\theta$ is a chain
isomorphism.
\end{proof}

\begin{remark}\label{r:qtoq-1}
We have canonical inclusions of chain complexes:
\begin{equation}\label{e:chain1}
C_*(\cu;\bR)\hookrightarrow (\ltwo C_*(\cu),\partial)
\hookrightarrow C_*^{lf}(\cu;\bR).
\end{equation}
So, using the isomorphism $\theta_{\bq^{-1}}$ of Lemma~\ref{l:qtoq-1} we get inclusions:
\begin{equation}\label{e:chain2}
C_*(\cu;\bR)\hookrightarrow L^2_{\bq^{-1}} C_*(\cu)
\hookrightarrow C_*^{lf}(\cu;\bR).
\end{equation}
Similarly, we have inclusions of cochain complexes:
\begin{equation}\label{e:chain3}
C^*_c(\cu;\bR)\hookrightarrow \ltwo C^*(\cu)
\hookrightarrow C^*(\cu;\bR).
\end{equation}
(Here $C_*^{lf}(\ )$ and $C^*_c(\ )$ stand for, respectively,
infinite cellular chains and
finitely supported cellular cochains.)  The second map in (\ref{e:chain1})
(or the second map in (\ref{e:chain2})) is obtained by dualizing the first
map in (\ref{e:chain3}).  Similarly, the second map in
(\ref{e:chain3}) is obtained
by dualizing the first map in (\ref{e:chain1}).

As was indicated in the Introduction and as will be explained further in  Section~\ref{s:ordinary}, for $\bq\in\car$,
the first maps in (\ref{e:chain1}) and (\ref{e:chain2})
induce monomorphisms with dense image
\begin{align*}
H_i(\cu;\bR)&\hookrightarrow \ch_i(L^2_{\bq^{-1}}C_*(\cu),\partial)\\
H_i(\cu;\bR)&\hookrightarrow \ltwo\ch_i(\cu).
\end{align*}
(The first  monomorphism agrees with one's intuition.)
Similarly, for $\bq\in\car^{-1}$,
the first map in (\ref{e:chain3}) induces a monomorphism with
dense image
\[
\can:H^i_c(\cu;\bR)\hookrightarrow \ltwo\ch^i(\cu).
\]

Dualizing we get isomorphisms:
\begin{align*}
\ltwo\ch^i(\cu)&\mapright{\cong} H^i(\cu;\bR) \text{ for $\bq\in
\car$ and}\\
\ltwo\ch_i(\cu)&\mapright{\cong} H_i^{lf}(\cu;\bR) \text{ for $\bq\in
\car^{-1}$}.
\end{align*}
All this is reminiscent of a well-known result of Cheeger-Gromov 
\cite{cg86} that if a discrete amenable group $A$ acts properly on a CW 
complex $X$, then the canonical map $L^2H^*(X)\to H^*(X;\bR)$ is 
injective. So, for $\bq\in \car$, weighted $L^2$-cohomology behaves as 
if $W$ were amenable.
\end{remark}

\rk{The Hilbert $\cn$-module structure on $\ltwo C^{\ast}(\cu)$}
Following \cite{dym},  realize $\ltwo$ as $L^2(W,\nu_{\bq})$, where
$\nu_{\bq}$ is the measure on $W$ defined by $\nu_{\bq}(w)=q_w$.  For each
subset $T$ of $S$, the
Hilbert $\cn$-submodule $A_T\subset \ltwo$, defined by (\ref{e:ATHT}),
is then identified with
$L^2(W,\nuq)^{W_T}$, the subspace of $L^2$ functions which are constant on
each right coset $wW_T$.

Since each cell of $\cu$ has the form $wc$ for some cell $c$ of $Z$ and
some $w\in W$, we have
    \begin{equation*}
    \ltwo C^i(\cu)=\bigoplus _{c\in Z^{(i)}}L^2(Wc,\mu_{\bq}),
    \end{equation*}
where the sum ranges over all $i$-cells $c$ of $Z$.  Moreover,
$L^2(Wc,\mu_{\bq})$ can be identified with $A_{S(c)}$ via the isometry
$\phi_c:L^2(Wc,\mu_{\bq})\to A_{S(c)}$ defined by
    \begin{equation*}\label{e:phi}
    \phi _c(f)=\sqrt{W_{S(c)}(q)}\Big(\sum _{u\in X_{S(c)}}f(uc)e_u
    a_{S(c)}\Big),
    \end{equation*}
where the summation is over all $(\emptyset,S(c))$-reduced elements $u$ and
where $a_{S(c)}\in \cn$ is the idempotent defined in Lemma~\ref{l:aT}.
So, we get an isometry
    \begin{equation*}
    \ltwo C^i(\cu)=\bigoplus_{c\in Z^{(i)}}L^2(Wc,\mu_{\bq})
    \mapright{\cong} \bigoplus_{c\in Z^{(i)}}A_{S(c)}.
    \end{equation*}
Since each $A_{S(c)}$ is a left $\R$-submodule of $\ltwo$,
this gives $\ltwo C^i(\cu)$ the structure of a Hilbert $\cn$-module as in
Definition~\ref{d:hilbert} (provided we assume, as we shall,
that $Z$ is a finite complex).
It also gives an isometric embedding
    \begin{equation}\label{e:Phi}
    \Phi:\ltwo C^i(\cu)\hookrightarrow \bigoplus_{c\in Z^{(i)}}\ltwo
    =C^i(Z)\otimes \ltwo.
    \end{equation}

It is shown in \cite[Lemma 3.2]{dym} that $\gd$ and $\partial^{\bq}$ are
maps of Hilbert $\cn$-modules.  (It is \emph{not} true that $\gd_\bq$ and
$\partial$ are maps of Hilbert $\cn$-modules; however, it is possible to
give $\ltwo C^{\ast}(\cu)$ and $\ltwo C_{\ast}(\cu)$ the structure of
Hilbert $\bR_{\bq^{-1}}[W]$-modules so that they are maps of Hilbert
$\bR_{\bq^{-1}}[W]$-modules.  To do this, one transports the
$\bR_{\bq^{-1}}[W]$-module structure from $L^2_{\bq^{-1}} C_{\ast}(\cu)$ 
via the isomorphism $\theta$ of Lemma~\ref{l:qtoq-1}.)
It follows that $\Ker \gd$, $\Ker \partial^{\bq}$,
$\ol{\Ima \gd}$ and $\ol{\Ima \partial ^{\bq}}$ are
Hilbert $\cn$-modules.  Hence, $\ltwo \ch ^i(\cu)$ (or $\ltwo \ch _i(\cu)$)
is also a Hilbert $\cn$-module.

\rk{Weighted (co)chain complexes for cellulations by Coxeter cells}
Let $\angt$ denote the Coxeter cell in $\gS$ corresponding to $W_T\in
W\cs$ (the
$W_T$-coset of the identity element).  Then $W\angt$, the $W$-orbit of
$\angt$, is the set of all Coxeter cells in $\gS$ of type $W_T$.  Define a
measure $\muq$ on $\gS_{cc}^{(i)}$ by $\muq (w\angt)=q_u$, where
$u=p_T(w)$ is the shortest element in $wW_T$.  Define the
\emph{$\bq$-weighted $L^2$-(co)chains} on $\gS$ (in dimension $i$) to be
the Hilbert space
    \begin{equation*}
    \ltwo C_i(\gS_{cc})=\ltwo C^i(\gS_{cc}):=L^2(\gS ^{(i)}_{cc},\muq).
    \end{equation*}
We have
    \begin{equation*}
    \ltwo C_i(\gS_{cc})=\bigoplus_{T\in \cs^{(i)}}L^2(W\angt,\muq).
    \end{equation*}

Choose arbitrarily a orientations for cells of the form $\angt$,
$T\in\cs$.  We  use the following orientation convention for the
remaining cells in
$W\angt$:  if $u\in X_T$ (i.e., if $u$ is $(\emptyset,T)$-reduced as
defined in Section~\ref{s:coxeter}), then orient
$u\angt$ so that
left translation by $u$ is an orientation-preserving map $\angt\to u\angt$.

As in (\ref{e:coboundary}), $\gd:\ltwo C^i(\gS_{cc})\to \ltwo
C^{i+1}(\gS_{cc})$ is the usual coboundary map.  Its adjoint
$\partial^{\bq}:\ltwo C_{i+1}(\gS_{cc})\to \ltwo C_i(\gS_{cc})$ is defined
similarly to (\ref{e:partialq}).

Next, we determine the formula for the restriction of
$\partial^\bq$ to the summand $L^2(W\langle U\rangle,\muq)$, 
where $U\in\cs^{(i+1)}$.  Any $w\in W$ can be uniquely decomposed as $w=uv$ with $u\in X_U$ and $v\in W_U$.
Suppose $T\in
\cs^{(i)}$ is obtained by deleting one element of $U$ and $w\in X_T$.  If $w\in X_T$, then $v\in W_U\cap X_T$.  For any $f\in L^2(W\langle
U\rangle,\muq)$, we have the following formula for $\partial^\bq$:
    \begin{equation}\label{e:partialT}
    \partial^\bq f (w\angt)=\geps_vq_v^{-1}f(u\langle U\rangle),
    \end{equation}
where $w=uv$ as above.

The group $W_T$ acts nontrivially on the cell $\angt$.  In fact, $v\in W_T$
is $\geps_v$ orientation-preserving.  Hence, the right
$\bR_{\bq}[W_T]$-action on $L^2(W\angt,\muq)$ is via the alternating 
character
$\gb_T$ of Definition~\ref{d:characters}. Therefore, $L^2(W\angt,\muq)$
can be identified with $H_T$.   A specific isometry 
$\psi:L^2(W\angt,\muq)\to
H_T$ can be defined by
    \begin{equation}\label{e:psi}
    \psi _T(f)=\sqrt{W_T(q^{-1})}\Big(\sum _{u\in X_T}f(u\angt)e_u \Big) 
h_T,
    \end{equation}
where $h_T$ is the idempotent of $\cn$ defined in Lemma~\ref{l:hT}.
So, we have an isometry:
    \begin{equation}\label{e:Coxchains2}
    \ltwo C_i(\gS_{cc})=\bigoplus_{T\in \cs^{(i)}}L^2(W\angt,\mu_{\bq})
    \mapright{\cong} \bigoplus_{T\in \cs^{(i)}}H_T.
    \end{equation}
Since each $H_T$ is a left $\R$-submodule of $\ltwo$,
this gives $\ltwo C^i(\gS_{cc})$ the structure of a Hilbert $\cn$-module.
It also gives an isometric embedding
\[
\Psi:\ltwo C_i(\gS_{cc})\hookrightarrow \bigoplus_{c\in \cs^{(i)}}\ltwo
=C_i(K)\otimes \ltwo.
\]

We use the isomorphism in (\ref{e:Coxchains2}) to
transport the Hilbert $\cn$-module structure from the right hand side of
(\ref{e:Coxchains2}) to $\ltwo C_i(\gS_{cc})$.
It is proved in \cite[Lemma 4.3]{dym} that $\gd$ and $\partial^\bq$ are
maps of Hilbert $\cn$-modules. We shall give the argument in
Lemma~\ref{l:chainruin} below.
Hence, we get reduced $L^2$-(co)homology groups:
\begin{equation*}
\ltwo \ch^i (\gS_{cc})=\Ker \gd^i/\ol{\Ima \gd^{i-1}} \text{ and }
\ltwo \ch_i (\gS_{cc})=\Ker \partial^{\bq}_i/\ol{\Ima \partial^{\bq}_i},
\end{equation*}
which are also  Hilbert $\cn$-modules.  It is proved in \cite[Section 
5]{dym}
that the (co)homology groups of $\ltwo C_{\ast}(\gS_{cc})$ are the
same as those of
$\ltwo C_{\ast}(\gS)$, i.e., $\ltwo H_{\ast}(\gS_{cc})\cong \ltwo
H_{\ast}(\gS)$, $\ltwo H^{\ast}(\gS_{cc})\cong \ltwo H^{\ast}(\gS)$ and
$\ltwo \ch^{\ast}(\gS_{cc})\cong \ltwo \ch^{\ast}(\gS)$.  (The point is
that the simplicial structure on $\gS$ is a subdivision of $\gS_{cc}$.)

The chain complex $(\ltwo C_{\ast}(\gS_{cc}),\partial^\bq)$ looks like 
this:
\begin{equation*}
\ltwo\ \longleftarrow \bigoplus_{s\in S}H_s\ \longleftarrow
\bigoplus_{T\in \cs^{(2)}}H_T \longleftarrow\cdots
\end{equation*}
(We shall describe the boundary maps explicitly in 
Lemma~\ref{l:chainruin} in
the next section.)

\rk{$\ltwo$-Betti numbers and the $\ltwo$-Euler characteristic}  Define
    \begin{equation*}
    c^i_\bq(\cu):=\dim_{\cn} \ltwo C^i(\cu),
    \end{equation*}
where $\dim_{\cn}$ denotes the von Neumann dimension defined by (\ref{e:defdim}).
For any cell $\gs\subset Z$, its stabilizer is the special subgroup 
$W_{S(\gs})$, where as before $S(\gs)=\{s\in S\mid \gs\subseteq Z_s\}$.  
So, 
the summand of $\ltwo C^i(\cu)$ corresponding to the orbit of an 
$i$-cell $\gs$ is isomorphic to $A_{S(\gs)}$.  Its dimension is
$1/W_{S(\gs)}$.  Hence,
\begin{equation}\label{e:ciq}
c^i_\bq(\cu)=\sum_{\gs\in Z^{(i)}}\frac{1}{W_{S(\gs)}}.
\end{equation}

The \emph{$i^{th}\ \ltwo$-Betti number}
of $\cu$ is defined by
    \begin{equation}\label{e:betti}
    b^i_\bq(\cu):=\dim_{\cn} \ltwo\ch^i(\cu).
    \end{equation}
A standard argument (cf., \cite[Theorem 3.6.1, p. 205]{eckmann}) gives:
    \begin{equation}\label{e:atiyah}
    \sum (-1)^ib^i_\bq(\cu)=\sum (-1)^ic^i_\bq(\cu).
    \end{equation}
(This is a version of Atiyah's Formula.)
We denote either side of (\ref{e:atiyah}) by $\chi_\bq(\cu)$ and call it the
\emph{$\ltwo$-Euler characteristic} of $\cu$.  

\begin{proposition}[Rationality of Euler characteristics]\label{p:rationaleuler}  
$\chi_\bq(\cu)=f(\bq)/g(\bq)$ where $f$ and $g$ are polynomials in $\bq$ 
with integral coefficients.
\end{proposition}

\begin{proof}
For each $T\in \cs$, we have the subcomplex $Z_T$ (resp. $\partial Z_T$) 
defined as the union of those cells $\gs$ such that $T\subseteq S(\gs)$ 
(resp. $T\subset S(\gs)$).  By (\ref{e:ciq}) and (\ref{e:atiyah}),  
\[
\chi_\bq(\cu)=\sum_{T\in\cs}\frac{\chi(Z_T)-\chi (\partial Z_T)}{W_T(\bq)}.
\]
\end{proof}

\begin{proposition}[{\cite[Corollary 3.4]{dym}}]\label{p:euler}
\[
\chi _\bq(\gS)=\frac{1}{W(\bq)}.
\]
\end{proposition}

\begin{proof}
We use the cellulation of $\gS$ by Coxeter cells.  If $T\in \cs$, then
\[
\dim_{\cn} L^2(W\angt,\muq)=\dim_{\cn} H_T=\frac{1}{W_T(\bq^{-1})}.
\]
Hence,
\begin{align}
c^i_\bq(\gS_{cc})&=\sum _{T\in \cs^{(i)}}
\frac{1}{W_T(\bq^{-1})}\text{ and}\notag\\
\chi_\bq(\gS)&=\sum _{T\in \cs}
\frac{\geps (T)}{W_T(\bq^{-1})}=\frac{1}{W(\bq)},\notag
\end{align}
where the last equality is by Lemma~\ref{l:growth}(iv).
\end{proof}

\begin{Remark}
The relationship between Euler characteristics (of groups acting on
buildings) and growth series of Coxeter groups was first
pointed out by Serre \cite{serre} (in the case where
with fundamental chamber $K$ is a simplex).  Serre showed that the
``Euler-Poincar\'e'' measure on the automorphism group of the building is
(suitably normalized) Haar measure multiplied by $1/W(\bq)$, where, as in
Section~\ref{s:buildings}, $\bq$
is the ``thickness vector'' of the building.
\end{Remark}

\rk{Cohomology in dimension $0$}  
The vertex set of $\gS$, with its cellulation by Coxeter cells, can be identified 
with $W$.  So, $\ltwo C^0(\gS_{cc})\cong \ltwo$.  A $0$-cochain is a cocycle 
if and only if it is the constant function on $W$.  If $c$ denotes the 
constant, then its norm, with respect to the inner product $\langle\ ,\ 
\rangle _\bq$ is $|c|\sum q_w$ and this is $<\infty$ 
if  and only if $\bq\in \car$ or $c=0$.  This 
proves the following result of \cite{dym}.

\begin{proposition}[\cite{dym}]\label{p:HO}
$\ltwo H^0(\gS)$ is nonzero if and only if $\bq\in\car$.  Moreover, when 
$\bq\in\car$, $b^0_\bq(\gS)=1/W(\bq)$.
\end{proposition}

\begin{remark}\label{r:H0}
It is easy to see that the space $\cu$ ($=\cu (W,Z)$) is connected if and only if $Z$ is connected and $Z_s\neq \emptyset$ for each $s\in S$. (This also follows from \cite[Theorem A]{d87} or \cite[Theorem 10.1]{d83}.)  Suppose 
these conditions hold.  An argument similar to the one in the previous paragraph then shows that 
$\ltwo H^0(\cu)$ is nonzero if and only if $\bq\in\car$ and when this is the 
case, $b^0_\bq(\cu)=1/W(\bq)$.
\end{remark}

\rk{The continuity of Betti numbers}  
\begin{theorem}\label{t:cont}
Suppose $(W,S)$ is a Coxeter system and that $\cu=\cu(W,Z)$ is as above.  
Then for each integer $i$, the function $\bq\to b^i_\bq(\cu)$ is 
continuous.
\end{theorem}

For the proof  we will need  the next two lemmas.

\begin{lemma}\label{l:hilbert}
Let $Y$ be a Hilbert space, $X$ a closed subspace of $Y$, $P_X$ the 
orthogonal projection onto $X$ and
 $y\in Y$ a unit vector.  Set 
\[
A(y):=\inf{\{\|x\|\mid x\in X,\, \langle 
x,y\rangle=1\}}.
\]
Then
\(
\langle P_X(y),y\rangle=A(y)^{-2}.
\)
(By convention,  $(+\infty)^{-2}=0$.)
\end{lemma}

\begin{proof}
Put $a:=\langle P_X(y),y\rangle$.  Since $\langle 
P_X(y),y\rangle=\|P_X(y)\|^2$, we see that $a\ge0$ with equality if and only if 
$X\perp y$.  Suppose first that $a=0$.  Then the left hand side of the 
formula in the lemma is $0$.  Since $X\perp y$,  
$\{\|x\|\mid x\in X,\, \langle x,y\rangle=1\}=\emptyset$, so 
$A(y)=+\infty$ and hence, the
right hand side is also $0$.

Suppose $a>0$.  Every $x\in X$ can be written as $bP_X(y)+x'$, where 
$x'\perp P_X(y)$.  
Then $\langle x,y\rangle=
b\langle P_X(y),y\rangle+\langle x',y\rangle=ba$. (Notice that for 
$x'\in X$,  $x'\perp P_X(y)$ if and only if $x'\perp y$.)  So, $\langle 
x,y\rangle=1$ implies
$b=\frac{1}{a}$.  Therefore,
\begin{align*}
A(y)&=\inf{\{\|x\|\mid x\in X,\, \langle x,y\rangle=1\}}\\
&=\inf{\{\|\frac{1}{a}P_X(y)+x'\|\mid x'\in X,\, x'\perp P_X(y)\}}\\
&=\|\frac{1}{a}P_X(y)\|=\frac{1}{a}\|P_X(y)\|=\frac{\sqrt{a}}{a}=\frac{1}{\sqrt{a}}.  
\end{align*}
So, $A(y)^{-2}=a$.
\end{proof}

Given $I$-tuples $\bq$ and $\bq'$ of real numbers, write $\bq\le\bq'$ if 
$q_i\le q'_i$ for each $i\in I$.  
\begin{definition}\label{d:increasing}
A function $f:\bR^I\to \bR$ is \emph{increasing} (resp. \emph{decreasing})
if $\bq\le \bq'$ implies $f(\bq)\le f(\bq')$ (resp. $f(\bq)\ge f(\bq')$).
\end{definition}

Let $|\ |$ denote the ``maximum norm'' on $\bR^I$ defined by $|\bq|:= 
\max\{|q_i|\}$.

\begin{definition}\label{d:lrcont}
A function $f:\bR^I\to \bR$ is \emph{left continuous} at $\bq_0$ if for 
any positive number $\geps$, there is a positive number $\gd$ so that 
if $\bq\le \bq_0$ and $|\bq_0-\bq|<\gd$, then 
$|f(\bq_0)-f(\bq)|<\geps $.  \emph{Right continuity} is similarly defined.
\end{definition}

\begin{lemma}\label{l:function}
If a decreasing function $f:\bR^I\to \bR$ is both left and right 
continuous, then it is continuous (and similarly, if $f$ is increasing).
\end{lemma}

\begin{proof}
Given a point $\bq_0\in \bR^I$ and a number $\geps>0$, choose $\gd$ small 
enough to work in the definitions of both left and right continuous at 
$\bq_0$.  Let $\bdelta$ denote the $I$-tuple with each component equal to 
$\gd$.  Assuming $f$ is decreasing, 
for any $\bq$ in a $\gd$-neighborhood of $\bq_0$, we have
\[
f(\bq_0)-\geps\le f(\bq_0+\bdelta)\le f(\bq)\le f(\bq_0-\bdelta)\le 
f(\bq_0)+\geps.
\]
\end{proof}

\begin{proof}[Proof of Theorem~\ref{t:cont}]
We have spaces of cochains, cocycles and coboundaries:
\begin{align*}
C^i_\bq &:=\ltwo C^i(\cu), & Z^i_\bq &:=\ltwo Z^i(\cu), & 
B^i_\bq &:=\ltwo B^i(\cu),\\
\intertext{as well as, spaces of chains, cycles and boundaries:}
C_i^\bq &:=C^i_\bq, & Z_i^\bq &:=\ltwo Z_i(\cu), & 
B_i^\bq &:=\ltwo B_i(\cu)\\
\intertext{($Z^i_\bq$ and $B^i_\bq$ are defined using the coboundary map $\gd$, while 
$Z_i^\bq$ and $B_i^\bq$ are defined using its adjoint $\partial^\bq$.)  We 
also have their von Neumann dimensions:}
c^i_\bq&:=\dim C^i_\bq,  &z^i_\bq&:=\dim Z^i_\bq, 
&a^i_\bq&:=\dim 
B^i_\bq,\\
c_i^\bq&:=c^i_\bq,  &z_i^\bq&:=\dim Z_i^\bq,  &a_i^\bq&:=\dim 
B_i^\bq,
\end{align*}
where, 
to simplify notation, we are  writing $\dim (\ )$ instead of $\dim_{\cn}(\ 
)$.  

We note that by formula (\ref{e:ciq}), $\bq \to c^i_\bq$ is a 
continuous decreasing function (since each $W_{S(\gs)}(\bq)$ is a 
polynomial with nonnegative coefficients).

\begin{step}
The function $\bq\to z^i_\bq$ is left continuous and decreasing.
\end{step}

\begin{proof}[Proof of Claim 1]
In (\ref{e:Phi}) we defined an isometric embedding $\Phi$ of 
$C^i_\bq$ into
$\bigoplus_{\gs\in Z^{(i)}}\ltwo$.  Let $e^\gs_1$ be the element of  
$\bigoplus_{\gs\in 
Z^{(i)}}\ltwo$ with $\gs$-component 
equal to $e_1$ and all other components equal $0$. Then 
\[
z^i_\bq=\sum_{\gs\in Z^{(i)}} \langle 
P_{\Phi(Z^i_\bq)}(e^\gs_1),e^\gs_1\rangle_\bq. 
\] 
Since 
$\Phi(Z^i_\bq)\subseteq\Phi(C^i_\bq)$, we have 
\[
\langle P_{\Phi(Z^i_\bq)}(e^\gs_1),e^\gs_1\rangle_\bq=
\langle P_{\Phi(Z^i_\bq)}P_{\Phi(C^i_\bq)}(e^\gs_1), 
P_{\Phi(C^i_\bq)}(e^\gs_1)\rangle_\bq.
\]
 All components of the vector $P_{\Phi(C^i_\bq)}(e^\gs_1)$ 
are 0, except the $\gs$-component,
which is equal to $P_{A_{S(\gs)}}(e_1)=a_{S(\gs)}=
\Phi((W_{S(\gs)}(\bq))^{-1/2}\delta_\gs)$,
where $\delta_\gs\in C^i_\bq$ is the function which is 1 on $\gs$ 
and 0 on all other cells.
Thus, 
\[
z^i_\bq=
\sum_{\gs\in Z^{(i)}} \frac{1}{W_{S(\gs)}(\bq)}
\langle P_{Z^i_\bq}(\delta_\gs),\delta_\gs\rangle_\bq.
\]
Since $1/W_{S(\gs)}(\bq)$ is continuous, 
we need to concentrate on $\langle 
P_{Z^i_\bq}(\delta_\gs),\delta_\gs\rangle_\bq$.  Set
 \(
z^\gs(\bq):=\inf{\{\|u\|_\bq\mid u\in Z^i_\bq,\, \langle u, 
\delta_\gs\rangle_\bq=1\}}. 
\)
By Lemma~\ref{l:hilbert}, it suffices to prove that each 
of the functions $z^\gs(\bq)$ is left continuous and increasing.  
Notice that if $\bq\le\bq'$, then $Z^i_\bq\supseteq Z^i_{\bq'}$. 
Moreover, if $u\in Z^i_{\bq'}$, then $\|u\|_\bq\le\|u\|_{\bq'}$, while
$\langle u,\delta_\gs\rangle_\bq=\langle u,\delta_\gs\rangle_{\bq'}$ 
(because both are equal to $u(\gs)$). 
It follows that $z^\gs$ is an increasing function.

Now suppose that $(\bq_n)$ is a sequence in $\bR^I$ converging 
to $\bq$ from below (that is, each $\bq_n\le \bq$) and that $u_n\in 
Z^i_{\bq_n}$ is a sequence such that, 
$\langle u_n,\delta_\gs\rangle_{\bq_n}=1$ (i.e., $u_n(\gs)=1$).  Further 
suppose that  
$\lim \|u_n\|_{\bq_n}=\xi$. We will show that 
$\xi\ge z^\gs(\bq)$. (This implies that $z^\gs$ is left continuous at 
$\bq$.)
Write $\cu$ as $\cu=\bigcup_{k=1}^\infty K_k$, where the $K_k$ are finite 
subcomplexes.
Assume that $\xi<+\infty$ (otherwise there is nothing to prove).
This implies that for every  $k$, the restrictions 
$u_{n}|_{K_k}$  are uniformly bounded.
Hence, by a diagonal argument, one can choose a subsequence $(u_{m})$ such 
that  the
$u_{m}|_{K_k}$ converge pointwise for each $k$. Let $u$ be the pointwise 
limit of $u_{m}$.
Then $u$ is a cocycle and $u(\gs)=1$ (because all  $u_{m}$ satisfy 
these conditions). 
Also, for each $k$ we have 
$\|u|_{K_k}\|_\bq=\lim{\|u_{m}|_{K_k}\|_{\bq_m}}\le\xi$.  
Therefore, $\|u\|_\bq\le\xi$. So, $z^\gs(\bq)\le\|u\|_\bq\le\xi$.
\end{proof}

\begin{step}
$\bq\to a^i_\bq$ is left continuous and so, $\bq\to z_i^\bq$ is 
left continuous.
\end{step}

\begin{proof}[Proof of Claim 2]
Since $c^{i-1}_\bq$ is continuous and since 
$a^i_\bq=c^{i-1}_\bq-z^{i-1}_\bq$, Claim 1 implies that $a^i_\bq$ is left 
continuous.  We have the Hodge decomposition:  $C^i_\bq=Z_i^\bq\osum 
B^i_\bq$.  So, $z_i^\bq=c^i_\bq-a^i_\bq$, which is left continuous.
\end{proof}

\begin{step}
$\bq\to z_i^\bq$ is right continuous and decreasing.
\end{step}

\begin{proof}[Proof of Claim 3]
This is a version of Claim 1 using cycles instead of cocycles.  
Basically,  
the argument in Claim 1 works provided we use the usual
boundary map $\partial$ instead of $\partial^\bq$.  To transfer this back 
into information about $Z_i^\bq$, we need to use the isometry $\theta$ 
from 
the proof of Lemma~\ref{l:qtoq-1}.  Set
\[
\wh{Z}^\bq_i:=\Ker (\partial_i:C_i^\bq\to C_{i-1}^\bq).
\]
As before, $C^\bq_i$,  $Z^\bq_i$ and  $\wh{Z}^\bq_i$
can be embedded into
$\bigoplus_{\gs\in Z^{(i)}}\ltwo (W)$, and $z^\bq_i=
\sum_{\gs\in Z^{(i)}} \langle 
P_{\Phi(Z^\bq_i)}(e^\gs_1),e^\gs_1\rangle_\bq$. 
Consider the isometry $\theta:\bigoplus \ltwo (W)\to
\bigoplus L^2_{\bq^{-1}}(W)$, given by $\theta(f)(w)=\bq_wf(w)$ on each 
component. By Lemma~\ref{l:qtoq-1}, $\theta$ restricts to a map 
$C^\bq_*\to C^{\bq^{-1}}_*$,
which intertwines $\partial^\bq$ and $\partial$.
Therefore, $\theta (Z^\bq_i)=\wh{Z}_i^{\bq^{-1}}$. Also,
$\theta(e^\gs_1)=e^\gs_1$, so 
$\theta(P_{\Phi(Z^\bq_i)}(e^\gs_1))=P_{\Phi(\wh{Z}_i^{\bq^{-1}})}(e^\gs_1)$, 
and $\langle P_{\Phi(Z^\bq_i)}(e^\gs_1),e^\gs_1\rangle_\bq
=\langle P_{\Phi(\wh{Z}_i^{\bq^{-1}})}(e^\gs_1), 
e^\gs_1\rangle_{\bq^{-1}}$.
(Note that the map $\Phi$ depends on $\bq$; thus, the maps on the left
hand sides  correspond to $\bq$, while those on the right hand sides 
correspond to $\bq^{-1}$.)
Now, the argument from Claim~1 can be repeated. We get that
$\langle P_{\Phi(\wh{Z}_i^{\bq^{-1}})}(e^\gs_1), 
e^\gs_1\rangle_{\bq^{-1}}$
is left continuous and increasing in $\bq^{-1}$. This implies that
$\bq\to \langle P_{\Phi(Z^\bq_i)}(e^\gs_1),e^\gs_1\rangle_\bq$ 
is right continuous and decreasing.
\end{proof}

\begin{step}
$\bq\to a_i^\bq$ is right continuous and so, $\bq\to z^i_\bq$ is 
right continuous.
\end{step}

\begin{proof}[Proof of Claim 4]
This follows from Claim~3 in the same way Claim~2 followed from Claim~1.
\end{proof}

\begin{step}
$z^i_\bq$, $z_i^\bq$, $a^i_\bq$ and $a_i^\bq$ are continuous in $\bq$
\end{step}

\begin{proof}[Proof of Claim 5]
The functions $z^i_\bq$ and $z_i^\bq$ are decreasing and left and right 
continuous; hence, by Lemma~\ref{l:function}, continuous.  Since 
$c^i_\bq$ is continuous, $a^i_\bq$ and $a_i^\bq$ are also continuous.
\end{proof}
To finish the proof of Theorem~\ref{t:cont} simply note that
 $b^i_\bq(\cu)=z^i_\bq-a^i_\bq$, which,  by Claim~5, is continuous.
\end{proof}

In view of Proposition~\ref{p:rationaleuler} and Atiyah's Conjecture (cf. 
\cite[Section 3.10]{eckmann} or \cite[Chapter 10]{luck}), it is natural to ask the following.

\begin{Question}
Is $\bq\to b^i_\bq(\cu)$ a piecewise rational function?  
\end{Question}

\section{Weighted $L^2$-homology of ruins}\label{s:ruins}\nopagebreak[4]
\rkn{Cosheaves}  Suppose $\gL$ is a simplicial complex
with vertex set $V$ and
that $\cs (\gL)$ is its face poset (including the empty face).
We  regard the poset $\cs (\gL)$ as a category in the usual way:  if
$\tau$ is a face of $\gs$, then there is a unique morphism $\gi
^{\tau}_{\gs}$ from $\tau$ to $\gs$ (which we can think of as being the
inclusion of vertex sets).

A \emph{cosheaf} on $\gL$ with values in a category $\cac$ is a 
contravariant
functor $F$ from $\cs (\gL)$ to $\cac$.  In the case of interest to
us, $\cac$ will be the category of Hilbert $\cn$-modules.

Now suppose that the simplicial complex $\gL$ is ordered (in other words,
suppose that its vertex set is totally ordered).  Then
for any $n\geq 0$, the vertices of an $n$-simplex form an ordered set
isomorphic to $\{0,1,\dots,n\}$ with its usual order.  For 
$0\leq i\leq n$ and any $n$-simplex $\gs$, the $i^{\mathrm{th}}$ face of
$\gs$ is defined to be the $(n-1)$-simplex spanned by all vertices
of $\gs$ except the $i^{\mathrm{th}}$.  If one writes
$\partial_i$ for $\iota^{\tau}_\gs$, where $\tau$ is the $i^{\mathrm{th}}$ 
face of $\gs$, then the relations between the morphisms become the 
familiar ``simplicial identities'' as in \cite[8.1]{weibel}.

A cosheaf $F$ of abelian groups on an ordered simplicial complex $\gL$ 
gives
rise to a chain complex $C_\ast(\cs (\gL);F)$ defined as follows: $C_n=0$ 
for
$n<0$, and for $n\geq 0$,
    \begin{equation*}
    C_n= \bigoplus_{\gs\in \cs^{(n)}(\gL)} F(\gs),
    \end{equation*}
where the indexing set is the set of $(n-1)$-simplices of $\gL$.  (The
indices on $C_{\ast}$ have been shifted up by one from the 
conventions in  \cite{dl}.)  Under the natural isomorphism
\[
\Hom(C_n,C_{n-1}) \cong \bigoplus_{\gs\in \cs^{(n)}(\gL),\,\,
\tau\in \cs^{(n-1)}(\gL)} \Hom(F(\gs),F(\tau)),
\]
the boundary map $\partial:C_n\rightarrow C_{n-1}$ corresponds to the
matrix $(\partial_{\gs\tau})$, where $\partial_{\gs\tau}=0$ unless $\tau$ is
a face of $\gs$, and is equal to $(-1)^iF(\iota_\gs^\tau)$ if
$\tau$ is the $i$th face of $\gs$.

\rk{Ruined chain complexes} We return to the situation where $(W,S)$
is a Coxeter system, $L$ is its nerve and  $\cs$ is the poset of
spherical subsets of $S$.
Let $T\in \cs$ and let $\Lk (T,L)$ denote
the link of $T$ in $L$.  (If $T=\emptyset$,  $\Lk (\emptyset,
L):=L$.)  We note that the face poset of $\Lk(T,L)$ is isomorphic to
$\cs_{\geq T}$.

Define a cosheaf $H_T$ of Hilbert $\cn$-modules on $\Lk(T,L)$  as the 
contravariant functor on $\cs_{\geq T}$ defined 
on objects by $U\to H_U$ where $H_U$ is defined by
(\ref{e:ATHT}).  If $U\leq V \in \cs_{\geq T}$, then the morphism 
$H(\gi^U_V):H_V\to H_U$  
is the natural inclusion.  Define the \emph{ruined chain complex}
$\ltwo C_{\ast}(H_T)$ by
    \begin{equation*}
    \ltwo C_{\ast}(H_T):=C_{\ast}(\cs_{\geq T};H).
    \end{equation*}
It looks like this:
    \begin{equation}\label{e:lookslike1}
    0\ \longleftarrow
    H_T\ \longleftarrow \bigoplus_{(T\cup
    \{s\})\in (\cs_{\geq T})^{(k+1)}}H_{T\cup\{s\}}\
    \longleftarrow\cdots,
    \end{equation}
where $k=\Card(T)$.
Similarly, by using the family of Hilbert $\cn$-modules
$(A_U)_{U\in \cs_{\geq T}}$,  we get a cosheaf $A$ on $\Lk (T,L)$ and a
chain complex $\ltwo C_{\ast}(A_T)$.

We denote the homology of $\ltwo C_{\ast}(H_T)$ and $\ltwo C_{\ast}(A_T)$
by $\ltwo H_{\ast}(H_T)$ and $\ltwo H_{\ast}(A_T)$, respectively.

\rk{The relationship between ruins and ruined chain complexes}
Recall that for any $U\subseteq S$ and $T\in \cs (U)$, 
$\gO (U,T)$ is the subcomplex of $\gS_{cc}$ consisting of all
closed Coxeter cells of type $T'$, with $T'\in \cs(U)_{\geq T}$.  

To simplify notation, the chain complex $\ltwo
C_*(\gO(U,T),\partial (\gO(U,T))$ will be denoted $\ltwo
C_*(\gO(U,T),\partial)$
and similarly for its homology.

Since the cell structure always will be  given by Coxeter cells,
we will omit the subscript $cc$ from our
notation.  We say  a Coxeter cell is \emph{type $T$}, $T\in \cs$, as
a shorthand for type $W_T$.

It follows from (\ref{e:sigmaU}) and the fact that $\cu(W_U,K(U))$
deformation retracts onto $\gS_U$ that
the $\cn [W]$-modules $\ltwo C_{\ast}(\gS
(U))$ and $\ltwo\ch_{\ast}(\gS(U))$ are induced from the
$\cn [W_U]$-modules $\ltwo C_{\ast}(\gS_{W_U})$ and
$\ltwo\ch_{\ast}(\gS_{W_U})$, respectively.  So, we can calculate von 
Neumann
dimensions over $\cn [W]$ by calculating with respect to
$\cn [W_U]$.

\begin{lemma}\label{l:chainruin}
\begin{enumeratei}
\item
There is a isomorphism of chain complexes of $\cn$-modules:
\[
\Psi':\ltwo C_\ast(\gS_{cc})\to \ltwo C_\ast(H_{\emptyset}),
\]
where $\ltwo C_\ast(H_{\emptyset})$
is the ruined chain complex associated to the cosheaf $H_{\emptyset}$ on 
$L$.
\item
Suppose
$T\in\cs^{(k)}$.  Then $\Psi'$ induces an
isomorphism of chain complexes of $\cn$-modules:
\[
\ltwo C_{*}(H_T)\mapright{\cong}\ltwo C_{*+k}(\gO (S,T),\partial).
\]
In particular, $\ltwo C_{m}(\gO (S,T),\partial)=0$ for $m<k$.
\end{enumeratei}
\end{lemma}

\begin{proof}
(i)  For each $T\in\cs$ modify the isometry $\psi_T$ of
(\ref{e:psi}) to another
Hilbert $\cn$-module isomorphism,
$\psi '_T:L^2(W\langle T\rangle,\muq)\to H_T$ as follows:
    \begin{equation}\label{e:psi'}
    \psi '_T(f):=\sqrt{W_T(q^{-1})}\psi_T(f)=
    W_T(q^{-1})\Big(\sum _{u\in X_T}f(u\angt)e_u \Big)h_T.
    \end{equation}
$\Psi '$ is defined to be the direct sum of the $\psi '_T$.
Suppose $U\in\cs_{>\emptyset}$ and $T\subset U$ is obtained by deleting
one element  from $U$.  Statement (i) follows immediately from the
next claim.

\begin{claim}
The following diagram commutes:
\[
    \begin{CD}
    L^2(W\langle U\rangle,\muq)@>\psi '_U>>H_U\\
    @V\partial_T^\bq VV@VViV\\
    L^2(W\langle T\rangle,\muq)@>\psi '_T>>H_T
    \end{CD}
\]
where $\partial_T^\bq $ denotes the $L^2(W\langle T\rangle,\muq)$-component
of $\partial ^\bq$ and $i$ is the natural inclusion.
\end{claim}

\begin{proof}[Proof of Claim]  Using (\ref{e:psi'}) and (\ref{e:partialT}), we get
\begin{align*}
\psi'_T(\partial_T^\bq f)&=W_T(\bq^{-1})\Big(\sum_{w\in 
X_T}(\partial_T^\bq f )
(w\angt) e_w\Big)h_T\\
&=W_T(\bq^{-1})\Big(\sum_{u\in X_U}\sum_{v\in W_U\cap X_T}\geps_vq_v^{-1}
f(u\langle U\rangle) e_ue_v\Big)h_T\\
&=W_T(\bq^{-1})\Big(\sum_{u\in X_U}f(u\langle
U\rangle)e_u\Big)\Big(\sum_{v\in W_U\cap X_T}\geps_vq_v^{-1}e_v\Big)h_T\\
&=W_U(\bq^{-1})\Big(\sum_{u\in X_U}f(u\langle U\rangle)e_u\Big)h_U\\
&=i(\psi'_U(f)),
\end{align*}
where the next to last equality is from the following formula for $h_U$, 
valid whenever $T\subseteq U$ and $\bq\in \car_U^{-1}$.
\[
h_U=\Big(\sum_{v\in W_U\cap X_T}\geps_vq_v^{-1}e_v\Big)h_T.
\]
(This formula holds since $W_U\cap X_T$ is a set of coset 
representatives for $W_U/W_T$ and since for any $v\in W_U\cap X_T$ 
and $w\in W_T$, we have $e_ve_w=e_{vw}$ and $q_vq_w=q_{vw}$.) 
\end{proof}

(ii)  Part (ii) of the lemma essentially  follows from part (i).  
Write $\gO$ for $\gO(S,T)$. The point is that the cells of the 
$W\langle T'\rangle$, $T'\in (\cs_{\ge T})^{(i+1)}$, are a basis for $\ltwo
C^i(\gO,\partial \gO)$.  Hence,
\[
   \ltwo C_i(\gO,\partial \gO)=\bigoplus_{T'\in (\cs_{\geq
   T})^{(i+1)}}L^2(W\langle T'\rangle,\muq)\cong\bigoplus H_{T'}
\]
and (i) shows that the $\partial^\bq$ maps are induced by the inclusions
$H_{T''}\hookrightarrow H_{T'}$, with $T'\subset T''$.
\end{proof}

\begin{Remark}\hspace{-1ex plus 1ex}
 The cochain complex $\ltwo C^i(\gO(S,T),\partial)$
is obtained by dualizing (\ref{e:lookslike1}):
    \begin{equation*}
    0\ \longrightarrow
    H_T\ \longrightarrow \bigoplus_{(T\cup
    \{s\})\in (\cs_{\geq T})^{(k+1)}}H_{T\cup\{s\}}\
    \longrightarrow\cdots,
    \end{equation*}
where the coboundary maps are induced by the orthogonal projections
$H_{T'}\to H_{T''}$, with $T'\subset T''$, and  $k=\Card(T)$.
\end{Remark}

The main result of this section as well as the results of Sections 9 
through 12 ultimately are based on the following key theorem from 
\cite{dym}.  

\begin{theorem}[{\cite[Theorem 10.3]{dym}}]\label{t:dym}
If $\bq\in\car$, then $\ltwo H_*(\gS)$ 
is concentrated in dimension $0$.
\end{theorem}

While the proof of this in \cite{dym} is straightforward, some 
technical estimates are involved.  In outline the argument goes as follows.
\begin{enumeratea}
\item
Using the $\text{CAT}(0)$-metric of \cite{moussong} it is proved, in 
\cite[Theorem 9.1]{dym}, that 
there is a chain contraction $H:C_*(\gS)\to C_{*+1}(\gS)$ and constants $C$ 
and $R$ such that for any simplex $\gs\subset \gS$, (i) the $L^{\infty}$-norm 
of $H(\gs)$ is $<C$ and  (ii) $H(\gs)$ is supported in an 
$R$-neighborhood of the geodesic connecting the central vertex of $K$ 
with $\gs$.
\item
It follows (\cite[Theorem 10.1]{dym}) that for $\bq\in\car^{-1}$, $H$ 
extends to a bounded linear map $H:\ltwo C_*(\gS)\to\ltwo C_{*+1}(\gS)$.
(Actually, in \cite{dym}, this is only proved for a single parameter, but the proof goes through without change in the case of a multiparameter 
$\bq$.)  Hence, for $\bq\in \car^{-1}$, $H$ is a chain contraction 
of $\ltwo C_*(\gs)$ with respect to the usual boundary map $\partial$.
\item
Finally, one uses the isometry $\theta$ of Lemma~\ref{l:qtoq-1} to 
transport $H$ to a chain contraction of $(\ltwo C_*(\gs),\partial^\bq)$ 
for $\bq\in\car$.
\end{enumeratea}

The main result of this section is the following generalization of Theorem~\ref{t:dym}.

\begin{theorem}\label{t:dimk}
Suppose $T\in \cs^{(k)}$.  If $\bq\in\car$, then $\ltwo H_{\ast}(\gO
(S,T),\partial)$  is concentrated in dimension $k$.  If $\bq\in 
\partial\car$,  the same holds for $\ltwo \ch_{\ast}(\gO(S,T),\partial)$. 
\end{theorem}

Note that the third sentence of the theorem follows from the second one 
and the continuity of the $b^i_\bq$ (Theorem~\ref{t:cont}).

In the special case $T=\emptyset$, we have $\gO (S,T)=\gS$ and so, 
Theorem~\ref{t:dimk} is Theorem~\ref{t:dym}.  We shall
use Theorem~\ref{t:dym} as the first step in an inductive proof.

Before beginning the proof,
note that we have an  excision isomorphism: 
	\begin{equation}\label{e:excision1}
	\ltwo C_{\ast}(\gO (U,T),\partial)\cong\ltwo C_{\ast}(\gS
	(U),\wh{\gO}(U,T)).
	\end{equation}
Also, for any $s\in T$ and $T':=T-s$, we have an excision isomorphism:
	\begin{equation}\label{e:excision2}
	\ltwo C_{\ast}(\gS (U-s),\wh{\gO}(U-s,T'))\cong
	\ltwo C_{\ast}(\wh{\gO}(U,T),\wh{\gO}(U,T')).
	\end{equation}

\begin{proof}[Proof of Theorem~\ref{t:dimk}]
Suppose $U\subseteq S$ and $T\in \cs^{(k)}(U)$.  We shall prove, by
induction on $k$ ($=\Card (T)$), that $\ltwo H_{\ast}(\gO (U,T),\partial)$
is concentrated in dimension $k$.  When $k=0$ this holds by 
Theorem~\ref{t:dym} (and the fact that  $\ltwo
C_{\ast}(\gS (U))$ is induced from $\ltwo C_{\ast}(\gS_{W_U})$).
Assume by induction, that our assertion holds for $k-1$, with $k-1\geq 0$.
By (\ref{e:excision1}), the assertion
is equivalent to showing that $\ltwo H_{\ast}(\gS (U),\wh{\gO}(U,T))$ is
concentrated in dimension $k$.  Choose an element $s\in T$ and set
$T':=T-s$, $\wh{\gO}:=\wh{\gO}(U,T)$, $\wh{\gO}':=\wh{\gO}(U,T')$.
Consider the long exact sequence of the
triple $(\gS (U),\wh{\gO},\wh{\gO}')$:
\[
\ltwo H^{\ast}(\gS (U),\wh{\gO}')\to
\ltwo H^{\ast}(\gS (U),\wh{\gO})\to
  \ltwo H^{\ast -1}(\wh{\gO},\wh{\gO}')
\]
By (\ref{e:excision2}), the right hand term excises to the homology of the
$(U-s,T')$-ruin, while the middle term is that of the $(U,T)$-ruin and the
left hand term is that of the $(U,T')$-ruin.  By induction, the left hand
and right hand terms are concentrated in dimension $k-1$.  So, the middle
term can only be nonzero in dimensions $k-1$ and $k$.  On the other hand,
by Lemma~\ref{l:chainruin}(ii), the middle term vanishes in
dimensions $<k$.
\end{proof}

Combining this theorem with Lemma~\ref{l:chainruin}, we get the following.

\begin{corollary}\label{t:acycliccomplex}
For any $\bq\in\car$ and any spherical subset $T$, $\ltwo H_{\ast}(H_T)$ is
concentrated in dimension 0.  Therefore, for any $\bq\in \ol{\car}$, the 
reduced homology $\ltwo \ch_{\ast}(H_T)$ is also concentrated in dimension 
$0$.
\end{corollary}

The meaning of this corollary is that, for $\bq\in \ol{\car}$, the family
of subspaces
$(H_T)_{T\in \cs}$ is ``in general position'' in $\ltwo$.

\section{The Decomposition Theorem}\label{s:decomp}

\begin{lemma}[Compare Lemma 1 in \cite{s}]\label{l:sol1}
Suppose we are given subsets
$U,V$ of $S$ and an $I$-tuple $\bq \in \car_V\cap \car_U^{-1}$
(so that $h_U$ and $a_V$ are both defined).
If $V\cap U\neq \emptyset$, then
$h_Ua_V=0$.
\end{lemma}

\begin{proof}
Let $s\in V\cap U$. Then $h_Ua_V = h_Uh_s a_s a_V =0$.
\end{proof}

We define some more subspaces of $\ltwo$:
 \begin{align*}
D_V&:=A_{S-V}\cap\Big(\sum _{U\subset V}A_{S-U}\Big)^\perp,\\
G_V&:=H_V\cap\Big(\sum _{U\supset V}H_U\Big)^\perp.
\end{align*}

\begin{lemma}\label{l:sumGT}
$$ \ol{\sum_{U\supseteq V} G_U} = H_V, $$
$$ \ol{\sum_{V\subseteq U} D_V} = A_{S-U}. $$
\end{lemma}

\begin{proof}
By definition of $G_V$, we have
$$H_V= G_V + \ol{\sum_{U\supset V} H_U},$$
and the first formula follows by induction on the size of $S-V$.
Similarly, 
$$A_{S-V}= D_V + \ol{\sum_{U\subset V} A_{S-U}},$$
and the second formula follows by induction on the size of $V$.
\end{proof}

\begin{lemma}\label{l:GUaS-V}
Suppose $\bq\in \car$ and $U\nsubseteq V$.  Then
$$G_Ua_{S-V}=0.$$
\end{lemma}

\begin{proof}
Since $G_T\subseteq H_T$ the assertion follows from
Lemma~\ref{l:sol1}.
\end{proof}

If $\bq \in \car$, then $H_V=0$ for all nonspherical
$V$ (because $V$ is spherical whenever $\car_V \cap\car_V^{-1}\neq\emptyset$).
So, for $V\notin \cs$, $G_V=0$, and for
$V=T\in \cs$, $G_T$ is the orthogonal complement of the image of
$\partial:\ltwo C_1(H_T)\to \ltwo C_0(H_T)=H_T$; hence,
$\ltwo \ch_0(H_T)=G_T$.

Denote by $R(\cn)$ the Grothendieck group of Hilbert $\cn$-modules.  If $F$
is such a Hilbert module, $[F]$ denotes its class in $R(\cn)$.  It follows
from additivity of dimension that the function $F\to \dim_{\cn}F$ induces a
homomorphism $\dim_{\cn}:R(\cn)\to \bR$.

\begin{corollary}\label{c:ac1}
For $\bq\in \ol{\car}$ and $T\in\cs$, the following formulas hold in
the representation group $R(\cn)$:
    \begin{align*}
    [G_T]&=\sum_{U\in \cs_{\geq T}}\geps (U-T)[H_U]\\
    [H_T]&=\sum_{U\in \cs_{\geq T}}[G_U].
    \end{align*}
\end{corollary}

\begin{proof}
Note that in $\ltwo C_{\ast}(H_T)$ the boundary maps are maps of
Hilbert $\cn$-modules.  Hence,
the first formula follows from Theorem~\ref{t:acycliccomplex} by taking the
Euler characteristics.  The second formula follows from this and
the M\"{o}bius Inversion Formula.
\end{proof}

\begin{corollary}\label{c:ac2}
Suppose $\bq\in \ol{\car}$ and $T\in \cs$.  Then
\( \dim_{\cn} G_T=W^T(\bq)/W(\bq) \).
\end{corollary}

\begin{proof}
By Lemma~\ref{l:HS} (iii), $\dim_{\cn}H_U=1/W_U(\bq^{-1})$.  So,
\[
\dim_{\cn} G_T=\sum_{U\in \cs_{\geq T}}\frac{\geps (U-T)}{W_U(\bq^{-1})}
=\frac{W^T(\bq)}{W(\bq)},
\]
where the first equality is by Corollary~\ref{c:ac1} and the
second by Lemma~\ref{l:growth}~(iii)(b).
\end{proof}

\begin{lemma}\label{l:as-t}
If $\bq\in \ol{\car}$ and $U \subseteq S$, then
$$\sum_{\substack{T\in \cs \\T\subseteq U}}G_{T}a_{S-U}$$
is a dense subspace of $A_{S-U}$ and a direct sum decomposition.
Moreover, if $T \in \cs$, then the right multiplication by
$a_{S-T}$ induces a weak isomorphism $G_T \to G_Ta_{S-T}$.
\end{lemma}

\begin{proof}
As in Section~\ref{s:growth}, $X_{S-U}$ denotes the set of 
$(\emptyset,S-U)$-reduced elements.  As in \cite[Ex. 26]{bourbaki}, 
$X_{S-U}$ is the disjoint union of the $W^T$, $T\subseteq U$.  
Hence, $X_{S-U}(\bq)=\sum_{T\subseteq U} W^T(\bq)$. Dividing this by 
$W(\bq)$ and using Lemma~\ref{l:growth}~(ii), we get
\[
\frac{1}{W_{S-U}(\bq)} =\sum_{T\subseteq U}\frac{W^{T}(\bq)}{W(\bq)}.
\]

By Lemma~\ref{l:sumGT},
$$\ltwo=\ol{\sum_{T\in\cs}G_{T}}.$$
Multiplying on the right by $a_{S-U}$ and using
Lemma~\ref{l:GUaS-V} we obtain:
$$A_{S-U}=\ol{\sum_{T\in\cs}G_{T}a_{S-U}}=\ol{\sum_{T\subseteq 
U}G_{T}a_{S-U}}.$$
By Lemma~\ref{l:AS}~(iii), 
$$\dim_{\cn} A_{S-U} = \frac{1}{W_{S-U}(\bq)} =
\sum_{T\subseteq U}\frac{W^{T}(\bq)}{W(\bq)}.$$ On the other hand,
$$\dim_{\cn} \ol{G_{T}a_{S-U}} \le
\dim_{\cn}G_{T} = \frac{W^{T}(\bq)}{W(\bq)}.$$ It follows that
each of the above inequalities  is an equality and hence, that $G_T$ is
weakly isomorphic to $G_Ta_{S-T}$ and the sum is direct.
\end{proof}

\begin{Remark} In what follows we will use the symbol $\biguplus$ to
denote the sum of submodules of $\ltwo$, once we have proved that the sum
is direct. 
\end{Remark}

Since $G_V=0$ for nonspherical $V$ and $\bq \in \ol{\car}$, we
can restate Lemma~\ref{l:as-t} as follows:
$$A_{S-U}=\ol{\biguplus_{V\subseteq U}G_V a_{S-U}}.$$
Letting  $U=S$, we get the following corollary.

\begin{corollary}\label{c:Gdecomp}
If $\bq\in \ol{\car}$, then
\[
\sum_{V\subset S}G_V
\]
is a dense subspace of $\ltwo$ and a direct sum decomposition.
\end{corollary}

The fact that the sum of $G_V$ is direct has the following two
corollaries.

\begin{corollary}\label{c:GTint}
Let $\ca$ and  $\cb$ be collections of subsets of $S$. If $\bq \in
\ol{\car}$, then
$$ \ol{\biguplus_{U \in \ca} G_{U}}\cap \ol{\biguplus_{U \in \cb} G_{U}}=
\ol{\biguplus_{U \in \ca \cap \cb} G_{U}}$$
\end{corollary}

\begin{corollary}\label{c:HTsumGT} If $\bq \in
\ol{\car}$ and  $V  \subseteq S$, then
$$ H_V=\ol{\biguplus_{U \supseteq V} G_{U}}.$$
\end{corollary}

\begin{lemma}\label{l:DTGT}
If $\bq\in \ol{\car}$ and  $V  \subseteq S$, then
$$D_V=\ol{G_{V}a_{S-V}}.$$
In particular, $D_V=0$ if $V\not\in \cs$.
\end{lemma}

\begin{proof}
Since, by definition,  $D_V\subseteq A_{S-V}$,
$D_V = D_V a_{S-V}$ and since 
$D_V\subseteq \Big(\sum _{U\subset V}A_{S-U}\Big)^\perp$, we have:
$$D_V  \subseteq \Big(\sum _{U\subset V}A_{S-U}\Big)^\perp a_{S-V}.$$
Using equations \eqref{e:intersectperp}, we compute:
$$\Big(\sum _{U\subset V}A_{S-U}\Big)^\perp = \bigcap_{U\subset V}
A_{S-U}^\perp =\bigcap_{U\subset V} \ol{\sum_{s\in S-U} H_s}.$$
By Corollary \ref{c:HTsumGT}, $H_s=\ol{\biguplus_{X \ni s} G_{X}}$.
Therefore,
$$\Big(\sum _{U\subset V}A_{S-U}\Big)^\perp =
\bigcap_{U\subset V} \ol{\sum_{s\in (S-U)}\  \biguplus_{X \ni s} G_{X}}=
\bigcap_{U\subset V} \ol{\biguplus_{X \nsubseteq U} G_{X}}.$$
Using Corollary~\ref{c:GTint} we obtain:
$$\Big(\sum _{U\subset V}A_{S-U}\Big)^\perp =
\ol{\biguplus_{X \nsubseteq U \  \forall U\subset V} G_{X}}=
\ol{\biguplus_{X\not\subset V} G_{X}}.$$
Thus,
$$D_V \subseteq \left( \ol{\biguplus_{X\not\subset V} G_{X}} \right) 
a_{S-V}=
\ol{\sum_{X\not\subset V} G_{X} a_{S-V}}.$$
By Lemma~\ref{l:GUaS-V}, the only nonzero term in the last sum is when 
$X=V$.
Therefore, $D_V\subseteq \ol{G_V a_{S-V}}$.

To prove the opposite inclusion,  note that, by Lemma~\ref{l:GUaS-V}, for 
all $U\subset V$,  $G_Va_{S-V}a_{S-U}= G_Va_{S-U} = 0$.
Therefore, since  $\Ker a_{S-U}=A_{S-U}^\perp$, we have
$G_Va_{S-V}\subseteq A_{S-U}^\perp$ for all $U\subset V$. Since
$G_Va_{S-V}\subseteq  A_{S-V}$, it follows from the definition of $D_V$ 
that $G_Va_{S-V}\subseteq D_V$.
\end{proof}

We shall need the following Decomposition Theorem. (Of course,
there is also a corresponding version with the $D_V$ replaced by
$G_V$.)

\begin{theorem}[The Decomposition Theorem]\label{t:decomp}
If $\bq \in \ol{\car} \cup \ol{\car^{-1}}$, then
$$\sum_{V\subseteq S} D_V$$
is direct and a dense subspace of $\ltwo$. Moreover, if $\bq \in 
\ol{\car}$,
then the only nonzero terms in the sum are those with $V \in \cs$, and if
$\bq \in \ol{\car^{-1}}$, then the only nonzero terms in the sum are those 
with
$S-V \in \cs$.
\end{theorem}

\begin{proof}
If $\bq \in \ol{\car}$, then we let $U=S$ in Lemma~\ref{l:sumGT} to obtain:
$$\ltwo=\ol{\sum_{V\subseteq S} D_V}.$$
The assertion follows, since by Lemma~\ref{l:DTGT}, all
nonspherical $V$ have $0$ contributions, and by Lemmas~\ref{l:DTGT},
\ref{l:as-t} and Corollary~\ref{c:Gdecomp}, the dimensions of the nontrivial
terms add up to $1$.

If  $\bq \in \ol{\car^{-1}}$, then the result follows from 
Corollary~\ref{c:Gdecomp} by applying the
$j$-homomorphism.
\end{proof}

\begin{corollary}\label{c:dirint}
Let $\ca$ be a collection of subsets of $S$ and
let $U \subseteq S$.

If $\bq \in \ol{\car} \cup \ol{\car^{-1}}$, then
$$D_U \cap \ol{\biguplus_{U \in \ca} D_{V}}=
\begin{cases}
0 &\text{ if } U\not\in \ca,\\
D_U  &\text{ if } U\in \ca.
\end{cases}$$
\end{corollary}

\begin{corollary}\label{c:sumDT}
If $\bq \in \ol{\car} \cup \ol{\car^{-1}}$ and $U\subseteq S$, then
$$ A_U=\ol{\biguplus_{V \subseteq S-U} D_{V}}.$$
\end{corollary}

\begin{corollary}[Compare Lemma~\ref{l:growth} and Corollary~\ref{c:ac2}]\label{c:ac2.1}
Suppose $T\in\cs$.
\begin{enumeratei}
\item
For $\bq\in \car$, 
\( \dim_{\cn} D_{T}=W^T(\bq)/W(\bq) \).
\item
For $\bq\in \car^{-1}$,
\( \dim_{\cn} D_{S-T}=W^T(\bq^{-1})/W(\bq^{-1}) \).
\end{enumeratei}
\end{corollary}

\begin{proof}
(i)  By Lemma~\ref{l:DTGT}, $a_{S-T}$ maps $G_T$ monomorphically onto a 
dense subspace of $D_T$.  So, \( \dim_{\cn} D_{T}=\dim_{\cn} 
G_{T}=W^T(\bq)/W(\bq) \), where the second equality is by 
Corollary~\ref{c:ac2}.

(ii)  For  $\bq\in \car^{-1}$, the following formulas hold in the representation 
ring $R(\cn)$,  
    \begin{align*}
	[A_T]&=\sum_{U\in \cs_{\geq T}}[D_{S-U}]\\
    	[D_{S-T}]&=\sum_{U\in \cs_{\geq T}}\geps (U-T)[A_U].
    \end{align*}
where the first formula is from Corollary~\ref{c:sumDT} and the second 
follows from the first by the M\"{o}bius Inversion Formula.  So, as in 
Corollary~\ref{c:ac2}, 
\[
\dim_{\cn} D_{S-T}=\sum_{U\in \cs_{\geq T}}\frac{\geps (U-T)}{W_U(\bq)}
=\frac{W^T(\bq^{-1})}{W(\bq^{-1})},
\]
where the second equality is Lemma~\ref{l:growth}~(iii)(b).
\end{proof}

In Section~\ref{s:solomon} we will need the following version of
Lemmas~\ref{l:as-t} and \ref{l:DTGT}.  Its proof is essentially the same 
as the proofs of these lemmas, except that we use Theorem~\ref{t:decomp} 
and its corollaries instead of the corresponding statements involving 
the $G_U$.

\begin{lemma}[Compare Lemmas~\ref{l:as-t} and \ref{l:DTGT}]\label{l:GTDT}
Suppose $\bq\in 
\ol{\car}$ 
and $U \subseteq S$.  Then
\begin{enumeratei}
\item
$$\sum_{\substack{T\in \cs \\T\subseteq U}}D_{T}h_U$$
is a dense subspace of $H_U$ and a direct sum decomposition.
Moreover, if $T \in \cs$, then the right multiplication by
$h_T$ induces a weak isomorphism $D_T \to D_Th_T$.
\item
$G_U=\ol{D_Uh_U}.$
\end{enumeratei}
\end{lemma}

\section{Decoupling cohomology}\label{s:decouple}
We retain notation from Sections~\ref{s:spaces} and \ref{s:weighted}, e.g., $Z$ is a finite CW complex, 
$(Z_s)_{s\in S}$ is a family of subcomplexes and $\cu = (W\times
Z)/\sim$. As in (\ref{e:Z^T}), given $U\subseteq S$, $Z^U$ denotes the union of mirrors $Z_s$, $s\in U$.

For any Hilbert $\cn$-submodule $E$ of $\ltwo$, define
    \begin{equation*}
    \ltwo C^i(\cu;E):=\Phi ^{-1} (C^i(Z)\otimes E),
    \end{equation*}
where $\Phi:\ltwo C^i(\cu)\hookrightarrow C^i(Z)\otimes \ltwo$ is the
monomorphism defined in (\ref{e:Phi}).  In other words,
\begin{equation}\label{e:phi-1}
\ltwo C^i(\cu;E)=\bigoplus_{c\in Z^{(i)}}(\phi_c)^{-1}(A_{S(c)}\cap E),
\end{equation}
where $S(c)$ is the subset of $S$ defined in (\ref{e:Sc}) and 
$\phi_c:L^2(Wc,\muq)\to A_{S(c)}$ is the isomorphism defined in
(\ref{e:phi}).

\begin{proposition}\label{p:DTcoeff}
Suppose $\bq\in \car \cup \car^{-1}$.
Then the map $\Phi$ restricts to an isomorphism of cochain
complexes:
\[
\ltwo C^{\ast}(\cu;D_{U})\mapright{\cong}C^{\ast}(Z,Z^{U})\otimes
D_{U}.
\]
\end{proposition}

\begin{proof}  (Compare \cite[the proof of Theorem B]{d87}.)
Let $c\in Z$ be an $i$-cell. By Corollary~\ref{c:sumDT}, 
$$A_{S(c)}=\ol{\biguplus_{V \subseteq S-S(c)} D_{V}}.$$
If $c\nsubseteq Z^{U}$, then $S(c) \subseteq S-U$ and therefore, by
Corollary~\ref{c:dirint}, $A_{S(c)} \cap D_U=D_U$ and so,
by (\ref{e:phi-1}), $\phi_c:\ltwo C^i(Wc;D_{U})\to D_{U}$ is
an isomorphism. If $c\subseteq Z^{U}$, then $S(c) \nsubseteq S-U$
and therefore, $A_{S(c)}\cap D_{U} = 0$ and so, by
(\ref{e:phi-1}), $\ltwo C^i(Wc;D_{U})=0$. Hence, a cochain in
$C^i(Z)\otimes D_U$ is in the image of the restriction of $\Phi$ if and 
only if 
it evaluates to $0$ on the orbit of every $i$-cell $c\subseteq
Z^{U}$.
\end{proof}

Suppose $\bq \in \car \cup \car^{-1}$. Let
$\Theta_U:C^{\ast}(Z,Z^{U})\otimes D_{U} \to \ltwo C^{\ast}(\cu;D_{U})
\hookrightarrow C^{\ast}(\cu)$ be the inverse of the
isomorphism of Proposition~\ref{p:DTcoeff}.  Define
\[
\Theta:\bigoplus _{U\subseteq S}C^{\ast}(Z,Z^U)\otimes
D_U\mapright{} \ltwo C^{\ast}(\cu)
\]
to be the sum of the $\Theta_U$.

\begin{proposition}\label{p:weakiso}
If $\bq\in \car \cup \car^{-1}$, then
\[
\Theta:\bigoplus _{U\subseteq S}C^{\ast}(Z,Z^U)\otimes
D_U\mapright{} \ltwo C^{\ast}(\cu)
\]
is a weak isomorphism of cochain complexes of Hilbert $\cn$-modules.
\end{proposition}

\begin{proof}
We have:
\[
\bigoplus_{U\subseteq S} \ltwo C^{\ast}(\cu;D_U)=
\ltwo C^{\ast}(\cu;\bigoplus_{U\subseteq S} D_U).
\]
By the Decomposition Theorem (Theorem~\ref{t:decomp}),
we have a weak isomorphism, $\ltwo C^{\ast}(\cu;\bigoplus D_U) \to \ltwo
C^\ast(\cu)$.
Combining this with the isomorphism of Proposition~\ref{p:DTcoeff}, the proposition follows.
\end{proof}

A weak isomorphism of chain complexes of Hilbert modules induces a weak
isomorphism on the level of reduced cohomology (\cite[Lemma 5]{dl}).  
Furthermore, 
if two Hilbert $\cn$-modules are weakly isomorphic, then they are
isometric (\cite[Lemma 2.5.3]{eckmann}).  So, we have the following corollary to
Proposition~\ref{p:weakiso}.

\begin{theorem}[{(Compare \cite{d87}, \cite [Theorem A]{d98})}]\label{t:decoupled}\ 
\begin{enumeratei}
\item
If $\bq\in \car$, then
\[
\ltwo \ch^{\ast}(\cu)\cong\bigoplus _{T\in \cs}H^{\ast}(Z,Z^T)\otimes D_T.
\]
\item
If $\bq\in \car^{-1}$, then
\[
\ltwo \ch^{\ast}(\cu)\cong\bigoplus _{T\in \cs}H^{\ast}(Z,Z^{S-T})\otimes
D_{S-T}.
\]
\end{enumeratei}
\end{theorem}

The special case $\cu=\gS$ is the following.

\begin{theorem}\label{t:decoupledsigma}\ 
\begin{enumeratei}
\item
\textup{(Theorem~\ref{t:dym} or \cite [Cor. 10.4]{dym})}
If $\bq\in \car$, then $\ltwo H^{\ast}(\gS)$ is concentrated in dimension
$0$ and
\[
\ltwo H^0(\gS)=\ltwo \ch^0(\gS)\cong A_S.
\]
   So, $b^0_\bq(\gS)=\chi_\bq(\gS)=1/W(\bq)=dim_{\cn}A_S$.
\item
If $\bq\in \car^{-1}$, then
\[
\ltwo \ch^{\ast}(\gS)\cong\bigoplus _{T\in \cs}H^{\ast}(K,K^{S-T})\otimes
D_{S-T}.
\]
So,
\[
b^i_\bq(\gS)=\sum_{T\in\cs}\frac{W^T(\bq^{-1})}{W(\bq^{-1})}b^i(K,K^{S-T}),
\]
where $b^i(K,K^{S-T})=\dim_\bR H^i(K,K^{S-T};\bR)$.
\end{enumeratei}
\end{theorem}

\noindent
(In the formula for $b^i_\bq(\gS)$ in Theorem~\ref{t:decoupledsigma}~(ii) we 
have used the  
formula for $\dim_{\cn}D_{S-T}$ from Corollary~\ref{c:ac2.1}.)

\section{A generalization of a theorem of Solomon}\label{s:solomon}
When $W$ is finite and $\bq=\bone$, L. 
Solomon \cite{s} proved some results  very similar to the Decomposition 
Theorem (Theorem~\ref{t:decomp}).  In this special case, formulas 
(\ref{e:aT5}) and (\ref{e:hT5}) for the idempotents $a_T$ and $h_T$ become
    \begin{align*}
    a_T&:=\frac{1}{\Card(W_T)}\sum_{w\in W_T}e_w \text{ and}\\
    h_T&:=\frac{1}{\Card(W_T)}\sum_{w\in W_T}\geps_w e_w
    \end{align*}
and we recognize $a_T$ and $h_T$ as the familiar elements
of  ``symmetrization'' and ``alternation'' in the group algebra
$\bR[W_T]$.

\begin{solomonthm}\cite{s}\label{t:solomon}
  Suppose $W$ is finite.   Then there 
are direct sum decompositions of the regular representation:
\begin{align*}
L^2(W)&=\sum_{T\subseteq S} L^2(W)a_Th_{S-T},\\
L^2(W)&=\sum_{T\subseteq S} L^2(W)h_{S-T}a_T.
\end{align*}
\end{solomonthm}

Our generalization of Solomon's Theorem is the following.

\begin{theorem}\label{t:gensolomon}
\begin{enumeratei}
\item
If $\bq\in\car$, then 
\[
\sum_{T\in\cs}\ol{\ltwo h_Ta_{S-T}} \text{ and } 
\sum_{T\in\cs}\ol{\ltwo a_{S-T}h_T}
\]
are direct sum decompositions and dense subspaces of $\ltwo $.
\item
If $\bq\in\car^{-1}$, then 
\[
\sum_{T\in\cs}\ol{\ltwo h_{S-T}a_T}\text{ and } 
\sum_{T\in\cs}\ol{\ltwo a_Th_{S-T}}
\]
are direct sum decompositions and dense subspaces of $\ltwo$.
\end{enumeratei}
\end{theorem}

This is an immediate consequence of Corollary~\ref{c:Gdecomp},
Theorem~\ref{t:decomp} and the following theorem.
\begin{theorem}\label{t:dgl}
Suppose $T\in\cs$.
\begin{enumeratei}
\item
If $\bq\in\car$, then $\ol{\ltwo a_{S-T}h_T}=G_T$ and 
$\ol{\ltwo h_Ta_{S-T}}=D_T$.  
\item
If $\bq\in\car^{-1}$, then $\ol{\ltwo a_Th_{S-T}}=G_{S-T}$ and 
$\ol{\ltwo h_{S-T}a_T}=D_{S-T}$.  
\end{enumeratei}
\end{theorem}

\begin{proof}
(i)  Suppose $\bq\in\car$.
By Lemma~\ref{l:as-t}, right multiplication by $a_{S-T}$ is a weak 
isomorphism from $G_T$ to $G_Ta_{S-T}$.  So, by Lemma~\ref{l:DTGT},
$\ol{\ltwo h_Ta_{S-T}}=\ol{G_Ta_{S-T}}=D_T$.  Similarly, by 
Lemma~\ref{l:GTDT},  
$\ol{\ltwo a_{S-T}h_T}=G_T$. 

(ii)  Applying the $j$-isomorphism to the two equations in (i), we get 
the two equations in (ii).
\end{proof}

\begin{remark} 
It seems probable that $\ol{\ltwo a_{S-U}h_U}=G_U$ and 
$\ol{\ltwo h_Ua_{S-U}}=D_U$ whenever $\bq\in \car_{S-U}\cap \car_U^{-1}$ 
(so that $h_U$ and $a_{S-U}$ are both defined).
\end{remark}

\begin{remark}\label{r:kl}
In \cite{kl} Kazhdan and Lusztig study the regular representation of the Hecke algebra $\R$ on itself and they generalize Solomon's Theorem in a slightly different direction.  $W$ can be infinite.  First, they define a basis $\{C_w\}_{w\in W}$ for $\R$, called the ``Kazhdan-Lusztig basis.''  It has many good properties.  Next they partition of $W$ into ``left cells.''  This partition is strictly finer than the partition of $W$ into the $W^T$, $T\in \cs$.  Given a left cell $Z$, they define a certain subquotient $\ci_Z/\ci'_Z$ of $\R$ such that $\{C_w\}_{w\in Z}$ projects to a basis for the subquotient.  If we sum these representations over all $Z\subseteq W^T$, we obtain a representation analogous to our $D_{S-T}$.  (Compare \cite{ddjo06}.) So, our Decomposition Theorem (Theorem~\ref{t:decomp}) is a partial generalization of Kazhdan-Lusztig theory to the Hecke - von Neumann algebra $\cn$.  It seems likely that left cells can be used to get a further direct sum decomposition of the $D_{S-T}$, although we do not yet know how to prove this.
\end{remark}

\section{Relationship with ordinary homology and cohomology with compact 
supports}\label{s:ordinary}

As in Section~\ref{s:spaces}, $Z$ is a CW complex which is a strict fundamental 
domain for a $W$-action on $\cu$ ($=\cu(W,Z)$).

\begin{theorem}\label{t:ordinary}
\begin{enumeratei}
\item
For $\bq \in \car$, the canonical map, $\can:H_*(\cu;\bR)\to 
\ltwo\ch_*(\cu)$, is an injection with dense image.
\item
For $\bq\in\car^{-1}$, the canonical map, $\can:H^*_c(\cu;\bR)\to 
\ltwo\ch^*(\cu)$, is an injection with dense image.
\end{enumeratei}
\end{theorem}

\begin{proof}
First, using \cite{d98}, we prove statement (ii) for cohomology with compact 
supports.  Given $w\in W$, let $\In (w)\in \cs$ be as in 
Section~\ref{s:coxeter} and
let $Q_w$ be the ``positive'' fundamental domain 
for $W_{\In(w)}$ on $\cu$, containing the chamber $Z$.  
(So, $Q_w\cong \cu/W_{\In(w)}$.) Let $q_w:\cu \to wQ_w$ denote the 
composition of projection onto $Q_w$ with translation by $w$. The composition 
of the map induced by $q_w$ with the excision isomorphism 
$C_*(wQ_w, w\ol{Q_w-Z})\cong C_*(wZ,wZ^{S-\In(w)})\cong 
C_*(Z,Z^{S-\In(w)})$ induces a chain map $p_w:C_*(\cu)\to C_*(Z,Z^{S-\In 
(w)})$.  For each $T\in \cs$, define $p_T^*:C^*(Z,Z^{S-T})\otimes 
\bR^{(W^T)}\to  C^*_c(\cu;\bR)$ by $c\otimes e_w \to p_w^*(c)$ and extending 
linearly.  In other words, $p_T^*(c\otimes e_w)=e_w\ta_T (c)$, where $\ta_T$ is 
defined by (\ref{e:aT1}) (and where $e_w\ta_T$ acts on $C^*_c(\cu;\bR)$ as 
an element of the group algebra $\bR[W]$, not as an element of the Hecke algebra).  
$p_T^*$ will also denote the induced map on cohomology.
It is proved in \cite{d98} that $\oplus p^*_T:\oplus H^*(Z,Z^{S-T})\otimes 
\bR^{(W^T)}\to H^*_c(\cu)$ is an isomorphism.  Computations similar to 
those in
Section~\ref{s:decouple} give an isomorphism 
$\gf_{S-T}:\ltwo\ch^*(\cu;G_{S-T})\to H^*(Z,Z^{S-T})\otimes G_{S-T}$.  It follows 
that we have a commutative diagram:
\[
    \minCDarrowwidth 0.5cm
    \begin{CD}
\bigoplus H^*(Z,Z^{S-T})\otimes \bR^{(W^T)}@>\oplus p_T^*>> H^*_c(\cu;\bR)\\
  @.  @VV\can V\\
      \smash{\Bigg{\downarrow}} \scriptstyle{g} && \ltwo\ch^*(\cu)\\
@. @VV\oplus \pi_{S-T} V\\
\bigoplus H^*(Z,Z^{S-T})\otimes G_{S-T}@<\oplus \gf_{S-T}<<
\bigoplus \ltwo\ch^*(\cu;G_{S-T})
\end{CD}
\]
Here $\pi_{S-T}$ is the coefficient homomorphism induced by orthogonal projection 
$\ltwo\to G_{S-T}$ and $g:=\oplus g_T$,
where $g_T:\bR^{(W^T)}\to G_{S-T}$ is induced by $e_w\to e_w\ta_Th_{S-T}$.  
In other words, up to weak isomorphisms of $\cn$-modules, the canonical map 
$H^*(\cu;\bR)\to \ltwo \ch^*(\cu)$ is identified with $g$.  By 
Theorem~\ref{t:dgl}~(ii), for $\bq\in \car^{-1}$,  $\{e_w\ta_Th_{S-T}\}_{w\in 
W^T}$ spans a dense subspace of $G_{S-T}$.  Hence, each $g_T$ is injective with 
dense image.  This proves (ii). 

The canonical map in (i) is induced by the composition of chain maps:
\[
(C_*(\cu;\bR),\partial)\hookrightarrow 
(L^2_{\bq^{-1}}C_*(\cu),\partial)\mapright{\cong}  \ltwo 
(C_*(\cu),\partial^{\bq}),  
\]
where the second map is the isomorphism of Lemma~\ref{l:qtoq-1}.  For 
each $T\in\cs$, let 
\[
\hat{h}_T:=\sum_{w\in W_T} \geps_w e_w 
\]
and let $\hat{h}_T:C_*(Z,Z^T)\otimes \bR^{(W^T)}\to C_*(\cu;\bR)$ be 
defined by $x\otimes e_w\to e_w\hat{h}_T (x)$.  The element $\tilde{h}_T$ defined in (\ref{e:1})) differs from $\hat{h}_T$ by inserting a $q_w^{-1}$ in front of each $e_w$. By \cite{d87}, 
the induced map in homology $\oplus(\hat{h}_T)_*:H_*(Z,Z^T)\otimes 
\bR^{(W^T)}\to H_*(\cu;\bR)$  is an isomorphism.  Hence, the map induced by $\oplus(\hat{h}_T)_*$ is also an isomorphism.  
We have a commutative diagram:
\[
    \minCDarrowwidth 0.5cm
    \begin{CD}
\bigoplus H_*(Z,Z^T)\otimes \bR^{(W^T)}@>\oplus (\tilde{h}_{T})_*>> H_*(\cu;\bR)\\
@.   @VV\can V\\
\smash{\Bigg\downarrow}\scriptstyle{d} && \ltwo\ch_*(\cu)\\
@. @VV\oplus \pi'_T V\\
\bigoplus H_*(Z,Z^T)\otimes D_T@<\oplus (\Theta_T)_*<<
\bigoplus \ltwo\ch_*(\cu;D_T)
\end{CD}
\]
Here $\pi'_T$ is the coefficient homomorphism induced by orthogonal projection 
$\ltwo\to D_T$ and $\Theta_T$ is the isomorphism of Section~\ref{s:decouple}.  
Also, $d:=\oplus d_T$, 
where $d_T:\bR^{(W^T)}\to D_T$ is induced by $e_w\to e_w\tilde{h}_Ta_{S-T}$.  
In other words, up to weak isomorphisms of $\cn$-modules, the canonical map 
$H^*(\cu;\bR)\to \ltwo \ch^*(\cu)$ is identified with $d$. By 
Theorem~\ref{t:dgl}~(i), each $d_T$ is injective with dense image.
\end{proof}

\section{$L^2$-cohomology of buildings}\label{s:buildings}
As in \cite{ronan}, a \emph{building} consists of the following data:
\begin{itemize}
\item
a set $\Phi$,
\item
a Coxeter system $(W,S)$,
\item
a collection of equivalence relations on $\Phi$ indexed by $S$.
\item
a function $\gd:\Phi\times \Phi\to W$.
\end{itemize}
This data must satisfy certain additional conditions which we
will explain below.  One condition is that for  $s\in S$, each
$s$-equivalence class contains at least two elements. 

The elements of $\Phi$ are called  \emph{chambers}.  Given $s\in S$,
two chambers $\gf$
and $\gf'$ are \emph{$s$-equivalent} if they are equivalent via the
equivalence relation corresponding to $s$.   If, in addition,
$\gf\neq\gf'$, they are \emph{$s$-adjacent}.  A \emph{gallery} is a
sequence $(\gf_0,\dots,\gf_n)$ of adjacent chambers; its \emph{type} is the
word $(s_1,\dots,s_n)$ in the letters of $S$, where $\gf_{i-1}$ and 
$\gf_i$ are
$s_i$-adjacent.  Given $T\subset S$, $(\gf_0,\dots,\gf_n)$  is a
\emph{$T$-gallery}
if each $s_i\in T$.
The gallery is \emph{reduced} if $w=s_1\cdots s_n$ is a
reduced expression.  

Another condition for $\Phi$ to be a building is that there exist a 
\emph{$W$-valued distance function}  
$\gd:\Phi\times \Phi\to W$.  This means that there is a reduced gallery of 
type 
$(s_1,\dots,s_n)$ from 
$\gf$ to $\gf'$ if and only if $s_1\cdots s_n$ is a reduced expression
for $\gd(\gf,\gf')$.  

The \emph{$s$-mirror} (or ``$s$-panel'') of a
chamber $\gf$
is the $s$-equivalence class containing $\gf$.  More generally, given a
subset $T\subset S$, the \emph{$T$-residue} of $\gf$ is the
$T$-gallery connected component containing $\gf$.  Each such $T$-residue
is naturally a building with associated Coxeter system $(W_T,T)$.
The residue is
\emph{spherical} if $T$ is a spherical.  

\begin{example}\label{ex:trees}
(\emph{Trees}).  Suppose $W$ is the infinite dihedral group (so that
$\Card(S)=2$).   Any  tree is
bipartite, i.e., its vertices can be labeled by the two elements of $S$ so
that the  vertices of any edge have distinct labels.  Suppose $T$ is 
a tree  with such a labeling and  suppose no  vertex of $T$ is of valence 
$1$.  Let
$\Phi$ be its set of edges.  Given $s\in S$, call two edges 
\emph{$s$-equivalent} if they meet at a vertex of type $s$. An 
$\{s\}$-residue
is the set of edges in the star of a vertex of type $s$.  A gallery in 
$\Phi$ corresponds to an edge path in $T$.  The type of the gallery is the 
word obtained by taking the types of the vertices crossed by the 
corresponding edge path.  This word is reduced if and only if the edge 
path does not backtrack.  Given two edges $\gf, \gf'$ of $T$, there is a 
(unique) minimal 
gallery connecting them.  The corresponding word represents an element of 
$w\in W$ and $\gd (\gf,\gf'):=w$.  Thus, every such tree $T$ defines a 
building of type $(W,S)$.  Not surprisingly, we will define the 
``geometric 
realization of a building'' so that for the building $\Phi$ corresponding 
to $T$, its geometric realization  will be $T$.  
\end{example}

A building
$\Phi$ of type $(W,S)$ has \emph{finite thickness} if
for each $s\in S$, each $s$-equivalence class is finite.
If $\Phi$ has finite thickness, then it follows from
the existence of a $W$-distance function that each of its spherical residues 
is finite.

Let us say that $\Phi$ is \emph{regular} if for each $s\in S$,  the
$s$-equivalence
classes have constant cardinality.  When finite, we denote this number
by $q_s+1$.  It is known (\cite{ronan}) that if $s$ and $s'$ are 
conjugate
in $W$, then $q_s=q_{s'}$.  Let $I$ be
the set of conjugacy classes of elements in $S$.  Then for any regular
building $\Phi$, the integers $q_s$ define an
$I$-tuple $\bq$ called the \emph{thickness vector} of $\Phi$.

A group $G$ of automorphisms of a building is \emph{chamber transitive} if
it acts transitively on $\Phi$.  When this is the case, we have $\Phi\cong
G/B$, where $B$ denotes the stabilizer of some given chamber $\zeta$.  If
$G_s$ denotes the stabilizer of the $s$-mirror containing $\zeta$, then
the chambers $g\zeta$ and $g'\zeta$ are $s$-equivalent if and only if $g$
and $g'$ belong to the same coset of $G_s$.   Obviously, if $G$ is chamber
transitive, then the building is regular.
For the
remainder of this section, we suppose that $\Phi$ has finite thickness and
that $G$ is chamber transitive.

Given a subset $T$ of $S$,  denote the stabilizer of the
$T$-residue containing $\zeta$ by $G_T$. Thus, $G_\emptyset=B$ and
$G_{\{s\}}=G_s$.  If $\Phi$ has finite thickness and $T\in\cs$, then the
number of elements in a $T$-residue is $\Card (G_T/B)$.  (This number is
known to be $W_T(\bq)$.)

Fix a chamber $\zeta\in\Phi$ and let $r$ (or $r_\zeta$) denote the
function $\Phi\to W$ defined by $\gf\to \gd(\zeta,\gf)$.  Since $\Phi\cong
G/B$, we can regard $r$ as a function from $G/B$ to $W$.
Since $B$ fixes $\zeta$,  $r:G/B\to W$ is $B$-invariant.
In other words, $r$ induces a map $\ol{r}:B\backslash G/B\to W$.

A \emph{Tits system} is a quadruple $(G,B,N,S)$, where $G$
is a  group, $B$ and $N$ are subgroups of $G$, $W:=N/N\cap B$, $S$ is a
subset of $W$ and where the conditions listed in 
\cite[pp. 15--26]{bourbaki}
are satisfied.  Given $w\in W$, put $C(w):=BwB$.  The conditions imply that
\begin{itemize}
\item
For each $s\in S$, $G_s:=B\cup C(s)$ is a subgroup of $G$.
\item
$(W,S)$ is a Coxeter system.
\item
There is a building with set of chambers $G/B$ such that two chambers $gB$
and $g'B$ are $s$-equivalent if and only if $gG_s=g'G_s$.
\item
Suppose $r:G/B\to W$ is defined by $gB\to \gd(B,gB)$ where $\gd$ is 
$W$-distance in the building.  Then the induced map
$\ol{r}:B\backslash G/B\to W$ is a bijection.
\end{itemize}
One says that the building \emph{comes from a $BN$-pair}.

\begin{definition}\label{d:ra}
The Coxeter system $(W,S)$ is \emph{right-angled} if $m_{st} =2$ or
$\infty$ for each pair $s,t$ of distinct elements in $S$.
\end{definition}

\begin{example}\label{ex:rabldg}
(\emph{Regular right-angled buildings}, \cite[pp. 112--113]{d98a}).
For any right-angled Coxeter
system $(W,S)$ (cf. Definition~\ref{d:ra})
and any $S$-tuple $\bq =(q_s)_{s\in S}$ of positive integers,
there is a regular building $\Phi$ of type $(W,S)$ with thickness vector
$\bq$.   In the case where $W$ is the infinite dihedral group this is
well-known: as in Example~\ref{ex:trees}, the building
is a (bipartite) tree with edge set $\Phi$,  it is ``regular'' in the sense
that  for each $s\in S$ there are exactly
$q_s+1$ edges meeting at each vertex of type $s$.

In the general case, the construction goes as follows.  For each
$s\in S$, choose a finite
group $\gG_s$ with $\Card (\gG_s)=q_s+1$ and let $\gG$ be the ``graph
product'' of the $(\gG_s)_{s\in S}$ where the graph is the $1$-skeleton
of $L$.   In other words,  $\gG$ is the quotient of the free product of
the $(\gG_s)_{s\in S}$ by the normal subgroup generated by all
commutators $[g_s,g_t]$ with $g_s\in \gG_s$, $g_t\in \gG_t$ and $m_{st}=2$.
As in \cite{d98a}, we get a building with $\Phi =\gG$ and with two
elements $g, g'\in \gG$ in an $s$-equivalence class if and only if
they determine same coset in $\gG/\gG_s$.  We leave the following two facts
as exercises for the reader:
\begin{itemize}
\item
Two regular right-angled buildings of a given type $(W,S)$ are isomorphic
if and only if they have the same thickness vector.
\item
Any regular right-angled building  comes from a $BN$-pair.
In other words, its full automorphism group $G$ is chamber transitive and
if $B$ denotes the stabilizer of a given chamber and $N$ the stabilizer of
some apartment containing that chamber, then there is a set of generators
$S$ for $W:=N/N\cap B$ so that $(G,B,N,S)$ is a Tits system.
\end{itemize}
\end{example}

\rk{Hecke algebras and functions on $B\backslash G/B$}  This paragraph
is taken from \cite[Ex. 22, pp. 56--57]{bourbaki}.

Suppose $G$ is a topological group and $B$ is a compact open
subgroup.  Let $C(G)$ denote the vector space of continuous real-valued
functions on $G$.  Let $\ga:G\to G/B$ and $\gb:G\to B\backslash G/B$ be
the natural projections.  Define subspaces $H\subseteq L\subseteq C(G)$ by
\[
L:=\ga^\ast \bR^{(G/B)}\text{ and } H:=
\gb^\ast \bR^{(B\backslash G/B)},
\]
where, as in Section~\ref{s:hecke}, for any set $X$, $\bR^{(X)}$ denotes 
the
vector space of finitely supported functions on $X$.

For each $gB\in G/B$, let $a_{gB}\in L$ be defined by $a_{gB}(x)=1$ for
$x\in gB$ and $a_{gB}(x)=0$ for $x\notin gB$.  Since $(a_{gB})$ is a basis
for $L$, there is a unique linear form on $L$ such that $a_{gB}\to 1$ for
all $gB\in G/B$.  We denote this form by $\gf\to \int \gf$ (since it
coincides with the Haar integral normalized by the condition that $\int
a_B=1$).

If $\gf\in L$ and $\psi\in H$, then for each $x\in G$, the function
$\theta_x:G\to
\bR$, defined by $\theta_x(y)=\gf(y)\psi(y^{-1}x)$, belongs to $L$.  The
function $\gf\ast\psi:x\to \int \gf(y)\psi(y^{-1}x)dy$ also belongs to
$L$.  Moreover, if $\gf\in H$, then $\gf\ast\psi\in H$.  The map
$(\gf,\psi)\to \gf\ast\psi$ makes $H$ into an algebra and $L$ into a
right $H$-module.  $H$ is called the \emph{Hecke algebra} of $G$ with
respect to $B$.

Next, suppose that $G$ is a chamber transitive automorphism group on a
building and that $r:G/B\to W$ is defined by taking the
$W$-distance from the chamber corresponding to $B$.  Let
$\gamma :=r\circ\ga:G\to W$ and $J:=\gamma^\ast (\bR^{(W)})\subseteq H$.

\begin{Remark}
If $(G,B,N,S)$ is a Tits system, then $\ol{r}:B\backslash G/B\to W$ is a
bijection and hence, $J=H$.
\end{Remark}

\begin{lemma}\label{l:subalgebra}
Suppose, as above, that a given  building
admits a chamber transitive automorphism group $G$ (so $G/B$ is the 
set 
of chambers).  Let $\bq$ be the thickness vector.  Then
\begin{enumeratei}
\item
$J$ is a subalgebra of $H$ and
\item
$J\cong \R$, the Hecke algebra of Section~\ref{s:hecke}.
\end{enumeratei}
\end{lemma}

\begin{proof}
Since $G$ is  chamber transitive,
$\gamma^*:{\bf R}^{(W)}\to J$ is an isomorphism of vector spaces. 
So, we only need to check
that $\gamma^*$ is an algebra homomorphism from ${\bf R_q}[W]$ to $H$.
Let $f_w=\gamma^*(e_w)$. Then $f_w$ is the characteristic function
of $\{g\in G\mid r(gB)=w\}$. In particular,  for each $s\in S$, $f_s$ is 
the characteristic function of 
$G_s-B$. We want to see that
\[
f_w*f_s=
\begin{cases}
f_{ws}, &\text{ if $l(ws)>l(w)$;}\\
q_sf_{ws}+(q_s-1)f_w, &\text{ if $l(ws)<l(w)$.}
\end{cases}
\]
By definition of convolution,
$$(f_w*f_s)(g)=
\int_Gf_w(x)f_s(x^{-1}g)dx=
\int_Gf_w(gu)f_s(u^{-1})du=
\int_{G_s-B}f_w(gu)du,$$
which is equal to the Haar measure of the set
$$U_g:=\{u\in G_s-B\mid r(guB)=w\}.$$
Let $C_0:=g_0B$ be the chamber which is $s$-adjacent to $gB$ 
and which is closest to $B$.   There are $q_s$ other chambers adjacent to
$gB$.  We list them as:  $C_1=g_1B,\ldots,C_{q_s}=g_{q_s}B$.  So, for $i>0$, 
$r(C_i)=r(C_0)s$. Notice that if 
$u\in G_s-B$, then 
$guB$ is $s$-adjacent to $gB$ and therefore, $guB$ is equal to some $C_i$.
So, if $r(guB)=w$, then $r(gB)=w$ or $ws$. In other words, 
if $r(gB)\not\in\{w,ws\}$, then $(f_w*f_s)(g)=0$. We now consider two 
cases.  Each case further divides into two subcases depending on whether 
$r(gB)=w$ or $ws$.

\rk{Case 1} $l(w)<l(ws)$. In this case $r(C_0)=w$ and $r(C_i)=ws$ for 
$i>0$. 

\item{a)} Suppose $r(gB)=w$. Then $gB=C_0$ and  $guB=C_i$ for  $i>0$, so that
$r(guB)=ws$. Thus, $U_g=\emptyset$ and $(f_w*f_s)(g)=0$.    

\item{b)} Suppose $r(gB)=ws$. Then $gB=C_k$, for some $k>0$, and
$$U_g=\{u\in G_s-B\mid guB=C_0\}=
(G_s-B)\cap g^{-1}g_0B.$$
Since $gB$ and $g_0B$ are $s$-adjacent and not equal,  
$g^{-1}g_0B\subseteq G_s-B$,
so that $U_g=g^{-1}g_0B$ has measure 1. Therefore, $(f_w*f_s)(g)=1$.
So, in Case 1, $f_w*f_s=f_{ws}$.

\rk{Case 2} $l(w)>l(ws)$. In this case $r(C_0)=ws$ and $r(C_i)=w$ for 
$i>0$.
\item{a)} Suppose $r(gB)=w$. Then $gB=C_k$ for some $k>0$. So, the set
$$U_g=\bigcup_{0<i}\{u\in G_s-B\mid guB=C_i\}=
\bigcup_{0<i\ne k}g^{-1}g_iB$$
has measure $q_s-1$.
\item{b)} Suppose $r(gB)=ws$. Then $gB=C_0$, and the set
$$U_g=\bigcup_{0<i}\{u\in G_s-B\mid guB=C_i\}=
\bigcup_{0<i}g^{-1}g_iB$$
has measure $q_s$.
So, in Case 2, $f_w*f_s=q_sf_{ws}+(q_s-1)f_w$.
\end{proof}

\rk{The geometric realization of a building}  Suppose $\Phi$ is a
building with associated Coxeter system $(W,S)$.  As in
Section~\ref{s:coxeter}, let $K$ be the geometric realization of $\cs$ and
$\gS$ the geometric realization of $W\cs$.  By (\ref{e:gscu}), 
$\gS=\cu(W,K)$,
where $\cu(W,K)=(W\times K)/\sim$ and where $\sim$ is the equivalence
relation defined in the beginning of Section~\ref{s:spaces}.
Following \cite[pp. 117--118]{d98a},
define the \emph{geometric realization} of $\Phi$ to be
    \begin{equation}\label{e:uphi}
    \cu(\Phi, K)=(\Phi \times K)/\sim,
    \end{equation}
where $(\gf,x)\sim (\gf',x')$ if and only if $x=x'$ and $\gf, \gf'$ belong
to the same $S(x)$-residue.  ($S(x)$ is defined in (\ref{e:Sz}).)

Since $K$ only involves the spherical subsets of $S$, $\cu(\Phi, K)$ only
involves the spherical residues of $\Phi$.  It follows that if $\Phi$ has
finite thickness, then $\cu (\Phi,K)$ locally finite.

We often write $X$ as a shorthand for $\cu(\Phi,K)$.

\rk{The von Neumann algebra of $G$}  Next suppose  $G$ is a chamber
transitive group of automorphisms of $\Phi$ and that $B$ is the stabilizer
of some fixed chamber $\zeta$.   $G$ acts as a group of homeomorphisms of
$X$, so give it the compact-open topology.  Then $B$ is a
compact open subgroup.  Let $\mu$ be Haar measure on $G$, normalized by 
the condition that $\mu(B)=1$.

We have the left regular representation of $G$ on $L^2(G)$.  The
\emph{von Neumann algebra} $\cN(G)$ consists of all $G$-equivariant bounded
linear endomorphisms of $L^2(G)$.

Any $\ga\in\cN(G)$ is represented by convolution with some distribution 
$f_\ga$.  This distribution need not be a function.  For example, if $\ga$ 
is the identity map on $L^2(G)$, then $f_\ga=\gd_1$ (the Dirac delta).
One would like to define the ``trace'' of  $\ga$ to be $f_\ga(1)$ 
whenever $f_\ga$ is a function.  However, since $f_\ga$ is well-defined 
only up to sets of measure $0$, we must proceed slightly differently.  

Suppose $\ga$ is a nonnegative self-adjoint element of $\cN(G)$.  Let  
$\gb$ be its square root.  If $f_\gb$ is a $L^2$ function, then put 
\[
\tr_{\cN(G)} \ga:= \norm{f_\gb}:=\left( \int_G f_\gb(x)^2 
d\mu\right)^{1/2}.
\]

This extends in the usual fashion to give a ``trace'' on
$(n\times n)$-matrices with coefficients in $\cN(G)$.  If $V$ is a closed
$G$-stable subspace of $\bigoplus L^2(G)$ and $\pi_V:\bigoplus L^2(G)\to 
\bigoplus
L^2(G)$ is orthogonal projection, then the \emph{von Neumann dimension}
of $V$ is defined by
    \[
    \dim_{\cN(G)}V:=\tr_{\cN(G)} \pi_V.
    \]

We identify $L^2(\Phi)=L^2(G/B)$ with the subspace of $L^2(G)$
consisting of the functions which are constant on each right coset $gB$, 
$g\in G$.
Orthogonal projection from $L^2(G)$ onto $L^2(G/B)$ is given by convolution
with the characteristic function of $B$.  In view of the assumption that
$\mu(B)=1$,
    \begin{equation*}
    \dim_{\cN(G)}L^2(G/B)=1.
    \end{equation*}

The map $r:G/B\to W$ defined by the $W$-distance from the base chamber
induces a bounded linear map $\ltwo (W)\to L^2(G/B)$ which we shall also
denote by $r$.  Since this map takes bounded elements of $\ltwo (W)$ to
bounded elements of $L^2(G/B)$, we get the following version of
Lemma~\ref{l:subalgebra}.

\begin{lemma}\label{l:vnsubalgebra}
The map $r:\ltwo (W)\to L^2(G/B)$ induces a monomorphism of von Neumann
algebras $r:\cn \to \cN(G)$.  (In particular, $r$ commutes with the $*$
anti-involutions on $\cn$ and $\cN(G)$.)
\end{lemma}

$L^2C^\ast(X)$ denotes the Hilbert space of square summable simplicial
cochains on $X$ and $\ch^\ast(X)$ denotes the subspace of harmonic
cocycles.  (Of course, the $\ch^\ast(X)$ are isomorphic to reduced
cohomology groups of the cochain complex $L^2C^\ast(X)$.)
Supposing $G$ is a chamber transitive automorphism group,
we have
\begin{equation*}
L^2C^i(X)=\bigoplus_{\gs\in K^{(i)}}L^2(G/G_\gs)\subset \bigoplus_{\gs\in
K^{(i)}} L^2(G),
\end{equation*}
where $G_\gs:=G_{S(\gs)}$ is the stabilizer of the $i$-simplex $\gs$.
($S(\gs)$ is the spherical subset defined in (\ref{e:Sc}).)
One then defines the \emph{$L^2$-Betti numbers of $X$ with respect to
$G$} by
    \begin{equation*}\label{e:bettiX}
    b^i(X;G)=\dim_{\cN(G)}\ch^i(X),
    \end{equation*}

The map $r:X\to \gS$ induces a map on cochains which we denote by the same
letter, i.e., $r$ is a cochain map from $\ltwo C^\ast(\gS)$ to 
$L^2C^\ast(X)$.
We also have ``transfer
maps'' on chains and cochains.  On the level of chains, the transfer
map sends a cell
$c$ of $\gS$ to $r^{-1}(c)/\Card (r^{-1}(c))$.  On the level of cochains,
the transfer map $t:L^2C^\ast(X)\to \ltwo C^\ast(\gS)$ is defined by
    \[
    t(f)(c):=\frac{1}{\Card (r^{-1}(c))}\sum f(c'),
    \]
where the sum is over all $c'\in r^{-1}(c)$.  (The orientations on the
$c'$ are induced from the orientation of $c$.)  Note that
    \(
    \Card (r^{-1}(c))=\muq (c),
    \)
where $\muq$ is the measure on $Wc$ defined in (\ref{e:muq})
(i.e., if $c=w\gs$ with $w$  $(\emptyset, S(\gs))$-reduced, then
$\muq (c)=q_w$).

\begin{Remark}
Suppose $X$ is the geometric realization of a building associated to a
Tits system $(G,B,N,S)$.  Then $\ltwo C^\ast(\gS)$ can be identified with
the $B$-invariant cochains $L^2C^\ast(X)^B$ and the map
$r:\ltwo C^\ast(\gS)\to L^2C^\ast(X)$ with the inclusion of the
$B$-invariant cochains.  The map $t:L^2C^\ast (X)\to \ltwo C^\ast(\gS)$ is
then identified with averaging over $B$.  In other words, if $\gS$ is
identified
with a subspace of $X$ via some section of $r:X\to \gS$, then
\[
t(f)(c)=\int_{x\in B} f(xc)d\mu.
\]
\end{Remark}

\begin{lemma}\label{l:randt}
\begin{enumeratei}
\item
$t\circ r=id:\ltwo C^i(\gS)\to\ltwo C^i(\gS)$.
\item
The maps  $r$ and $t$ are adjoint to each other.
\item
These maps take harmonic cocycles to harmonic cocycles.
\end{enumeratei}
\end{lemma}

\begin{proof}
Statement (i) is obvious.

(ii) For $f\in \ltwo C^i(\gS)$ and $f'\in L^2C^i(X)$, we have
\begin{align*}
\langle r(f),f'\rangle &=\sum_{c'\in X^{(i)}} [r(f)(c')][f'(c')]
=\sum_{c\in \gS^{(i)}}\sum_{c'\in r^{-1}(c)}f(r(c'))f'(c')\\
&=\sum_{c\in \gS^{(i)}}f(c)\sum_{c'\in r^{-1}(c)}f'(c')\\
&=\sum_{c\in\gS^{(i)}} \Card (r^{-1}(c))[f(c)][t(f')(c)]
=\sum_{c\in\gS^{(i)}} \muq (c)[f(c)][t(f')(c)]\\
&=\langle f,t(f')\rangle.
\end{align*}

(iii) Since $r:\ltwo C^\ast(\gS)\to L^2C^\ast(X)$
is induced by the simplicial map $r:X\to \gS$, it takes
cocycles  to cocycles. We must show it also takes cycles to cycles.
If $c'\in X^{(i-1)}$ and $d'\in X^{(i)}$ and if the incidence number
$[c':d']$ is nonzero, then it is equal to $[r(c'):r(d')]$.  Hence,
\begin{equation*}
\partial (r(f))(c')=\sum [c':d']f(r(c'))=\sum [c:d]\frac{\muq (c)}{\muq
(d)}f(c)=\partial^\bq (f)(c),
\end{equation*}
where $c=r(c')$, $d=r(d')$ and the last equality comes from the definition
given in
equation (\ref{e:partialq}).  So, $\partial^\bq (f)=0$ implies that 
$\partial (r(f))=0$.  
Since $t$ is the adjoint of $r$, it also must take cocycles to
cocycles and cycles to cycles.
\end{proof}

Consider the diagram:
    \[
    \minCDarrowwidth 0.5cm
    \begin{CD}
    \ltwo C^\ast(\gS)@>r>> L^2C^\ast(X)\\
    @VpVV       @VVPV\\
    \ltwo \ch^\ast(\gS)@>r>> \ch^\ast (X)
    \end{CD}
    \]
where $p$ and $P$ denote the projections onto harmonic cocycles.

\begin{lemma}\label{l:rp}
$P\circ r=r\circ p$.
\end{lemma}

\begin{proof}
Let $x\in \ltwo C^\ast(\gS)$.  It is enough to show that $P\circ
r(x)-r\circ p(x)$ is
orthogonal to any harmonic cocycle $h\in \ch^\ast (X)$.  We have:
$\langle P\circ r(x),h\rangle=\langle r(x),P(h)\rangle=\langle 
r(x),h\rangle$.
Hence,
\[
\langle P\circ r(x)-r\circ p(x),h\rangle=\langle r(x-p(x)),h\rangle=
\langle x-p(x),h\rangle=0,
\]
where the second and third equalities follow, respectively,
from parts (ii) and (iii) of  Lemma~\ref{l:randt}.
\end{proof}

\begin{theorem}\label{t:bettibuilding}
Suppose $\Phi$ is a building with a
chamber transitive automorphism group $G$ and with thickness vector
$\bq$.  Then the $L^2$-Betti numbers of $X$ ($=\cu (\Phi,K)$)  equal  the
$\ltwo$-Betti numbers of $\gS$, i.e.,
\[
b^i(X;G)=b^i_\bq(\gS).
\]
\end{theorem}

\begin{Remark}
This theorem is proved in \cite[Fact 3.5]{dym} in the case where the
building comes from an $BN$-pair.  Here we use Lemma~\ref{l:rp} to weaken
the hypothesis to the case of an arbitrary chamber transitive
group $G$.  The key technique of \cite{dym} of  integrating over $B$
is  replaced by the use of the transfer map $t$.
\end{Remark} 

\begin{proof}[Proof of Theorem~\ref{t:bettibuilding}]
For each simplex $\gs$ in the fundamental chamber $K$, consider the
commutative diagram:
    \[
    \minCDarrowwidth 0.5cm
    \begin{CD}
    \ltwo (W)@>r>> L^2(G/B)\\
    @VVV    @VVV\\
    \ltwo (W/W_{S(\gs)})@>r>> L^2(G/G_{S(\gs)})
    \end{CD}
    \]
where $S(\gs):=\{s\in S\mid \gs\subset K_s\}$,
where $W_{S(\gs)}$ and $G_{S(\gs)}$ are the isotropy subgroups of $\gs$ in 
$W$
and $G$, respectively, where the vertical maps are orthogonal projections 
and
where $r$ ($=r^\ast$) is the map induced by $r:G/B\to W$.  Let $e_B\in
L^2(G/B)$ denote the characteristic function of $B$ and let $e_\gs$ be its
orthogonal projection in $L^2(G/G_{S(\gs)})$.   
($e_\gs$ is the characteristic
function of $G_{S(\gs)}$ renormalized to have norm $1$.)  We note that
$e_B$ is the image of the basis vector $e_1\in L^2(W)$ under $r$ and
$e_\gs$ is the image of $a_{S(\gs)}$.  We have the commutative 
diagram:
\[
\minCDarrowwidth 0.5cm
\begin{CD}
\bigoplus \ltwo (W) @. & \mathclap{\xrightarrow{\hspace{5em}r\hspace{5em}}} & &&  
\bigoplus L^2(G)\\  
@VVV    @.  @.  @VVV\\
\bigoplus \ltwo (W/W_{S(\gs)}) @=\ltwo C^i(\gS)@>r>> L^2C^i(X)
@= \bigoplus L^2(G/G_{S(\gs)})\\
@. @VpVV    @VVPV @.\\
{}@.\ltwo \ch^i(\gS)@>r>> \ch^i(X),
\end{CD}
\]
where the sums are over all $\gs\in K^{(i)}$.  Let $\be\in\bigoplus
L^2(G/G_{S(\gs)})$
denote the vector $(e_\gs)_{\gs\in K^{(i)}}$ and let $\ba\in\bigoplus\ltwo
(W/W_{S(\gs)})$ be the vector
$(a_{S(\gs)})_{\gs\in K^{(i)}}$.  (So, $r(\ba)=\be$.) Using
Lemma~\ref{l:randt}, we get
\begin{align*}
b^i(X;G)&:=\dim_{\cN(G)}\ch^i(X)\\
&=\langle P(\be),\be\rangle =\langle Pr(\ba),r(\ba)\rangle =
\langle rp(\ba),r(\ba)\rangle\\
&=\langle p(\ba), tr(\ba)\rangle=\langle p(\ba),\ba\rangle=\dim_{\cn}\ltwo
\ch^i(\gS)\\
&:=b^i_\bq(\gS).
\end{align*}
\end{proof}

\rk{The Decomposition Theorem for $L^2(G/B)$}  As above, $G$ is a
chamber transitive automorphism group of a building $\Phi$.  For each
$T\in \cs$, let
\begin{equation*}
\wh{A}_T:=L^2(G/G_T)=L^2(G)^{G_T}
\end{equation*}
be the subspace of
$L^2(G/B)$ consisting of the square summable functions on $G$ which
are constant
on each coset  $gG_T$.   Set
    \begin{equation*}
    \wh{D}_{S-T}:=\wh{A}_T\cap \Big(\sum _{U\in
    \cs_{>T}}\wh{A}_U\Big)^\perp.
    \end{equation*}
$\wh{D}_{S-T}$ is a closed $G$-stable subspace in the regular 
representation.
(It corresponds to the $\cn$-module $D_{S-T}$ defined in
Section~\ref{s:decomp}.)

\begin{theorem}[The Decomposition Theorem for $L^2(G/B)$]\label{t:decompG}
Suppose $G$ is a
chamber transitive automorphism group of a building $\Phi$ and $B$ is the
stabilizer of a chamber.  If the thickness vector $\bq$ lies in
$\car^{-1}$, then
\[
\sum_{T\in \cs} \wh{D}_{S-T}
\]
is a dense subspace of $L^2(G/B)$ and a direct sum decomposition.
\end{theorem}

Given a module $M$ and a collection of submodules $(M_\ga)_{\ga\in 
\mathcal{A}}$,  the statement that  $(M_\ga)_{\ga\in \mathcal{A}}$ gives 
a direct sum decomposition of $M$ can be interpreted as a statement about 
chain complexes as follows.  Set
\[
C_1:=\bigoplus_{\ga\in \mathcal{A}} M_\ga \text{ and } 
C_0:=M,
\]
where $\bigoplus$ means external direct sum.  Let $\partial:C_1\to C_0$ 
be the natural map.  This gives a chain complex, $C_*:=\{C_0,C_1\}$, with 
nonzero terms only in degrees $0$ and $1$.  The statement that the 
internal sum $\sum M_\ga$ 
is direct is equivalent to the statement that $\partial$ is injective, 
i.e., that $H_*(C_*)$ vanishes in dimension $1$.  The statement that the 
$M_\ga$ span $M$ is equivalent to the statement that $\partial$ is onto, 
i.e., that $H_*(C_*)$ vanishes in dimension $0$.  Similarly, if $M$ and 
the $M_\ga$ are Hilbert spaces, then the statement that $M_\ga$ is dense 
in $M$ is equivalent to the statement that the reduced homology 
$\ch_*(C_*)$ vanishes in dimension $0$.

\begin{proof}[Proof of Theorem~\ref{t:decompG}]
The map $r$ from Lemma~\ref{l:vnsubalgebra} takes $A_T$ to $\wh{A}_T$ and 
$D_{S-T}$ to $\wh{D}_{S-T}$.  Define chain complexes $\wh{C}_*=\{\wh{C}_0,\wh{C}_1\}$ and $C_*=\{C_0,C-\}$ by 
\begin{align*}
\wh{C}_1&:=\bigoplus_{T\in\cs}\wh{D}_{S-T} \text{ and }  
&\wh{C}_0&:=L^2(G/B)\\
C_1&:=\bigoplus_{T\in\cs}D_{S-T} \text{ and } &C_0&:=\ltwo (W),
\end{align*}
where the boundary maps $\wh{C}_1\to \wh{C}_0$ and $C_1\to C_0$ are the natural 
maps.  
By the Decomposition Theorem for $\ltwo$ (Theorem~\ref{t:decomp}), 
$\ch_*(C_*)$ vanishes identically.  So, by
the proof of Theorem~\ref{t:bettibuilding}, $\ch_*(\wh{C}_*)$ has 
dimension $0$ with respect to $\cN (G)$ and hence, also vanishes 
identically.  The theorem then follows from the previous paragraph.
\end{proof}

\rk{Decoupling cohomology}  As in Section~\ref{s:spaces}, suppose we 
are given a finite CW complex $Z$
and a family of subcomplexes $(Z_s)_{s\in S}$.  
As in (\ref{e:uphi}), given a building $\Phi$,  define its
\emph{$Z$-realization} to be
    \begin{equation*}
    \cu(\Phi, Z)=(\Phi \times Z)/\sim,
    \end{equation*}
where $(\gf,x)\sim (\gf',x')$ if and only if $x=x'$ and $\gf, \gf'$ belong
to the same $S(x)$-residue.

The proof of Theorem~\ref{t:decoupled} goes through to give the following
two results.

\begin{theorem}\label{t:decoupled'}
Suppose $\Phi$ is a building with a chamber transitive automorphism group
$G$ and that its thickness vector $\bq$ lies in $\car^{-1}$.  Then there
is an isomorphism of orthogonal $G$-representations:
\[
\ch^{\ast}(\cu(\Phi,Z))\cong\bigoplus _{T\in \cs}H^{\ast}(Z,Z^{S-T})\otimes
\wh{D}_{S-T}.
\]
\end{theorem}

\begin{corollary}[{Compare \cite{dm} and \cite[Cor. 8.2 and Prop. 8.5]{dj}}]\label{t:decoupledsigma'}
Suppose $\Phi$ is a building with a chamber transitive automorphism group
$G$ and that its thickness vector $\bq$ lies in $\car^{-1}$.  Then, for
$X=\cu(\Phi,K)$, there is an isomorphism of orthogonal $G$-representations:
\[
\ch^{\ast}(X)\cong\bigoplus _{T\in \cs}H^{\ast}(K,K^{S-T})\otimes
\wh{D}_{S-T}.
\]
\end{corollary}

\section{The case where $L$ is a sphere}\label{s:sphere}
A simplicial complex $\gL$ is a \emph{generalized homology $m$-sphere}\!
(for short, a $GHS^m$)
if it is a homology $m$-manifold having the same homology
as $S^m$.  This is equivalent to the condition that, for each $T\in \cs 
(\gL)$,
$\Lk (T,\gL)$ has the same homology as $S^{m-\Card (T)}$.

Similarly, a pair$(\gL,\partial \gL)$ is a \emph{generalized homology
$m$-disk} (for short, a $GHD^m$)
if it is an acyclic homology $m$-manifold with boundary.

From now on, when we say that a complex is a generalized
homology sphere or disk or that it is a homology manifold,
\emph{we only require that it be one with respect to
homology  with real coefficients}.  (This is all that
is needed to insure that Poincar\'e duality holds for the (weighted) 
$L^2$-cohomology of various related complexes.)

If the nerve $L$ of $(W,S)$ is homeomorphic to $S^{n-1}$, then $\gS$ is a
contractible $n$-manifold.   If $L$ is PL-homeomorphic to 
$S^{n-1}$, then each face $K_T$ of the fundamental chamber $K$ is a 
PL-disk of codimension $\Card(T)$.  Similarly, if $L$ is a $GHS^{n-1}$, 
then $\gS$
is a contractible homology $n$-manifold and each $K_T$ is a contractible 
$GHD^{n-\Card(T)}$.  (See \cite{d98a,d02}.)

For the remainder of this section suppose that $L$ is a $GHS^{n-1}$.

\rk{Poincar\'e duality}  It is proved in \cite{dym} that
$\ltwo\ch_\ast(\gS)$ satisfies Poincar\'e duality, where the duality
changes $\bq$ to $\bq^{-1}$.  We repeat the argument below.

For each $T\in \cs$ and $w\in W$, the subcomplex $wK_T$ is the ``dual 
cell'' to the Coxeter cell $w\angt$ (defined in Sections~\ref{s:spaces} 
and \ref{s:weighted}).  (Strictly speaking, $wK_T$ is not a cell unless 
$\Lk(T,L)$ is a PL-sphere; however, since $(K_T,\partial K_T)$ is a 
$GHD^{n-\Card (T)}$, the $wK_T$ behave homologically as if they were dual 
cells.)  The  
chain complex obtained by partitioning $\gS$ into these ``dual cells'' 
is denoted $\ltwo C_*(\gS_{ghd})$ in \cite{dym}.   It is naturally 
identified with the cochain complex $\ltwo C^{n-*}(A_\emptyset)$ associated to the 
cosheaf $A$ on $L$, defined in Section~\ref{s:ruins}.  By 
Lemma~\ref{l:chainruin}~(ii), $\ltwo C_*(\gS_{cc})$ 
is identified with the chain complex $\ltwo C_*(H_\emptyset)$ associated 
to the cosheaf $H$ on $L$.  It is proved in \cite{dym} that the chain 
complexes $\ltwo C_*(\gS_{ghd})$ and $\ltwo C_*(\gS_{cc})$ are both chain 
homotopy equivalent to $\ltwo C_*(\gS)$, 
the chain complex defined via the  
standard simplicial structure on  $\gS$.  (This simplicial structure is a 
common subdivision of $\gS_{ghd}$ and $\gS_{cc}$.)  Hence, all three  
complexes have the same homology.  The map $\ltwo C^{n-*}(\gS_{ghd})\to 
L^2_{\bq^{-1}} C_*(\gS_{cc})$, induced by $wK_T\to w\angt $ is a chain 
isomorphism.  (When viewed as a map $\ltwo C_*(A_\emptyset) \to L^2_{\bq^{-1}} 
C_*(H_\emptyset)$, it is induced by the $j$-isomorphism of 
Section~\ref{s:h-neumann}.)  So, we have proved the following. 

\begin{proposition}[{\cite[Theorem 6.1]{dym}}]\label{p:duality}
Suppose the nerve $L$ of $(W,S)$ is a 
$GHS^{n-1}$.  Then there is $j$-equivariant isomorphism from the Hilbert 
$\cn$-module $\ltwo\ch_k(\gS)$ to the Hilbert $\cN_{\bq^{-1}}$-module 
$L^2_{\bq^{-1}}\ch_{n-k}(\gS)$ (where $j$ is the isomorphism of 
Section~\ref{s:h-neumann}).  Hence, $b^k_\bq(\gS)=b^{n-k}_{\bq^{-1}}(\gS)$.
\end{proposition}

\begin{Remark} The same type of
Poincar\'e duality (exchanging $\bq$ with $\bq^{-1}$) 
holds  for $\cu(W,Z)$, whenever $Z$ is compact and 
$\cu(W,Z)$ is a homology  manifold.  In other words, it holds provided that, 
for each $T\in \cs$, 
$(Z_T,\partial Z_T)$ is a compact homology manifold  with 
boundary (see \cite{d83,d98}).
\end{Remark}

\begin{corollary}[\cite{cd}]\label{c:reciprocity}
Suppose the nerve $L$ of $(W,S)$ is a $GHS^{n-1}$.
Then the growth series of $W$ is $(-1)^n$-reciprocal, i.e.,
    \[
    \frac{1}{W(\bq)}=\frac{(-1)^n}{W(\bq^{-1})}.
    \]
\end{corollary}

\begin{proof}
Take the alternating sums of the dimensions on both sides of the equation 
of Proposition~\ref{p:duality}.  By Proposition~\ref{p:euler}, 
the left hand side gives  $\chi_\bq(\gS)$ and
the right hand side $(-1)^n\chi_{\bq^{-1}}(\gS)$.  
\end{proof}    

The next result is proved in \cite{dym} as a corollary of 
Proposition~\ref{p:duality}.     
(It is also a consequence of Theorem~\ref{t:decoupledsigma}.)

\begin{corollary}[{\cite[Cor. 10.4]{dym}}]\label{p:GHS}\hspace{-1ex plus 1ex}
Suppose the nerve $L$ of $(W,S)$ is a $GHS^{n-1}\!$.
\begin{enumeratei}
\item
If $\bq\in \ol{\car}$, then $\ltwo\ch_\ast(\gS)$ is concentrated in 
dimension 
$0$;
moreover,
    \[
    \ltwo\ch_0(\gS)\cong A_S,
    \]
where $A_S$ is the representation of $\R$ on $\bR$ via the
symmetric character $\ga_S$ of Definition~\ref{d:characters}.
\item
If $\bq\in \ol{\car^{-1}}$, then $\ltwo\ch_\ast(\gS)$ is concentrated in 
dimension
$n$ and
    \[
    \ltwo\ch_n(\gS)\cong H_S,
    \]
where $H_S$ is the representation of $\R$ on $\bR$ via the
alternating character $\gb_S$ of Definition~\ref{d:characters}.
\end{enumeratei}
\end{corollary}

\begin{Remark}
If $L$ is a $GHS^{n-1}$, then $(K,\partial K)$ is a $GHD^n$ where $\partial K:=K^S$.  
Since $H^1(K;\zz/2)=0$, $K$ is orientable.
So, we can choose orientations for the
$n$-simplices of $K$ so that their sum is a relative cycle, $\xi_K$.  Its
homology
class $[K]\in H_n(K,\partial K)$ is the \emph{fundamental class} of $K$.
By Theorem~\ref{t:decoupledsigma}, $\ltwo\ch_n(\gS)$ is spanned by
$[K]h_S$.  (This was proved in \cite{dymara}.)
A representative for this class is obtained by taking the
fundamental cycle $\xi_K$ and then harmonizing it to $\xi_Kh_S$.
\end{Remark}

\begin{example}\label{ex:dim1}
($\dim L=1$.)  Suppose $L$ is a $k$-gon.  In other words, suppose we
are given a
Coxeter matrix on a set $S$,  so that
its nerve $L$ is a circle and so that $\Card (S)=k$.  
This means, first of all, that the $1$-skeleton
of $L$ is a $k$-gon.  When $k=3$, for $L$ to be equal
to its $1$-skeleton,  a further condition is needed.  Suppose
$S=\{s_1,s_2,s_3\}$, $m_{ij}:=m_{s_is_j}$ and $\ga_{ij}:=\pi/m_{ij}$, where
$\{i,j,k\}=\{1,2,3\}$.  The condition is that
$\ga_{12}+\ga_{23}+\ga_{13}\leq \pi$.  When this holds,
the $W$-action on $\gS$ is  isomorphic to the action of a
group of isometries on the Euclidean or hyperbolic plane generated by the
reflections across the edges of a $k$-gon.  (The Euclidean case occurs only
when $k=4$ and $W$ is right-angled or when $k=3$ and
$\ga_{12}+\ga_{23}+\ga_{13}=\pi$.)

If $\bq\in \car$, then $L^2_q \ch_\ast(\gS)$ is concentrated in dimension
$0$; if $\bq\in \car^{-1}$, it is concentrated in dimension $2$; if
$\bq\notin \car\cup\car^{-1}$, then it is concentrated in dimension $1$
(since it vanishes in dimensions $0$ and $2$).  In each case, the nonzero
Betti number is given by $\pm\chi _\bq$. 
\end{example}


\begin{corollary}\label{c:euclidean}
Suppose that $W$ is a Euclidean reflection group, i.e., that it can be
represented as a cocompact group generated by isometric reflections on
$\bR^n$.  Suppose further that $\bq\geq\bone$.  Then $\ltwo \ch_\ast(\gS)$ is
concentrated in the top dimension, $*=n$.  (It is $0$ if $\bq=\bone$.)
\end{corollary}

\begin{proof}
By Proposition~\ref{p:subexp} (or Remark~\ref{r:Euclidean}), when $t$
is a single
indeterminate, the reciprocal of the radius of convergence of $W(t)$ is
$1$.  It  follows that $\{\bq\mid \bq>\bone\}\subseteq \car^{-1}$.

Since $\gS\cong \bR^n$, it follows from \cite[Theorem B]{d98} 
that $L$ is a $GHS^{n-1}$.  In fact,
$L$ is a triangulation of  $S^{n-1}$.
(When $(W,S)$ is irreducible, $L$ is isomorphic to
boundary complex of an $n$-simplex, by \cite[Prop. 8, p. 90]{bourbaki};
when it is not irreducible it is a join of such complexes.)
So, the result is a consequence of  the previous proposition.
\end{proof}

Combining this with Corollary~\ref{t:decoupledsigma'}, we get the following (known) result.

\begin{corollary}\label{c:eucbldg}
Suppose that $X$ is a  Euclidean building with a chamber transitive 
automorphism group.  
Then its reduced $L^2$-cohomology is concentrated in the top
dimension.
\end{corollary}

\rk{A generalization of the Singer Conjecture}  
Proposition~\ref{p:GHS} states that when $L$ is a $GHS^{n-1}$, $\ltwo 
\ch_\ast(\gS)$ is concentrated in dimension $0$ for $\bq\in \ol{\car}$ and in
dimension $n$ for $\bq\in \ol{\car^{-1}}$.  What about the intermediate range,
$\bq\notin \ol{\car}\cup\ol{\car^{-1}}$?  By Remark~\ref{r:H0}, in this range,
$\ltwo \ch_0(\gS)=0$ and by Poincar\'e duality, $\ltwo \ch_n(\gS)=0$.  
For
$\bq=\bone$,  $\ltwo \ch_\ast(\gS)$ is the ordinary reduced $L^2$-homology
$\ch_\ast(\gS)$.  In this case, the Singer Conjecture predicts that
$\ch_\ast(\gS)$ vanishes except in dimension $\frac{n}{2}$.
There is considerable evidence for this version of the Singer Conjecture, at 
least in the case where $(W,S)$ is right-angled.    
For example,
it holds for $n=\dim \gS\leq 4$ and, when $L$ is a 
barycentric subdivision, for  $n=6,8$.  (See \cite{do, do04}.)

This suggests that the following generalization of the Singer Conjecture
for Coxeter groups should hold for weighted $L^2$-homology.

\begin{conjecture}[The Generalized Singer Conjecture]\label{conj:singer}\ 
Suppose $L$ is a 
$GHS^{n-1}$.
If  $\bq\leq \bone$ and $k> \frac{n}{2}$, then $\ltwo\ch_k(\gS)=0$.
\end{conjecture}

By Poincar\'e duality, this  is equivalent to the conjecture that if 
$\bq\geq \bone$ and $k< \frac{n}{2}$, then $\ltwo\ch_k(\gS)=0$.

In Section~\ref{s:rightsphere} we prove Conjecture~\ref{conj:singer} (as 
Theorem~\ref{t:11.1.1}) in the case where $W$ is right-angled and $n\leq 4$.

To further simplify the discussion, suppose $\bq=q$, a single 
indeterminate.
By Corollary~\ref{c:reciprocity}, 
the roots of $\chi_q$ ($=1/W(q)$) are symmetric
about 1, i.e., if $q$ is a root, then so is $q^{-1}$.

At one point, the following scenario (which is stronger than
Conjecture~\ref{conj:singer}) seemed plausible:
\begin{enumeratea}
\item
$\chi_q$ has exactly $n$ positive real roots (counted with
multiplicity) and
\item
$L^2_q\ch_\ast(\gS)$ is always concentrated in a single dimension.  The
dimension jumps each time $q$ passes a root of $\chi_q$ and the size of
the jump is the multiplicity of the root.
\end{enumeratea}

In fact, both (a) and (b) are false.  Gal \cite{gal} has given
counterexamples to (a) in dimensions $\geq 6$.  We shall explain
why (b) is false in dimensions $n\geq 4$ in Section~\ref{s:failure} below.

\section{Properties of weighted $L^2$-homology in the right-angled case}
\label{s:rightangled}
The usual $L^2$-cohomology of $\gS$ is the case $\bq=\bone$.
In \cite{do} the first and fourth authors studied this case
when $(W,S)$ was right-angled.  
(Recall  Definition~\ref{d:ra}:  $(W,S)$ is \emph{right-angled} 
if $m_{st}=2$ or 
$\infty$ for all pairs $\{s,t\}$ of distinct elements in $S$.)  Much of 
$\cite{do}$ extends in a straightforward fashion from $\bq=\bone$ to the 
case of a general $\bq$.  The purpose of this section is to rewrite parts of 
\cite{do} in the  general case.  

If $(W,S)$ is right-angled, then its nerve $L$ is a flag complex.  (A 
simplicial   complex $\gL$ is a \emph{flag complex} if any finite set 
of vertices in $\gL$ which are pairwise connected by edges span a simplex 
of $\gL$.)  Conversely, given any finite flag complex $L$, there is a 
right-angled Coxeter group $W_L$ with nerve $L$.  (The set of generators 
$S$ for $W_L$ is the vertex set of $L$ and $m_{st}=2$ if and only if 
$\{s,t\}$ spans an edge of $L$.)  For further explanations, see 
\cite{d83,d02,do}.

In this section, as well as in Sections~\ref{s:rightsphere} and 
\ref{s:failure},
 \emph{all simplicial complexes will be flag complexes and
all subcomplexes will be full subcomplexes.}  Given a finite flag
complex $L$, let $\Sigma_L$ be the
complex on which $W_L$ acts.   As usual,
$\bq$ is an $I$-tuple of positive real numbers.
For each $i\in\nn$, we have a Hilbert
$\cn$-module, $\ltwo\ch_i(\gS_L)$.  Similarly, to each pair
$(L, A)$, we can associate the Hilbert $\cn$-module, $\ltwo\ch_i
(\Sigma_L, W_L\Sigma_A)$.

We introduce some useful notation which reflects this situation.

\rk{Notation}
\begin{align}
\H_i(L)&:= \ltwo\ch_i(\gS_L)\quad
\Hdown^i(L):= \ltwo\ch^i(\gS_L)\label{e:not1}\\
\H_i(A) &:= \ltwo\ch_i(W_L\gS_A)\label{e:not2}\\
\H_i(L, A)&:= \ltwo\ch_i(\gS_L,W_L\gS_A)\label{e:not3}\\
b^i_\bq(A)&:= \dim_{\cn}(\H_i(A))\label{e:not4}\\
b^i_\bq(L,A)& := \dim_{\cn}(\H_i(L,A))\label{e:not5}\\
\chi_\bq(L)&:= \sum (-1)^ib^i_\bq(A).\label{e:not6}
\end{align}
The notation in (\ref{e:not2}) and (\ref{e:not4}) will not lead to
confusion, since 
$\ltwo\ch_i(W_L\gS_A)$
is the induced representation from $\ltwo\ch_i(\gS_A)$ and therefore,
    \(
    b^i_\bq(W_L\gS_A)=b^i_\bq(\gS_A),
    \)
where the left hand side of this equation denotes a dimension calculated with 
respect to $\cn (W_L)$ while the right hand side is with respect to $\cn 
(W_A)$.

\rk{Basic algebraic topology}    The next theorem is a compilation of 
properties of $\H_i(L, A)$ 
which were proved  in \cite{do} for the case $\bq=\bone$.
\pagebreak
\begin{theorem}[{Compare \cite[Section 7.2]{do}}]\label{t:compile}\  
\begin{enumeratea}
\item\label{a}
\textup{(Exact sequence of the pair)}   The sequence
$$\to \H_i(A)\to\H_i (L)\to\H_i(L, A)\to$$ is weakly exact.
\item\label{b}
\textup{(Excision)}  Let $T$ be a set
of vertices of $A$ such that the open star of any vertex in $T$ is
contained in the interior of $A$.  Then
\[
\H_i(L, A)\cong \H_i(L-T, A-T).
\]
\item\label{c}
\textup{(Mayer--Vietoris sequence)}  Suppose $L=L_1\cup L_2$
and $A=L_1\cap L_2$, where $L_1$ and  $L_2$ (and therefore, $A$) are
full subcomplexes of $L$.  Then
$$\to\H_i(A)\to\H_i(L_1)\oplus\H_i(L_2)\to\H_i(L)\to$$
is weakly exact.
\item\label{d}
With $L_1$, $L_2$ and $A$ as in (c),
\[
\H_i(L, A)\cong \H_i(L_1, A)\oplus\H_i(L_2, A).
\]
\item\label{e}
\textup{(The K\"unneth Formula:  the Betti numbers of a join)}
\[
b^k_\bq(L_1\ast L_2) = \sum_{i+j=k}b^i_\bq (L_1)b^j_\bq (L_2).
\]
\item\label{f}
\textup{(Atiyah's Formula)}
\[
\chi_\bq(L)=\sum_{T\in \cs}\prod_{s\in T}\frac{-q_s}{1+q_s}
=\frac{1}{W_L(\bq)}.
\]
\item\label{g}
\textup{($0$-dimensional homology, \cite{dym})}
\[
b^0_\bq (L)=
    \begin{cases}
    0 &\text{ if $\bq\notin\car$,} \\
    \frac{1}{W(\bq)} &\text{ if $\bq \in \car$.}
    \end{cases}
\]
and
\[
b^i_\bq(L)=0 \text{ for $i>0$, $\bq\in \car$.}
\]
\item\label{h}
\textup{(Pseudomanifolds, \cite[Theorem 10.3]{dymara})}
   Suppose $L$ is a  $(n-1)$-dimensional pseudomanifold.
Then $\gS_L$ is an $n$-dimensional pseudomanifold and, since the
$1$-skeleton of $\gS_L$ is the Cayley graph of $W_L$,  each
component of the complement of codimension $2$ skeleton of $\gS_L$ is
infinite.  So, if  $\bq \notin\car^{-1}$,  then $b^n_\bq (L)=0$.
(If $\bq \in\car^{-1}$ and in
addition, $L$ is orientable and the complement of its codimension
$2$-skeleton is connected, then $b^n_\bq (L)=1/W(\bq^{-1})$.
\item\label{i}
\textup{(The empty set)}  Since $\gS_\emptyset$ is a point,
\[
b^i_\bq(\emptyset )=
    \begin{cases}
    1 &\text{ if $i=0$, } \\
    0 &\text{ if $i\not= 0$.}
    \end{cases}
\]
\item\label{j}
\textup{(A $k$-simplex)} Given a  $k$-simplex $\gs$,
$W_{\gs}\cong (\zz_2)^{k+1}$ and $\gS_\gs= [-1, 1]^{k+1}$.
Hence, for $\bq=q$, a single indeterminate:
    \[
    b^i_\bq(\gs ) =
    \begin{cases}
    \left(\frac{1}{1+q}\right)^{k+1}&\text{ if $i=0$,}\\
    0 &\text{ if $i\not= 0$.}
    \end{cases}
    \]
\item\label{k}
\textup{(The Betti numbers of a disjoint union)}
Suppose $L$ is the disjoint union of $L_1$ and $L_2$.  Then,
for $i\ge 2$,
\[
b^i_\bq(L)=b^i_\bq(L_1)+b^i_\bq(L_2).
\]
For $\bq\notin \car_{L_1}\cup 
\car_{L_2}$, 
    \[
    b^1_\bq(L)=b^1_\bq(L_1)+b^1_\bq(L_2)+1.
    \]
\end{enumeratea}
\end{theorem}

\begin{proof}
Properties (\ref{a}) through (\ref{e}) follow from general principles as
in \cite{do}.  Property (\ref{f}) is Proposition~\ref{p:euler}; (\ref{g})
is proved in Section~\ref{s:weighted} as Proposition~\ref{p:HO}; (\ref{h})
is proved in \cite{dymara} (it also follows from
Theorem~\ref{t:decoupledsigma});  properties (\ref{i}) and (\ref{j}) are 
special cases of (\ref{g}).
Property (\ref{k}) follows from (c) (the
Mayer-Vietoris sequence); the last sentence of (\ref{k}) follows  after noting that
$L_1\cap L_2=\emptyset$ has nonzero Betti number,
$b^0_\bq (\emptyset )=1$ and that, by (\ref{g})
$b^0_\bq (L_1)=b^0_\bq (L_2)=0$.
\end{proof}

In the next proposition we assume that $I$ is a singleton so that $\bq$ is 
a single parameter $q$.  
We extend some simple calculations  of \cite{do}
from  $q=1$ to the case where $q$ is arbitrary.

\begin{proposition}[{Compare \cite[Section 7.3]{do}}]\label{p:compile}
Suppose $\bq=q$, a positive real number.
\begin{enumeratea}
\item
\textup{(The Betti numbers of $k$ points)}  Let $P_k$
denote
the disjoint union of $k$ points.  If $k\ge 2$, then
\[
b^0_q (P_k) =
    \begin{cases}
    \frac{1-(k-1)q}{1+q}  &\text{ if $q<\frac{1}{k-1}$,}\\
    0              &\text{ if $q\ge\frac{1}{k-1}$.}
    \end{cases}
\]
\[
b^1_q (P_k) =
    \begin{cases}
    0   &\text{ if $q<\frac{1}{k-1}$,}\\
    \frac{(k-1)q-1}{1+q}    &\text{ if $q\ge\frac{1}{k-1}$.}
    \end{cases}
\]
In particular,
\[
b^0_q(S^0)=b^i_q(P_2)=
    \begin{cases}
    \frac{1-q}{1+q} &\text{ if $q<1$,}\\
    0       &\text{ if $q\ge 1$,}
    \end{cases}
\]
\[
b^1_q(S^0)=b^1_q (P_2) =
\begin{cases} 0   &\text{ if $q<1$, }\\
\frac{q-1}{1+q}    &\text{ if $q\ge 1$.}
\end{cases}
\]
\item
\textup{(The Betti numbers of a suspension)}  The ``suspension'' of $L$ is
defined by $SL:=S^0\ast L$.  Then
\[
b^i_q(SL)=
\begin{cases} \frac{1-q}{1+q}b^i_q(L) &\text{ if $q<1$,}\\
   \frac{q-1}{1+q}b^{i-1}_q(L) &\text{ if $q\ge 1$,}
\end{cases}
\]
for all $i$.
\item
\textup{(The boundary complex of an $n$-octahedron)} Let
\[
O_n:=\underbrace{S^0\ast \cdots \ast S^0}_n.
\]
Then
\[
b^0_q(O_n)=
    \begin{cases}
    \left(\frac{1-q}{1+q}\right)^n &\text{ if $q<1$,}\\
    0 &\text{ if $q\ge 1$,}
    \end{cases}
\]
\[
b^i_q(O_n)=0, \text{ for $1\leq i\leq n-1$ and for all $q$,}
\]
\[
b^n_q(O_n)=
    \begin{cases}
    0 &\text{ if $q\le 1$,}\\
    \left(\frac{q-1}{1+q}\right)^n  &\text{ if $q> 1$.}
    \end{cases}
\]
\item
\textup{(The Betti numbers of a cone)}
    \begin{align*}
    b^i_q(CL)&=\frac{1}{q+1}b^i_q(L),\\
    b^{i+1}_q(CL, L)&=\frac{q}{1+q}b^i_q(L).
    \end{align*}
Moreover, the sequence of the pair $(CL, L)$ breaks up into short exact
sequences:
    \[
    0\to \Hq_{i+1}(CL, L)\to\Hq_i(L)\to\Hq_i(CL)\to 0.
    \]
\end{enumeratea}
\end{proposition}

\begin{proof}  Since $\gS_{P_k}$ is $1$-dimensional, $b^i_q(P_k)=0$ for
$i>1$.  By Theorem~\ref{t:compile}(\ref{g}), $\H_i(P_k)$ is concentrated in
one dimension.
Since $\chi_q(P_k) =\frac{1-(k-1)q}{1+q}$, the calculation in (a) follows.

The calculations of Betti numbers in (b), (c) and (d) follows immediately 
from part (a) and Theorem~\ref{t:compile}(\ref{e}).   The proof of the last
sentence of (d) is similar to \cite[Lemma 7.3.3]{do}. $\Hq_i(L)$ means 
$L^2_q\ch_i(\gS_L \times \zz_2)$ ($=L^2_q\ch_i(\gS_L)\otimes 
L^2_q(\zz_2)$). 
$\Hq_i(CL)$ can be identified with the subspace $L^2_q\ch_i(\gS_L)\otimes 
A_s$, where $s$ is the generator corresponding to the cone point.
The map $\Hq_i(L)\to \Hq_i(CL)$ is then identified with orthogonal projection 
onto this subspace.    
\end{proof}

\section{$W$ is right-angled and $L$ is a sphere}\label{s:rightsphere}
In the
right-angled case, Conjecture~\ref{conj:singer} can be attacked using
the techniques of \cite{do}.  In this case, the arguments
of \cite{do} are sufficient to prove the conjecture  for
$n\leq 4$.  We give the details below.

\rk{Poincar\'e duality}
If a pair $(D, \partial D)$ of
flag complexes is a generalized homology disk, then $\Sigma_D$ is a polyhedral 
homology manifold with boundary (its boundary being 
$W_D\Sigma_{\partial D}$). Hence, it satisfies a relative version of
Poincar\'e duality.

\begin{proposition}[{Compare \cite[Section 7.4]{do}}]\label{p:GHD}\ 
\begin{enumeratei}
\item
If $L$ is a $GHS^{n-1}$, then $b^i_\bq(L)=b^{n-i}_{\bq^{-1}}(L)$.
\item
If $(D, \partial D)$ is a $GHD^{n-1}$, then $b^i_\bq
(D, \partial D)=b^{n-i}_{\bq^{-1}}(D)$.
\item
If $(D, \partial D)$ is a $GHD^{n-1}$, then the
homology and cohomology sequences of the pair $(D, \partial D)$ are
isomorphic under Poincar\'e duality in the sense that the following
diagram commutes up to sign,
\[
\minCDarrowwidth 0.48cm
\begin{CD}
@>>> \H_{i+1}(D, \partial D)  @>>> \H_i(\partial D)  @>>>  \H_i(D) @>>>
\H_i(D, \partial D)  @>>>\\ && \Big\updownarrow\cong &&
\Big\updownarrow\cong && \Big\updownarrow\cong &&
\Big\updownarrow\cong  \\ {}  @>>>  \Hqminus^{n-i-1}(D)  @>>>
\Hqminus^{n-i-1}(\partial D)  @>>>  \Hqminus^{n-i} (D, \partial D)  @>>>
\Hqminus^{n-i}(D)  @>>>
\end{CD}
\]
where the vertical isomorphisms are given by Poincar\'e duality.
\end{enumeratei}
\end{proposition}

Suppose that $L=D_1\cup D_2$ and $M=D_1\cap D_2$.
Also suppose that $L$ is a $GHS^{n-1}$ and that $(D_1, M)$ and
$(D_2, M)$ are $GHD^{n-1}$'s.  By Theorem~\ref{t:compile}(d),
$\Hqminus^i(L, M)\cong\Hqminus^i(D_1, M)\oplus\Hqminus^i(D_2, M)$.
Similarly
to Proposition~\ref{p:GHD}(iii), the homology Mayer-Vietoris sequence of
$L=D_1\cup D_2$
is isomorphic, via Poincar\'e duality, to the exact sequence of
the pair $(L, M)$ in cohomology.  In other words, the following
diagram commutes up to sign,
\[
\minCDarrowwidth 0.5cm
\begin{CD}
@>>> \H_{i+1}(L) @>>> \H_i(M) @>>> \H_i(D_1) \oplus\H_i(D_2) @>>>\\
&&\Big\updownarrow\cong   &&\Big\updownarrow\cong
&&\Big\updownarrow\cong  \\ @>>> \Hqminus^{n-i-1}(L)
@>>> \Hqminus^{n-i-1}(M) @>>> \Hqminus^{n-i}(D_1, M)
\oplus\Hqminus^{n-i}(D_2, M) @>>>\end{CD}
\]
where the first row is the Mayer-Vietoris sequence, 
the second is the exact 
sequence of the pair and the vertical isomorphisms are given by
Poincar\'e duality.  We record the special case of this where $n=2k+1$
and $i=k$ as the following lemma.

\begin{lemma}[{Compare \cite[Lemma 7.4.6]{do}}]\label{l:connecting}
With hypotheses as above, suppose $n=2k+1$.  Then
the map $i_\ast\colon \Hqq_k(M)\to\Hqq_k(L)$ induced by the inclusion is
dual  (under Poincar\'e duality) to the connecting homomorphism
$\partial_\ast\colon\H_{k +1}(L)\to\H_k(M)$ in the Mayer-Vietoris
sequence.
\end{lemma}

\begin{proof}
In this special case, the previous diagram becomes the following:
\[
\begin{CD}
\H_{k+1}(L) @>\partial_\ast>> \H_k (M) @>>>
\H_k(D_1)\oplus\H_k(D_2)\\ \Big\updownarrow\cong
&&\Big\updownarrow\cong   &&\Big\updownarrow\cong \\ \Hqminus^k(L)
@>i^\ast>> \Hqminus^k(M) @>>> \Hqminus^{k+1}(D_1,
M)\oplus\Hqminus^{k+1}(D_2, M)
\end{CD}
\]
\end{proof}

\rk{Vanishing Conjectures}
We now consider several conjectures, $\I(n)$,
$\III(n)$, $\III'(n)$
and $\V(n)$, concerning the reduced
$\ltwo$-homology of $\Sigma_L$, where $L$ is
a generalized homology sphere.
(The notation $\I(n)$, $\III(n)$, $\V(n)$, is taken from \cite{do};
the ``$n$'' refers to the dimension of
$\Sigma_L$, so that $\dim L=n-1$.)

\begin{I(n)}  If $L$ is a $GHS^{n-1}$ and $\bq\leq \bone$, then 
$b^i_\bq(L)=0$  for $i>n/2$.   
\end{I(n)}

Given $(D, \partial D )$, a generalized homology disk,  denote by $\wh{D}$ 
(or $L$)
the $GHS$ formed by gluing on $C(\partial D)$ (the cone on $\partial D$) 
to $D$ along $\partial D$.  
If $v$ denotes the cone point, then $\partial D=L_v$ (the link of
$v$ in $\wh{D}$) and $C(\partial D)=CL_v$.  Conversely, given a $GHS$, call it
$L$, and a vertex $v$, we obtain a $GHD$, with $D=L-v$ (the full subcomplex of 
$L$ spanned by the vertices $\neq v$) and with $\partial D=L_v$.

Next we consider a seemingly weaker version of $\I(2k+1)$.

\begin{III(2k+1)}  Suppose $(D, L_v)$ is a $GHD^{2k}$ and  
$\wh{D}=D\cup CL_v$ as above.   If  $\bq\le \bone$, 
then, in the Mayer-Vietoris sequence, the
map
\[
j_\ast\oplus h_\ast\colon\H_k(L_v)\to\H_k(D)\oplus\H_k (CL_v)
\]
is a monomorphism.
\end{III(2k+1)}

By Lemma~\ref{l:connecting}, $\III(2k+1 )$ is equivalent to the following.

\begin{III'(2k+1)}  Suppose $(D, L_v)$ is a $GHD^{2k}$ and $\wh{D}=D
\cup CL_v$ as above.   If $\bq\ge \bone$, then the map
$i_\ast\colon\H_k(L_v)\to \H_k(L)$, induced by the inclusion, is the
zero homomorphism.
\end{III'(2k+1)}

The following is a stronger version of $\I(n)$.

\begin{V(n)}  Suppose $L$ is a $GHS^{n-1}$ and $A$ is any full
subcomplex.
\begin{itemize}
\item If $n=2k$ is even and  $\bq\le \bone$, then $b^i_\bq(L,A)=0$
for all $i>k$.
\item If $n=2k+1$ is odd and  $\bq\le \bone$, then $b^i_\bq(A)=0$ for
all $i>k$.
\end{itemize}
\end{V(n)}

By \cite{dymara}, $\I(1)$ and $\I(2)$ hold.

Next we list some obvious implications among these
conjectures.

\begin{lemma}[{Compare \cite[Section 8]{do}}]\label{l:implications}\ 
\begin{enumeratea}
\item
$\I(2k+1)\implies \III(2k+1 )$.
\item
$\V(n)\implies \I(n)$
\item
$\V(2k)$ implies that for any full
subcomplex $A$ of $L$ (a $GHS^{2k-1}$), we have
\[
b^i_\bq(A)=0 \text{  for all $i>k$ and $\bq\le \bone$. }
\]
\end{enumeratea}
\end{lemma}

\begin{proof}
(a) is obvious:  if $\I(2k+1)$ holds, then the $\H_\ast(L)$
terms in the Mayer-Vietoris sequence all vanish, so the
map $j_\ast\oplus h_\ast$ in $\III(2k+1)$ is a weak isomorphism.

(b) If $n=2k$, take
$A=\emptyset$ to get $b^i_\bq(L)=0$ for $i>k$.  If $n=2k+1$, take $A=L$,
to get $b^i_\bq(L)=0$ for $i>k$.

(c) Assume $\V(2k)$ holds.  By (b), $b^i_\bq(L)=0$ for
$i>k$.  Hence, in the exact sequence of the pair,
\[
\H_{i+1}(L, A)\to\H_i(A)\to\H_i(L),
\]
the first and third terms vanish for all $i>k$.
\end{proof}

\begin{lemma}\label{l:IIIimpliesIII}
$\III(2k+1)\implies \III(2l+1)$ for all $l\leq k$.
\end{lemma}

\begin{proof}
The proof is the same as in \cite[8.8.1 on p. 41]{do}.  Suppose $(D,L_v)$ is a
$GHD^{2l}$, with $l<k$.  Let $A$ be the join of $k-l$ copies of an
$m$-gon, $m\ge 5$ and assign to $A$ a thickness vector $\bq=\bone$.  If
$\III(2l+1)$ fails for $D$, then $\III(2k+1)$ fails for $D*A$ (the join of 
$D$
and $A$).
\end{proof}

\rk{Inductive Arguments}
We describe the program of \cite{do} for proving Conjecture
$\V(n)$.  The idea is to use a double induction:  first, induction on
the dimension $n$ and second, depending on the parity of $n$,
induction either on the number of vertices of $A$ or on the number of
vertices in $L-A$. In this section we always assume $\bq\le \bone$.

   As in \cite{do}, we set up some notation for the
induction on the number of vertices.  Suppose $A$ and $B$ are full
subcomplexes of $L$, the vertex sets of which differ by only one
element, say $v$.  In other words, $B=A-v$, for some
$v\in {\mathcal S}^{(1)}(A)$.
Let $A_v$ and $L_v$ denote the link of $v$ in $A$ and $L$,
respectively.  Thus, $A=B\cup CA_v$ and $CA_v\cap B=A_v$.  We note
that $L_v$ is a $GHS$ of one less dimension than $L$ and that $A_v$ is
a full subcomplex of $L_v$.

\begin{lemma}[{Compare \cite[Lemma 9.2.1]{do}}]\label{l:VimpliesV}
$\V(2k-1)\implies \V(2k)$.
\end{lemma}

\begin{proof}  Suppose $\V(2k-1)$ holds.  Let $(L, A)$ be as in $\V(2k)$
and
let $B=A-v$.  Assume, by induction on the number of vertices in $L-A$,
that $\V(2k)$ holds for $(L, A)$.  (The case $A=L$ being trivial.)  We
want to prove it also holds for $(L, B)$, i.e., that $b^i_\bq(L, B)=0$ for
$i>k$.  Consider the exact sequence of the triple $(L, A, B)$:
\[
\to\H_i(A, B)\to\H_i(L, B)\to\H_i(L, A)\to.
\]
Suppose $i>k$.  By
inductive hypothesis, $b^i_\bq(L, A)=0$.  By excision
(Theorem~\ref{t:compile}(b)),
$b^i_\bq(A, B)=b^i_\bq(CA_v, A_v)$.  By Theorem~\ref{t:compile}(e),
\[
b^i_\bq (CA_v, A_v)=\frac{q_v}{q_v+1}b^{i-1}_\bq(A_v).
\]
Since $\V(2k-1)$ holds for
$(L_v, A_v)$ and since $i-1>k-1$, $b^{i-1}_\bq(A_v)=0$.  So,
$0=b^i_\bq(CA_v, A_v) =b^i_\bq(A, B)$.  Consequently, $b^i_\bq(L, B)=0$.
\end{proof}

Essentially the same argument proves the following lemma.

\begin{lemma}[{Compare \cite[Lemma 9.2.2]{do}}]\label{l:Vmfld}
Assume that $\V(2k)$ holds.  Suppose that a flag
complex $L$ is a polyhedral homology manifold of dimension $2k$ and
that $A$ is a full subcomplex.  Then $b^i_\bq(L, A)=0$ for $i>k+1$and
$\bq\le \bone$.
\end{lemma}

\begin{proof}  We proceed as in the previous proof.  If $B=A-v$, then
\[
b^i_\bq(A, B)=b^i_\bq (CA_v, A_v)=\frac{q_v}{q_v+1}b^{i-1}_\bq (A_v).
\]
Since
we are assuming $\V(2k)$ holds, Lemma~\ref{l:implications}(c) implies that
$b^{i-1}_\bq (A_v)=0$ for $i>k+1$.  Hence, if we assume by induction
that the lemma holds for $(L, A)$, then it also holds for $(L, B)$.
\end{proof}

\begin{lemma}[{Compare \cite[Lemma 9.2.3]{do}}]\label{l:IIIimpliesV}
$\left[\V(2k)\text{ and }\III(2k+1)\right]\implies \V(2k+1).$
\end{lemma}

\begin{proof}
   Assume $\V(2k)$ and $\III(2k+1)$ hold.  Let $(L, A)$ be as in $\V(2k
+1)$ and let $B=A-v$.  Assume, by induction on the number of vertices
in $B$, that $\V(2k+1)$ holds for $B$.  (The case $B=\emptyset$ being
trivial.)  We want to prove that it also holds for $A$, i.e., that
$b^i_\bq(A)=0$ for $i>k$.

First suppose that $i>k+1$.  Consider the Mayer-Vietoris sequence for
$A=B\cup CA_v$:
\[
\H_i(B)\oplus\H_i(CA_v)\to\H_i(A)\to\H_{i-1} (A_v).
\]
By $\V(2k)$ and Lemma~\ref{l:implications}(c), $b^{i-1}_\bq(A_v)=0$ (since
$i-1>k$) and hence,
$b^i_\bq(CA_v)=0$ (by Theorem~\ref{t:compile}(c)).  By inductive 
hypothesis,
$b^i_\bq(B)=0$, and consequently, $b^i_\bq(A)=0$.

For $i=k+1$, we compare the Mayer-Vietoris sequence of $A=B\cup CA_v$
with that of $L=D\cup CL_v$ (where $D=L-v$):
\[
\begin{CD}
& & & &\H_{k+1}(L_v, A_v)\\ & & & & @VVV\\ 0 @>>>
\H_{k+1}(A) @>>> \H_k(A_v) @>{j'_\ast \oplus h'_\ast}>>
\H_k(B)\oplus\H_k(CA_v)\\ & & & & @V{f_\ast}VV  @VVV \\ & & & &
\H_k(L_v) @>{j_\ast\oplus  h_\ast}>> \H_k (D)\oplus
\H_k(CL_v)
\end{CD}
\]
By $\V(2k)$, $b^{k+1}_\bq(L_v, A_v)=0$; hence,
$f_\ast$ is injective.  By $\III(2k+1)$, $j_\ast\oplus h_\ast$ is
injective.  Hence,  $j'_\ast\oplus h'_\ast$ is injective and therefore,
$b^{k+1}_\bq(A)=0$.
\end{proof}

One of the main results of \cite{do} has the
following analog.

\begin{theorem}[{Compare \cite[Theorem 9.3.1]{do}}]\label{t:maindo}
Statement $\III(2k-1)$ implies that $\V(n)$ holds
for all $n \leq 2k$.
\end{theorem}

\begin{proof}
By Lemma~\ref{l:IIIimpliesIII}, $\III(2k-1)$ implies $\III(2l -1)$,
for all $l\le
k$.  Suppose, by induction on $n$, that $\V(n-1)$ holds for some $n\le
2k$.  If $n-1$ is odd, then by Lemma~\ref{l:VimpliesV}, $\V(n-1)$ implies
$\V(n)$.  If $n-1$ is even, then by Lemma~\ref{l:IIIimpliesV}, $\V(n-1)$ 
and
$\III(n)$ imply $\V(n)$.
\end{proof}

\rk{The conjecture in dimension $3$}
We begin with a discussion of triangulations of $S^2$. (Details can be
found in \cite[Section 10.2]{do}.)

For $j=1, 2$, suppose that $L_j$ is a flag
triangulation of $S^2$ and that $s_j$ is a vertex of valence $4$ in
$L_j$.  Choose an identification of the link of $s_1$ with that of
$s_2$.  (They are both $4$-gons.)  Define a new triangulation
$L_1\square L_2$ of $S^2$ by gluing together the $2$-disks $L_1-s_1$
and $L_2-s_2$ along their boundaries.

Conversely, suppose $C$ is an empty $4$-circuit in
$L$.  Then $C$ separates $L$ into two $2$-disks, $D_1$ and $D_2$.  Let
$L_1$ and $L_2$ denote the result of capping off $D_1$ and $D_2$,
respectively (where ``capping off'' means adjoining a cone on the
boundary).  Then $L=L_1\square L_2$.

\begin{lemma}[{Compare \cite[Lemma 10.2.7]{do}}]\label{l:10.2.7}
For $\bq\le \bone$, $b^2_\bq (L_1\square
L_2)=b^2_\bq (L_1)+b^2_\bq(L_2)$.
\end{lemma}

\begin{proof}  This follows from the Mayer-Vietoris sequence and
Proposition~\ref{p:compile}(c).
\end{proof}

Andreev \cite{andreev1, andreev2} determined the possible 
fundamental
polytopes for any reflection group on $\bH^3$ of cofinite volume. 
The right-angled case of the Andreev's Theorem is the following.

\begin{theorem}[Andreev's Theorem]\label{t:10.3.1}
Suppose that $L$ is a flag triangulation of $S^2$ and that
\begin{description}
\item{(i)} $L$ has no empty $4$-circuits, and

\item{(ii)} $L$ is not the suspension of a $4$- or $5$-gon.
\end{description}
Let $T$  denote the set of valence $4$ vertices of $L$ and let
$K_{[S-T]}$  be the dual of the cellulation $[S-T]$ of $S^2$ obtained by
replacing stars of  vertices of $T$ by square cells.
Then $K_{[S-T]}$ can be realized as an ideal, right-angled convex
polytope in $\bH^3$.  (The ideal vertices correspond to the square faces
of $[S-T]$,  i.e., to the vertices of valence $4$ in $S$.)  The
resulting hyperbolic reflection group is the right-angled Coxeter group
$W_{S-T}$.
\end{theorem}

Next we show that $\III'(3)$ is true for right-angled reflection groups on
$\bH^3$.

\rk{Equidistant hypersurfaces}  
Suppose a Coxeter group $W$ acts by reflections on
hyperbolic $(2k+1)$-space
$\bH^{2k+1}$ with a fundamental polytope $K$ of finite volume.

Let $\bH^{2k}$ be a wall.  We claim  that the
map $\ltwo\ch_k(\bH^{2k}) \to \ltwo\ch_k(\bH^{2k+1})$, induced by 
inclusion, is
the zero map for $\bq\ge \bone$.

Our argument  uses weighted $L^2$-de~Rham cohomology
theory. We
will show that the map $\ltwo\ch^k(\bH^{2k+1}) \to
\ltwo\ch^k(\bH^{2k})$, induced
by restriction of forms, is the zero map.  To define these terms we first 
need
a ``weight function'' on $\bH^{2k+1}$ which we can then use to define a
new inner product on the vector space  $C^{\infty}$ $j$-forms on
$\bH^{2k+1}$.

Given any measurable nonnegative function $f:H^{2k+1}\to [0,\infty)$, one 
can
modify the volume form on $\bH^{2k+1}$ by multiplying by $f$ and then
define a new norm on $C^\infty$ $j$-forms $\go$ by
\[
\normf{\go}^2=\int_{\bH^{2k+1}} \norm{\omega}^2_p\, f(p)\,dV,
\]
where $\norm{\go}_p$ denotes the pointwise norm.
$\normf{\go}$ is called the \emph{$L^2_f$-norm} of $\go$.

Let $K$ be a fundamental polytope for $W$ on $\bH^{2k+1}$.
As usual, $\bq$ is an $I$-tuple of positive real numbers.  For
any point $p$ in $\bH^{2k+1}$, put $f(p)=q_w$ when $p\in wK$. Of course,
this expression is ambiguous for $p\in w\,\partial K$.  Nevertheless,
choose some
convention to remove the ambiguity, for example,
that $w$ is the element of minimum word length with $p\in wK$.
Then $f$ is the \emph{word length weight function} on $\bH^{2k+1}$.
It is a sort of
step function in that it is constant on the interior of each chamber.

When $K$ is compact, the arguments of \cite{dod77} go through to show
that  the cellular weighted $L^2$-cohomology of $\gS$ can be calculated 
using
weighted de Rham cohomology, i.e.,
\[
\ltwo \ch^\ast(\gS)\cong \ltwo\ch^k(\bH^{2k+1}),
\]
where the right hand side is defined using $L^2_f$ forms with $f$ the
word length weight function defined above.  When $K$ is not compact but
has finite volume we can reach the same conclusion by using \cite{cg85}.

Next let $\bH^{2k}$ be a supporting wall of $K$ (i.e., $\bH^{2k}$ is a wall
determined by a codimension one face of $K$).  Put coordinates $(x,y)$
on $\bH^{2k+1}$ by letting $y\in \bR$ be the oriented distance from $p$
to the nearest point $x\in \bH^{2k}$.
Let $N_y$ be the hypersurface in $\bH^{2k+1}$ consisting of the points
of (oriented) distance $y$ from $\bH^{2k}$. Let $p_y\colon N_y\to
\bH^{2k}$ be the projection which takes a point in $N_y$ to the
closest point in $\bH^{2k}$. Then $p_y$ is a homothety. Let
$\phi_y\colon \bH^{2k} \to N_y$ be its inverse. Also, let  $i\colon
\bH^{2k} \to \bH^{2k+1}$ and  $i_y\colon N_y\to \bH^{2k+1}$ be the
inclusions. Thus, $i$ and $i_y\circ\phi_y$ are properly homotopic.

Let $g(x,y)=f(x,0)$. Note that $f(x,y)\ge g(x,y)$.

Let $\omega$ be a closed $L^2_f$-$k$-form on $\bH^{2k+1}$.  We claim
that the restriction $i^*(\omega)$ of $\omega$ to $\bH^{2k}$
represents the zero class in reduced  $L^2_f$-cohomology.
Suppose, to the contrary, that
$[i^*(\omega)]\neq 0$. Then   $\normg{i^*(\omega)}\geq
\normg{[i^*(\omega)]}\geq 0$, where  $\normg{[i^*(\omega)]}$ denotes the
norm of the harmonic representative of the cohomology class
$[i^*(\omega)]$.  Since $\phi_y$ is a conformal diffeomorphism, it
follows  that it preserves norms of middle-dimensional forms:
$\normg{\phi^*_y(i^*_y(\omega)}=\normg{i^*_y(\omega)}$.  Since $i$ and
$i_y\circ\phi_y$ are properly homotopic, $[\phi^*_y(i^*_y(\omega)]=
[i^*(\omega)]$, so it follows that $\normg{i^*_y(\omega)}\ge
\normg{[i^*(\omega)]}$.  Now, since $i^*_y(\omega)$ is just a
restriction of $\omega$, we have a pointwise inequality
$\norm{\omega}_x \ge \norm{i^*_y(\omega)}_x$. Therefore, using
Fubini's Theorem, we obtain
\begin{multline*}
\normg{\omega}^2=\int_{\bH^{2k+1}} \norm{\omega}^2_x\, g(x,y)\,dV =
\int_\bR \int_{N_y} \norm{\omega}^2_x \, g(x,y)\,dA \, dy \ge \\
\int_\bR \int_{N_y} \norm{i^*_y(\omega)}^2_x \,g(x,y)\, dA \, dy  =
\int_\bR \normg{i^*_y(\omega)}^2 \, dy \ge
\int_\bR   \normg{[i^*(\omega)]}^2 \, ds =
\infty.
\end{multline*}
Since $\normf{\omega}\ge\normg{\omega}$,
this contradicts our assumption that the $L^2_f$-norm of $\omega$ is finite 
and thereby completes the proof.

In dimension $3$ we get the following.

\begin{theorem}\label{t:10.1.3}
Suppose that $L$ is a flag
triangulation of $S^2$ satisfying the conditions of the Andreev's
Theorem. Then $\III'(3)$ is true for $L$.
\end{theorem}

\begin{proof}
If $v\in T$, then, by Proposition~\ref{p:compile}(c), $b^1_\bq(L_v)=0$
so $\III'(3)$ is automatic.  If
$v\not\in T$, then the result follows from the Andreev's Theorem
and the previous paragraphs.
\end{proof}

\begin{theorem}\label{t:10.4.1}
$\I(3)$ is true:  if $L$ is a triangulation of the $2$-sphere
as a flag complex, then
\[
b^i_\bq (L)=0  \text{  for  $i>1$ and $\bq\le \bone$. }
\]
\end{theorem}

\begin{proof}  If $L$ is the suspension of a $4$- or $5$-gon, then the
theorem
follows from Proposition~\ref{p:compile}(b).  If $L$ is not the suspension 
of a
$4$-gon or a $5$-gon and if it has no empty $4$-circuits,
then the theorem follows from
Theorem~\ref{t:10.1.3}, Lemma~\ref{l:IIIimpliesV} and the fact that $I(2)$
holds (\cite{dymara}).

In every other case, $L$ has an empty $4$-circuit which we can use to
decompose it as, $L=L_1\square L_2$, as before.  Since $L_1$
and $L_2$ each have fewer vertices than does $L$, this process must
eventually terminate.  So, the theorem follows from Lemma~\ref{l:10.2.7}.
\end{proof}

Since $\I(3 )$ is
true, Theorem~\ref{t:maindo} (together with
Lemma~\ref{l:implications}(a)) yields
the  following.

\begin{theorem}[{Compare \cite[Theorem 11.1.1]{do}}]\label{t:11.1.1}
$\V(n)$ is true for $n\le 4$.
\end{theorem}

If $L$ is a flag triangulation of $S^3$, then $\V(4)$, Poincar\'e duality
and \cite{dym} imply:
\begin{description}
\item for $\bq\in \ol{\car}$, $\H_\ast(L)$ is concentrated in dimension $0$,
\item for $\bq\le 1$ and $\bq\notin \car$,
$\H_\ast(L)$ is concentrated in dimensions $1$ and $2$,
\item for $\bq > 1$ and $\bq\notin \car^{-1}$, $\H_\ast(L)$ is
concentrated in dimensions $2$ and $3$,
\item for $\bq\in \ol{\car^{-1}}$, $\H_\ast(L)$ is concentrated in dimension 
$4$.
\end{description}

\section{\!Failure of concentration in the intermediate
range}\label{s:failure}
In this section $I$ is a singleton (so that $q$ is a single
parameter)  and $W$ is right-angled.   We retain the notation and conventions 
of Section~\ref{s:rightangled}.  

\rk{The $h$-polynomial} Combinatorialists have
associated two polynomials to a finite simplicial complex $L$: its
``$f$-polynomial,'' $f_L(t)$, and its ``$h$-polynomial,''
$h_L(t)$.  The first is defined by
    \begin{equation}\label{e:fpoly}
         f_L(t):=\sum_{T\in \cs(L)}t^{\Card(T)}=\sum_{i=0}^n f_{i-1}t^i,
    \end{equation}
where $f_m$ is the number of $m$-simplices of $L$, $f_{-1}=1$ and $\dim
L=n-1$.  The second one is defined by
    \begin{equation}\label{e:hpoly}
    h_L(t):=(1-t)^nf_L\left(\frac{t}{1-t}\right).
    \end{equation}

If a Coxeter system $(W,S)$ is right-angled,  then for each
spherical subset $T$, $W_T\cong (\zz/2)^{\Card(T)}$. So,
$W_T(t)=(1+t)^{\Card (T)}$.  Hence,
\begin{equation}\label{e:rightangled}
\frac{1}{W_T(t)}=\left(\frac{1}{1+t}\right)^{\Card
(T)}\text{ and }
\frac{1}{W_T(t^{-1})}=\left(\frac{t}{1+t}\right)^{\Card (T)}.
\end{equation}

\begin{proposition}\label{p:hpoly}
Suppose $(W,S)$ is a right-angled Coxeter system and that its
nerve $L$ is $(n-1)$-dimensional.  Then
\[
\frac{1}{W(t)}=\frac{h_L(-t)}{(1+t)^n}.
\]
\end{proposition}

\begin{proof}
By Lemma~\ref{l:growth}~(\ref{c:Wt-1}) and (\ref{e:rightangled}),
\begin{align*}
\frac{1}{W(t)}&=\sum_{T\in\cs}\frac{\geps (T)}{W_T(t^{-1})}=
\sum_{T\in\cs}\left(\frac{-t}{1+t}\right)^{\Card (T)}\\
&=f_L\left(\frac{-t}{1+t}\right)=\frac{h_L(-t)}{(1+t)^n}.
\end{align*}
\end{proof}

In the next proposition we record some properties of  $h_L(t)$.

\begin{proposition}\label{p:hpolyfacts}
Suppose $L$ is a $GHS^{n-1}$.  Let $h_L(t)=\sum h_it^i$ be its
$h$-polynomial.  Then
\begin{enumeratei}
\item
$h_L$ is a polynomial of degree $n$.  The constant term $h_0$ is 1.
\item
$h_L(t)=t^nh_L(t^{-1})$.  (This means that the coefficient sequence
$(h_0,\dots,h_n)$ is palindromic.  It also implies that $t\to t^{-1}$ is a
symmetry of the set of  roots of $h_L$.)
\item
Each $h_i\geq  0$.
\item
If $L$ is also assumed to be $3$-dimensional and a flag complex, then all
four roots of $h_L(t)$ are real.
\end{enumeratei}
\end{proposition}

Statements (i),(ii) and (iii) are  well-known; proofs can be found in
\cite{bronsted}.  Statement (iv) is proved in \cite{boros}.  We give a
simple argument for it below.

\begin{proof}[Proof of (iv)]
Put $\hat{h}(t)=h_L(-t)$.  By Proposition~\ref{p:hpoly},
$1/W_L(t)=\hat{h}(t)/(1+t)^n$.  
The Flag Complex Conjecture is that for $n-1=2k-1$, 
$(-1)^k/W_L(1)\geq 0$,
i.e., $(-1)^k\hat{h}(1)\geq 0$.  (See \cite{do}, \cite{d02}.) Let $\rho$
be the radius of convergence of $W_L(t)$.  By Lemma~\ref{l:smallroot},
$\rho$ is a root of $\hat{h}$ and it is the smallest root in absolute
value.  By (ii), $\rho^{-1}$ is also a root of $\hat{h}$ and it is the
largest in absolute value.

Now suppose $\dim L=3$.    To prove (iv), it suffices to show the four
roots of $\hat{h}$ are positive reals.
The Flag Complex Conjecture is known to hold in
this dimension (by \cite{do}), i.e., $\hat{h}(1)\geq 0$.  We know that 
$\rho$
and $\rho^{-1}$ are roots and also
that $\hat{h}(t)>0$ for $t\in [0,\rho)$ or $t\in (\rho^{-1},\infty)$.  If
the other two roots of $\hat{h}$ don't lie in $[\rho,\rho^{-1}]$, then
$\hat{h}$ must be negative on that interval, contradicting the fact that
$\hat{h}(1)\geq 0$.
\end{proof}

For any full subcomplex $A$ of $L$,  set
    \(
    r_A:=\rho_A^{-1},
    \)
where, as before, $\rho_A$ is the radius of convergence of $W_A(t)$.  Since
$W_A$ is a subgroup of $W_L$, $\rho_A\geq\rho_L$; hence, $r_A\leq r_L$.

Next, suppose that $M$ is a  $GHS^{n-2}$ and a full subcomplex of $L$  (so,
$M$ is a homology submanifold of codimension one in $L$).  Then $M$
separates $L$ into two generalized homology $(n-1)$-disks, say, $A$ and
$B$.  Thus, $\partial A=\partial B =M$ and $L=A\cup B$.  Let $CM$ denote
the cone on $M$.   Let  $\wh{A}$ (resp. $\wh{B}$) denote the result of
gluing $CM$ onto $A$ (resp. $B$) along $M$.

\begin{lemma}\label{l:estimate}
With hypotheses as above, suppose $q<\min\{r_L,r_{\wh{A}},r_{\wh{B}}\}$ and
$q>r_M$.  Then
\[
b^{n-1}_q(L)\geq \frac{q-1}{q+1}b^{n-1}_q(M)>0.
\]
\end{lemma}

\begin{proof}
Since $q>1$, by Proposition~\ref{p:compile}(d), we have
\begin{equation}\label{e:estimate1}
b^k_q(CM)=b^0_q(\mathrm{point})b^k_q(M)=\frac{1}{1+q}b^k_q(M).
\end{equation}
By Remark~\ref{r:H0}, since $q<r_L$, $\Hq_n(L)\cong\Hq_0(L)=0$.  By
Proposition~\ref{p:GHS}, since $q>r_M$, $\Hq_k(M)=0$ for $k\neq n-1$.
Hence, the Mayer-Vietoris sequence (Theorem~\ref{t:compile}(c))
for $L=A\cup B$ gives a weakly exact
sequence:
\[
0\mapright{}\ \Hq_{n-1}(M)\mapright{}\
\Hq_{n-1}(A)\oplus\Hq_{n-1}(B)\mapright{}\ \Hq_{n-1}(L)\mapright{} \ 0.
\]
So,
    \begin{equation}\label{e:estimate2}
    b^{n-1}_q(L)=b^{n-1}_q(A)+b^{n-1}_q(B)-b^{n-1}_q(M).
    \end{equation}
A similar Mayer-Vietoris sequence for $\wh{A}=A\cup CM$ gives
    \[
    b^{n-1}_q(\wh{A})=b^{n-1}_q(A)+b^{n-1}_q(CM)-b^{n-1}_q(M),
    \]
which we rewrite as
    \begin{align}
    b^{n-1}_q(A)&=b^{n-1}_q(\wh{A})-b^{n-1}_q(CM)+b^{n-1}_q(M)\notag\\
    &=b^{n-1}_q(\wh{A})+\frac{q}{1+q} b^{n-1}_q(M)\label{e:estimate3},
    \end{align}
where the second equality is from (\ref{e:estimate1}).  Similarly,
    \begin{equation}\label{e:estimate4}
    b^{n-1}_q(B)=b^{n-1}_q(\wh{B})+\frac{q}{1+q} b^{n-1}_q(M),
    \end{equation}
Combining (\ref{e:estimate2}), (\ref{e:estimate3}) and
(\ref{e:estimate4}), we get
\begin{align*}
b^{n-1}_q(L)&=b^{n-1}_q(\wh{A})+\frac{q}{1+q}b^{n-1}_q(M)+b^{n-1}_q(\wh{B})
+\frac{q}{1+q}b^{n-1}_q(M)-b^{n-1}_q(M)\\
&=b^{n-1}_q(\wh{A})+b^{n-1}_q(\wh{B})+\frac{q-1}{1+q}b^{n-1}_q(M)\\
&\geq\frac{q-1}{1+q}b^{n-1}_q(M)>0,
\end{align*}
where the last inequality holds because $q>1$ and $b^{n-1}_q(M)>0$ (since
$q>r_M$).
\end{proof}

\begin{lemma}[Failure of concentration in dimension 4]\label{l:concentration}
Suppose that
$L$ is a triangulation of $S^3$ as a flag complex, that a full
subcomplex $M$ is isomorphic to the boundary complex of an
octahedron and that $M$ divides $L$ into two $3$-disks $A$ and $B$
nontrivially, i.e., neither $A$ nor $B$ is a cone on $M$.
Suppose further that $\chi_1(L)\neq 0$. Let $p$ be the second
largest root of $h_L(-t)$, and let $r=\min\{p,
r_{\wh{A}},r_{\wh{B}}\}$. Then $r>1$ and for $1<q<r$, $b^2_q (L)$
and $b^3_q (L)$ are both nonzero.
\end{lemma}
\begin{proof}
We want to use Lemma~\ref{l:estimate} for $n=4$.
Since $W_M$ is the product of 3 copies of the infinite dihedral group,
its growth series is given by
    \begin{equation*}
    W_M(t)=\left(\frac{1+t}{1-t}\right)^3.
    \end{equation*}
So, $\rho_M=1=r_M$.

Suppose $r_{\wh{A}}=1$.  Then $r_A=1$ and by
Proposition~\ref{p:subexp}, $W_A$ splits as $W_0\times W_1$, where
$W_0$ is a Euclidean reflection group and $W_1$ is finite.  Since
$M=\partial A$, the only possibility is $W_0=W_M$ and $W_1=\bZ/2$,
i.e., $A=CM$, which we have ruled out by hypothesis.   Similarly,
for $B$. Thus, $\min\{r_{\wh{A}},r_{\wh{B}}\}>1$. Since
$\chi_1(L)\neq 0$, $p\neq 1$ and by \cite{do}, $\chi_1(L)
>0$. So, $\chi_q(L)$ is positive on the interval $(p^{-1},p)$ and
therefore, also on the subinterval $(1,r)$.

By Lemma~\ref{l:estimate}, for $1<q<\min\{r_L,r_{\wh{A}},r_{\wh{B}}\}$,
$b^3_q(L)>0$.  On the interval $(p^{-1},p)$ we  have
$b^4_q(L)=0=b^0_q(L)$ as well as $\chi_q(L)>0$ and this
forces $b^2_q(L)>0$.  Therefore, for $1<q<r$, $\Hq_\ast (L)$ is
nonzero in dimensions 2 and 3.
\end{proof}

\begin{example}\label{ex:existence}
(\emph{Existence}).  
Here we show that there is a flag triangulation $L$ of $S^3$ together with a full
subcomplex $M\subset L$ so that the conditions  of
Lemma~\ref{l:concentration} are satisfied.  
Let $P_m$ denote an $m$-gon (i.e., a
triangulation of $S^1$ with $m$ vertices).  Let $I_l$ denote the
triangulation of an interval with $l$ vertices.  Let $A_{k,m}$ denote a
triangulation of the annulus $S^1\times [0,1]$ such that its two boundary
components are $P_k$ and $P_m$ and such that there are no interior
vertices. (This does not determine the triangulation, but it does determine
the number of $i$-simplices in $A_{k,m}$ for $i=0,1,2$.) Form the
suspension $SA_{k,m}:=S^0\ast A_{k,m}$.  It has
two boundary components:  $SP_k$ and $SP_m$.
Fill in $SP_m$ with $I_4\ast P_m$ to get a triangulation $A$ of $D^3$,
i.e.,
\[
A:=SA_{k,m}\cup_{SP_m} (I_4\ast P_m).
\]
If $k=4$, then $\partial A=SP_4$, which is  the boundary complex $M$ of an
octahedron.  Hence, we can double $A$ along its boundary to get a
triangulation $L$ of $S^3$ (so, $B=A$).

By Theorem~\ref{t:compile}~(f),
    \begin{equation*}
    \chi_q(L)=\frac{1}{W(q)}=f_L\left(\frac{-q}{1+q}\right),
    \end{equation*}
where $f_L$ was defined in (\ref{e:fpoly}).
This formula is the basic method used for computing Euler
characteristics.  It gives
    \begin{align*}
    \chi_q(P_m)&=1-\frac{mq}{(1+q)}+\frac{mq^2}{(1+q)^2}=
    \frac{1-(m-2)q+q^2}{(1+q)^2},\\
    \chi_q(I_4)&=1-\frac{4q}{(1+q)}+\frac{3q^2}{(1+q)^2}=
    \frac{1-2q}{(1+q)^2}.\\
    \end{align*}
We compute the number of
simplices in $A_{k,m}$.   Each triangle of $A_{k,m}$ has exactly one of its
edges on the boundary and each interior edge is on the boundary of two
triangles.  Hence, there are $k+m$ triangles in $A_{k,m}$ and $k+m$
interior edges.
So, $f_0(A_{k,m})=k+m$, $f_1(A_{k,m})=2(k+m)$,
$f_2(A_{k,m})=k+m$ and
\begin{align*}
\chi_q(A_{k,m})&=1-\frac{(k+m)q}{(1+q)}+\frac{2(k+m)q^2}{(1+q)^2}
-\frac{(k+m)q^3}{(1+q)^3}\\
&=(1-(k+m-3)q+3q^2+q^3)/(1+q)^3.
\end{align*}
Therefore,
\begin{align*}
\chi_q(SA_{k,m})&=\chi_q(S^0)\chi_q(A_{k,m})=
\frac{1-q}{1+q}\chi_q(A_{k,m})\\
&=(1-(k+m-2)q+(k+m)q^2-2q^3-q^4)/(1+q)^4,\\
\chi_q(I_4\ast P_m)&=\chi_q(I_4)\chi_q(P_m)\\
&=(1-mq+(2m-3)q^2-2q^3)/(1+q)^4,\\
\chi_q(SP_m)&=\chi_q(S^0)\chi_q(P_m)\\
&=(1-(m-2)q+(m-2)q^3-q^4)/(1+q)^4
\end{align*}
So,
\begin{align*}
\chi_q(A)&=\chi_q(SA_{k,m})+\chi_q(I_4\ast P_m)-\chi_q(SP_m)\\
&=(1-(k+m)q+(k+3m-3)q^2-(m+2)q^3)/(1+q)^4
\end{align*}
Taking $k=4$, $\chi_q(A)=(1-(m+4)q+(3m+1)q^2-(m+2)q^3)/(1+q)^4$; hence,
\begin{align*}
\chi_q(\wh{A})&=\chi_q(A)-\left(\frac{q}{1+q}\right)\chi_q(M)=
\chi_q(A)-\left(\frac{q}{1+q}\right)\left(\frac{1-q}{1+q}\right)^3\\
&=(1-(m+5)q+(3m+4)q^2-(m+5)q^3+q^4)/(1+q)^4.
\end{align*}
When $m=10$, the numerator is
\[
h_{\wh{A}}(-q)=1-15q+34q^2-15q^3+q^4,
\]
which has roots $.08$, $.48$, $2.10$ and $12.34$ (rounded off to two 
decimal places).  Similarly,
\begin{equation*}
\chi_q(L)=2\chi_q(A)-\chi_q(M)
=(1 - 26q + 62q^2 - 26q^3+q^4 )/(1+q)^4,
\end{equation*}
which has roots
 $.04$, $.48$, $2.08$ and $23.40$.
So, the numbers in Lemma~\ref{l:concentration} are
$r_{\wh{A}}=r_{\wh{B}}=12.34$ and $r=p=2.08$.  In particular,
since $r>2$, the right-angled building with $q=2$ has nonvanishing
$L^2$-homology in dimensions 2 and 3.
\end{example}

\section{Remarks about other groups}\label{s:final}
Suppose $\gG$ is a countable discrete group and $|\ |$ is a ``norm'' on 
it, i.e., $|\ |$ is a function from $\gG$ to $[0,\infty)$ such that 
$|\ga\gb|\leq |\ga|+|\gb|$.  For example, $|\ |$ might be defined by 
$|\gamma|=l(\gamma)$ where $l:\gG\to \bZ$ is word length with respect to a 
finite set of generators $S$.  Suppose further that $\gG$
acts properly and cellularly on a CW complex $X$
and that a subcomplex $D\subset X$ is a 
``fundamental domain'' in the sense that it contains at least one 
cell from  each $\gG$-orbit of cells.  
Given a cell $\gs\subset X$, define $d(\gs)$, its \emph{distance 
from $D$}, by $d(\gs):=\min \{l(\gamma)\mid \gs\subset \gamma D\}$.

As before, given a positive real number $q$, define an inner product 
$\langle\ ,\ \rangle_q$ on 
$\bR^{(\gG)}$ by $\langle e_\gamma,e_{\gamma'}\rangle_q:=q^{|\gamma|} 
\gd_{\gamma\gamma'}$.  Let $L^2_q(\gG,|\ |)$ be its completion.  
Similarly, define an inner product on    
compactly supported cellular $i$-cochains, $C^i_c(X)$, by 
$\langle e_\gs,e_{\gs'}\rangle_q:=q^{l(\gs)} 
\gd_{\gs\gs'}$ and let $\Ltwo C^i(X)$ be its completion.  Using 
the usual coboundary operator $\gd$, we get  the 
weighted $L^2$-cohomology  spaces, $\Ltwo H^*(X)$.  Let $\partial^q_i$ 
denote the adjoint of $\gd:\Ltwo C^{i-1}(X)\to \Ltwo C^i(X)$.  The  
$\partial^q_i$ give us a chain complex and allow us to define the
weighted $L^2$-homology,  $\Ltwo H_*(X)$. 

The infinite sum
\[
\gG(t):= \sum_{\gamma\in \gG}t^{|\gamma|}.
\]
converges for $t$ in a some neighborhood of $0$ in $[0,\infty)$.  
$\gG(t)$ is the \emph{growth function} of $(\gG,|\ |)$.
It is a power series if $|\ |$ is integer-valued (e.g., if it is given by a 
word length).  Let $\car$ be the region of convergence of $\gG(t)$.  
Suppose $X$ is connected.
The argument in the proof of Proposition~\ref{p:HO} shows that any $0$-cocycle 
is constant and that if $q\in \car$, the only constant which is  square summable 
is $0$.  Hence, $\Ltwo H^0(X)\cong \bR$ if $q\in \car$ and is $0$ if 
$q\notin \car$.

$\gG$ acts on these vector spaces; however, it does not act 
via isometries.

The usual boundary map $\partial$ gives us a different chain complex 
structure (on the same underlying Hilbert spaces $\Ltwo C^i(X)$).

As in Lemma~\ref{l:qtoq-1}, we have the isometry $\theta:\Ltwo C_i(X)\to 
L^2_{1/q}C_i(X)$ 
defined by $e_\gs\to q^{d(\gs)}e_\gs$ intertwines $\partial^q$ with 
$\partial$.  Hence, it induces an isomorphism  $\theta _*:H_*(\Ltwo 
C_*(X),\partial)\to L^2_{1/q} H_*(X)$.

As in Remark~\ref{r:qtoq-1}, we have natural inclusions of cochain 
complexes:
\[
C^i_c(X)\hookrightarrow \Ltwo C^i(X)\hookrightarrow C^i(X).
\]
There is also a version for chain complexes (using ordinary boundary 
map, $\partial$):
\[
C_i(X)\hookrightarrow \Ltwo C_i(X)\hookrightarrow C_i^{lf}(X),
\]
where $C_i^{lf}(X)$ denotes the infinite cellular chains on $X$.  Using 
the isometry $\theta$, we get a monomorphism
of chain complexes
\begin{equation}\label{e:chain4}
C_i(X)\hookrightarrow L^2_{1/q}C_i(X)\mapright{\theta} \Ltwo C_i(X),
\end{equation}
where the boundary maps in the first two terms are the usual ones and 
where the boundary map in the third term is $\partial^q$.  We then have 
the following version of  Theorem~\ref{t:ordinary}.
\pagebreak
\begin{conjecture}\ 
\begin{enumeratei}
\item
For $q\in \car$, the canonical map $\Ltwo H ^i(X)\to H^i(X;\bR)$ is a 
monomorphism.  Moreover, the map $H_i(X;\bR)\to \Ltwo H_i(X)$, induced 
by (\ref{e:chain4}), is a monomorphism with dense image.
\item
For $q^{-1}\in \car$, the canonical map $H^i_c(X;\bR)\to \Ltwo H^i(X)$ is a 
monomorphism with dense image.
\end{enumeratei}
\end{conjecture}

Quite possibly it will be necessary to add  more hypotheses for this 
conjecture to be true.  For example, we might need to assume that the 
$\gG$-action is cocompact and that the norm is given by a word length with 
respect to a set of generators  induced by the choice of 
fundamental domain $D$.

The missing feature from this picture is that for a general group $\gG$ 
there 
is no analog of the Hecke algebra  and no analog of the Hecke - von 
Neumann algebra $\cn$.  So, in the general situation we don't know how to 
define  ``dimension'' and we don't have weighted $L^2$-Betti numbers.
Nevertheless, in some situations 
it is still possible to assign a ``dimension'' to these weighted 
$L^2$-cohomology spaces and obtain weighted $L^2$-Betti numbers.  The 
condition that is needed for these numbers to be well-defined is that the 
$\gG$-action on $X$ has a strict fundamental domain.


\begin{thebibliography}{88}
\bibitem{andreev1}
\textbf{E.M. Andreev}, \emph{On convex polyhedra in
Loba\u{c}evski\u{i} space},  Math. USSR SB. {10} (1970), 413--440.

\bibitem{andreev2}
\textbf{E.M. Andreev}, \emph{On convex polyhedra of finite volume in
Loba\u{c}evski\u{i} space},  Math. USSR SB. {12} (1970), 255--259.

\bibitem{boros}
\textbf{D. Boros}, \emph{$\ell^2$-homology of low dimensional buildings},
Ph.D. thesis, The Ohio State University, 2003.

\bibitem{bueler}
\textbf{E. Bueler}, \emph{The heat kernel weighted Hodge Laplacian on noncompact 
manifolds}, Trans. Amer. Math. Soc. {351} (1999), 683--713.

\bibitem{bourbaki}
\textbf{N. Bourbaki}, \emph{Lie Groups  and Lie Algebras, Chapters 4--6},
Springer-Verlag,  2002.

\bibitem{bh}
\textbf{M. Bridson and A. Haefliger}, \emph{Metric Spaces of Non-positive
Curvature}, Springer-Verlag, Berlin, Heidelberg and New York, 1999.

\bibitem{bronsted}
\textbf{A. Br{\o}nsted},
\emph{An Introduction to Convex Polytopes}, GTM {90},
Springer-Verlag, Berlin, Heidelberg and New York, 1983.

\bibitem{brown}
\textbf{K.S. Brown}, \emph{Buildings}, Springer-Verlag, Berlin, Heidelberg and New 
York,
1989.

\bibitem{cd}
\textbf{R. Charney and M.W. Davis}, \emph{Reciprocity of growth functions of Coxeter
groups}, Geom. Dedicata {39} (1991), 373--378.

\bibitem{cg85}
\textbf{J. Cheeger and M. Gromov}, \emph{Bounds on the von Neumann dimension of
$L^2$-cohomology and the Gauss--Bonnet theorem for open manifolds}, J.
Differential Geom. {21} (1985), 1--34.

\bibitem{cg86}
\textbf{J. Cheeger and M. Gromov}, \emph{$L_2$-cohomology and group cohomology},
Topology {25} (1986), 189--215.

\bibitem{d83}
\textbf{M.W. Davis}, \emph{Groups generated by reflections and aspherical manifolds
not covered by Euclidean space}, Ann. of Math. {117} (1983),
293--325.

\bibitem{d87}
\textbf{M.W. Davis}, \emph{The homology of a space on which a reflection group
acts}, Duke Math. J. {55} (1987), 97--104.

\bibitem{d98}
\textbf{M.W. Davis},  \emph{The cohomology of a Coxeter group with group ring
coefficients}, Duke Math. J. {91} (1998), 297--314.

\bibitem{d98a}
\textbf{M.W. Davis}, \emph{Buildings are $CAT(0)$} in \emph{Geometry and Cohomology in
Group Theory} (eds., P. Kropholler, G. Niblo and R. St\"{o}hr), LMS
Lecture Note Series {252}, pp. 108--123,  Cambridge Univ.
Press, Cambridge, 1998.

\bibitem{d02}
\textbf{M.W. Davis},  \emph{Nonpositive curvature and reflection groups} in
\emph{Handbook
of Geometric Topology}  (eds., R.J. Daverman and R.B. Sher), pp. 373-422,
Elsevier, Amsterdam, 2002.

\bibitem{davisbook}
\textbf{M.W. Davis}, \emph{The Geometry and Topology of Coxeter Groups}, London Math. Soc. Monograph Series, Vol. 32, Princeton Univ. Press, Princeton, in press.

\bibitem{ddjo06}
\textbf{M.W. Davis, J. Dymara, T. Januszkiewicz and B. Okun}, \emph{Cohomology of Coxeter groups with group ring coefficients: II}, Algebraic \& Geometric Topology {6} (2006), 1289--1318.

\bibitem{dl}
\textbf{M.W. Davis and I.J. Leary}, \emph{The $\ell^2$-cohomology of Artin 
groups},  J.
of London  Math. Soc. {68} (2003), 493--510.

\bibitem{dm}
\textbf{M.W. Davis and J. Meier}, \emph{The topology at infinity of Coxeter groups
and buildings}, Comment. Math. Helv. {77} (2002), 746--766.

\bibitem{do}
\textbf{M.W. Davis and B. Okun}, \emph{Vanishing theorems and conjectures for the
$\ell^2$-homology of right-angled Coxeter groups}, Geometry \& Topology
{5} (2001), 7--74.

\bibitem{do04}
\textbf{M.W. Davis and B. Okun},  \emph{$\ell^2$-homology of right-angled Coxeter groups
associated  to barycentric subdivisions}, Topology and Its Applications
{140} (2004), 197--202.

\bibitem{dix}
\textbf{J. Dixmier}, \emph{Les $C^{\ast}$-alg\'ebres et leur rep\'esentations},
Gauthier-Villars, Paris, 1964.

\bibitem{dix81}
\textbf{J. Dixmier}, \emph{Von Neumann Algebras}, North Holland, 1981.

\bibitem{dod77}
\textbf{J. Dodziuk}, \emph{de Rham--Hodge theory for $L^2$-cohomology of infinite
coverings}, Topology {16} (1977), 157--165.

\bibitem{dymara}
\textbf{J. Dymara}, \emph{$L^2$-cohomology of buildings with fundamental class}, Proc. Amer. Math. Soc. {132} (2004), 1839--1843.

\bibitem{dym}
\textbf{J. Dymara},   \emph{Thin buildings}, Geometry \& Topology {10} (2006), 667--694.

\bibitem{dj}
\textbf{J. Dymara and T. Januszkiewicz}, \emph{Cohomology of buildings and their
automorphism groups}, Invent. Math. {150} (2002), 579--627.

\bibitem{eckmann}
\textbf{B. Eckmann}, \emph{Introduction to $\ell_2$-methods in topology: reduced
$\ell_2$-homology, harmonic chains, $\ell_2$-Betti numbers}, Israel J.
Math. {117} (2000), 183--219.

\bibitem{gal}
\textbf{{\'S}. Gal}, \emph{Real root conjecture fails for five-
  and higher-dimensional spheres}, Discrete Comput. Geom. 34 (2005) 269--284

\bibitem{cgon}
\textbf{C. Gonciulea}, \emph{Virtual epimorphisms of Coxeter groups onto free 
groups}, Ph.D. thesis, The Ohio State University, 2000.


\bibitem{kl}
\textbf{D. Kazhdan and G. Lusztig}, \emph{Representations of Coxeter groups and
Hecke algebras}, Invent. Math. {53} (1979), 165--184.

\bibitem{luck}
\textbf{W.  L\"{u}ck}, \emph{$L^2$-Invariants:  Theory and Applications to Geometry 
and $K$-theory},
Springer-Verlag, Berlin, Heidelberg and New York, 2002.

\bibitem{l83}
\textbf{G. Lusztig}, \emph{Left cells in Weyl groups}, in \emph{Lie Group
Representations I} (eds., R.L.R. Herb and J. Rosenberg),  Springer
Lecture Notes in Math. {1024}, pp. 99--111, Springer-Verlag, Berlin,
Heidelberg and New York, 1983.

\bibitem{mv}
\textbf{G. Margulis and E.B. Vinberg}, \emph{Some linear groups virtually having a 
free quotient}, J. Lie Theory {10} (2000), 171--180.

\bibitem{moussong}
\textbf{G. Moussong},  \emph{Hyperbolic Coxeter groups}, Ph.D. thesis, The Ohio 
State University, 1988.


\bibitem{ronan}
\textbf{M. Ronan}, \emph{Lectures on Buildings}, Academic Press, San Diego, 1989.

\bibitem{serre}
\textbf{J-P. Serre}, \emph{Cohomologie des groupes discrets} in \emph{Prospects in
Mathematics}, Annals of Math. Studies {70}, pp. 77--169, Princeton
Univ. Press, Princeton, 1971.

\bibitem{s}
\textbf{L. Solomon}, \emph{A decomposition of the group algebra of a finite Coxeter
group}, J. of Algebra {9} (1968), 220--239.

\bibitem{steinberg}
\textbf{R. Steinberg}, \emph{Endomorphisms of linear algebraic groups}, Mem. Amer.
Math. Soc. {80} (1968).

\bibitem{tits}
\textbf{J. Tits}, \emph{Le probl\'em des mots dans les groupes de Coxeter} in
\emph{Symposia Mathematica (INDAM, Rome, 1967/68)}, vol. 1, pp. 175-185,
Academic Press, London, 1969.

\bibitem{weibel}
\textbf{C. Weibel}, \emph{An Introduction to Homological Algebra}, Cambridge 
Studies in Advanced Mathematics {38}, Cambridge University Press, 
Cambridge, 1994.

\end{thebibliography}

\end{document}